\documentclass[a4paper,10pt]{amsart}
\usepackage[margin=1.25in]{geometry}
\usepackage{amssymb,bm,amsfonts,amsmath, psfrag,multicol,amsthm}
\usepackage{tensor}
\usepackage[dvipsnames]{xcolor}
\usepackage[nosort]{cite}
\usepackage{cite}
\usepackage{ytableau}
\usepackage{tikz}
\usetikzlibrary{decorations.pathreplacing,calligraphy,positioning, shapes.geometric}
\usetikzlibrary{arrows,patterns.meta}
\usetikzlibrary{decorations.pathreplacing}
\usepackage{graphicx}
\usepackage{ifpdf}
\usepackage{dashrule}
\usepackage{booktabs}
\usepackage{pdfpages}
\usepackage{epsfig}
\usepackage{epstopdf}
\usepackage{mathrsfs}
\usepackage{stmaryrd}
\usepackage{enumerate}
\usepackage{soul}
\usepackage[colorlinks,urlcolor=RedViolet,linkcolor=red,anchorcolor=green,citecolor=blue]{hyperref}
\usepackage{longtable}
\usepackage{booktabs}
\usepackage{threeparttable}
\usepackage{array}
\usepackage{multirow}
\usepackage{boldline}

\allowdisplaybreaks
\numberwithin{equation}{section}
\newtheorem{thm}{Theorem}[section]

\newtheorem{definition}{Definition}[section]
\newtheorem{prop}{Proposition}[section]

\newtheorem{exam}{Example}[section]

\newtheorem{example}{Example}[section]
\newtheorem{remark}{Remark}[section]

\def\P{{\mathcal P}}
\def\D{{\mathcal D}}
\def\G{{\mathcal G}}
\def\B{{\mathfrak B}}
\def\qed{\hfill \rule{4pt}{7pt}}
\def\pf{\noindent {\it{Proof.} \hskip 2pt}}
\def\ij{\langle i,j \rangle}
\def\royalblue{\textcolor{royalblue}}
\def\red{\textcolor{red}}
\def\des{\operatorname{des}}
\def\Des{\operatorname{Des}}
\def\ides{\operatorname{ides}}
\def\iDes{\operatorname{iDes}}
\def\asc{\operatorname{asc}}
\def\iasc{\operatorname{iasc}}

\def\ob{{\overline{B}}}
\def\ub{{\underline{B}}}

\definecolor{babyblue}{rgb}{0.54, 0.81,1}
\definecolor{airforceblue}{rgb}{0.36, 0.54, 0.66}
\definecolor{gainsboro}{rgb}{0.86, 0.86, 0.86}
\definecolor{afblue}{rgb}{0.36, 0.54, 0.66}
\definecolor{purple}{rgb}{0.62,0.13,0.94}
\definecolor{lightsalmon}{rgb}{0.92,0.57,0.45}
\definecolor{royalblue}{rgb}{0.25,0.41,0.87}
\definecolor{bondiblue}{rgb}{0.0, 0.58, 0.71}
\definecolor{orange}{rgb}{1.0, 0.6, 0.4}
\def\blue{\textcolor{blue}}
\def\red{\textcolor{red}}
\def\orange{\textcolor{orange}}
\def\gainsboro{\textcolor{gainsboro}}
\def\afblue{\textcolor{afblue}}
\def\purple{\textcolor{purple}}
\def\LightSalmon{\textcolor{LightSalmon}}
\def\royalblue{\textcolor{royalblue}}
\def\babyblue{\textcolor{babyblue}}
\def\bondiblue{\textcolor{bondiblue}}

\newcommand{\cpfthm}[1]{\noindent{\emph{Combinatorial Proof of Theorem #1.}\hskip 2pt}}
\newcommand{\cpfcor}[1]{\noindent{\emph{Combinatorial Proof of Corollary #1.}\hskip 2pt}}
\newcommand{\cpflem}[1]{\noindent{\emph{Combinatorial Proof of Lemma #1.}\hskip 2pt}}

\begin{document}

\title[The descents and inverse descents on the hyperoctahedral group]{On the Combinatorics of descents and inverse descents \\in the hyperoctahedral group}

\author[X. Gao]{Xiaoqin Gao}

\author[F.Z.K. Li]{Frank Z.K. Li}

\author[L. Wan]{Lingli Wan}

\author[J.Y.X Yang]{Jane Y.X. Yang}

\address{School of Science, Chongqing University of Posts and Telecommunications, Chongqing 400065, People's Republic of China}
\email{xqgao7@163.com, zkli@cqupt.edu.cn, llwan0@163.com, yangyingxue@cqupt.edu.cn}

\date{\today}

\keywords{hyperoctahedral group, inverse descents, signed permutation grids,
signed involutions, recurrence formulas}

\subjclass[2010]{05A05, 05A19}

\maketitle

\begin{abstract}
The elements in the hyperoctahedral group $\B_n$ can be treated as signed permutations with the natural order $\cdots<-2<-1<0<1<2<\cdots$,
or as colored permutations
with the $r$-order $-1<_r-2<_r\cdots<_r0<_r1<_r2<_r\cdots$.
For any $\pi\in\B_n$,
let $\des^B(\pi)$ and $\ides^B(\pi)$ be the number of descents and inverse descents in $\pi$ under the natural order,
and let $\des_B(\pi)$ and $\ides_B(\pi)$ be the number of descents and inverse descents in $\pi$ under the $r$-order.
In this paper, by investigating signed permutation grids under both the natural order and the $r$-order,
we give combinatorial proofs for six recurrence formulas of the joint distribution
of descents and inverse descents over
the hyperoctahedral group $\B_n$,
the set in involutions of $\B_n$ denoted by $\mathcal{I}_n^B$,
and the set of fixed-point free involutions in $\B_n$ denoted by $\mathcal{J}_n^B$,
respectively.
Some of these six formulas are new,
and some reveal the combinatorial essences
of the results obtained
by Visontai, Moustakas and Cao-Liu through algebraic approaches such as quasisymmetric functions.
Furthermore, from these formulas,
we conclude that $(\des^B,\ides^B)$ and $(\des_B,\ides_B)$ are equidistributed
over both $\B_n$ and $\mathcal{I}_n^B$,
but not on $\mathcal{J}_n^B$.

\end{abstract}


%
%
%
%
%
%


\section{Introduction}

The main purpose of this paper is to present combinatorial proofs
for six recurrence formulas of descents and inverse descents over the
hyperoctahedral groups.
For any positive integer $n$,
let $[n]=\{1, 2, \ldots, n\}$,
and $\Omega_n=\{0,1, -1, 2, -2, \ldots, n, -n\}$.
Following the notations in \cite[Section 2]{AdinAthan2017},
we denote by $\mathfrak{B}_n$ the set of all signed permutations of length $n$,
which are bijective maps $\pi: \Omega_n\rightarrow \Omega_n$ with $\pi(a)=b$ implying $\pi(-a)=-b$ for any $a\in\Omega_n$.
By this restriction,
we use $\pi=\pi_1\pi_2\cdots \pi_n$ with
$\pi_i=\pi(i)$ for $i\in[n]$
to indicate the signed permutations such that $\overline{a}$ stands for $-a\in\Omega_n$ with $a>0$.
For example, we write $\overline{3}16\overline{5}24\in \B_6$.

The set $\B_n$ is known as the \emph{hyperoctahedral group},
a \emph{Coxeter group} with type $B$ (see \cite[Section 8.1]{Bjorner2005}),
which can be regarded as
the group of symmetries of the $n$-dimensional cube (see \cite[Part III]{Petersen-2015}).
In addition, the hyperoctahedral group $\B_n$ can be also
viewed as the \emph{wreath product} $\mathbb{Z}_r\wr\mathfrak{S}_n$ with $r=2$, where $\mathfrak{S}_n$ is the symmetric group,
and $\mathbb{Z}_r$  is the cyclic group of order $r$.
The wreath product $\mathbb{Z}_r\wr\mathfrak{S}_n$ can be realized as
the $r$-colored permutation group $\mathfrak{S}_{n,r}$,
which allows us to treat signed permutations as $2$-colored permutations.
We refer the reader to \cite[Section 2.1.7]{Athan2018} and \cite[Section 1.2]{Moustakas-2021}
to check more information on wreath product and colored permutations.

Since the hyperoctahedral group can be realized by two different permutations,
such as signed permutations or colored permutations,
it is necessary to choose a suitable total order on $\mathbb{Z}$ when studying statistics
over permutations.
We use the natural order
\begin{equation}\label{norder}
\cdots <-2<-1<0<1<2<\cdots
\end{equation}
on $\mathbb{Z}$ when considering the hyperoctahedral group $\B_n$ as a Coxeter group of type $B$.
On the other hand, by taking $\B_n$ as a colored permutation group,
it is convenient to use the $r$-order
\begin{equation}\label{rorder}
-1<_r -2 <_r \cdots <_r 0 <_r 1 <_r 2 <_r\cdots
\end{equation}
on $\mathbb{Z}$,
where $r$ stands for the lexicographic order on the set $\mathbb{Z}_{>0}\{-,+\}$ (see \cite[Equation (2.1)]{AdinAthan2017}).
For any signed permutation $\pi\in\B_n$, we define
\[
\Des^B(\pi)=\{i\in[0,n-1]\mid \pi_i >\pi_{i+1}\},
\]
and
\[
\Des_B(\pi)=\{i\in[0,n-1]\mid \pi_i >_r \pi_{i+1}\},
\]
where $[a,b]$ denotes the set $\{a, a+1,\ldots,b\}$ for integers $a\leq b$,
and note that $\pi(0)=0$.
We call the indices $i\in\Des^B(\pi)$ the \emph{$\overline{B}$-descent} in $\pi$,
and the indices $i\in\Des_B$ the \emph{$\underline{B}$-descent}.
Without causing confusion,
we refer to $\overline{B}$-descents and $\underline{B}$-descents collectively as \emph{descents} in $\pi$.
Let
\[
\des^B(\pi)=|\Des^B(\pi)|, \qquad \des_B(\pi)=|\Des_B(\pi)|,
\]
then the statistics $\des^B$ and $\des_B$ record
the number of descents in Coxeter groups (see \cite[Section 13.1]{Petersen-2015}) and colored permutations (see \cite[Section 2]{Steing1994}), respectively.
Actually, in \cite[Page 218]{AdinBrenti2001}, Adin et al. remark that
$\des^B$ and $\des_B$ are equally distributed over $\B_n$,
and it is easy to check this fact via bijections.
For any positive integer $n$, the $\B_n$-Eulerian polynomial $B_n(t)$ is defined as
\begin{equation}\label{BnEulerian}
B_n(t)=\sum_{\pi\in\B_n}t^{\des^B(\pi)}=\sum_{\pi\in\B_n}t^{\des_B(\pi)}.
\end{equation}
It is proved that  $\B_n$-Eulerian polynomials share same properties with
classical Eulerian polynomials, such as $\gamma$-positivity,
where combinatorial interpretations for the $\gamma$-coefficients in $B_n(t)$
can be found in \cite[Theorem 2.10]{Athan2018}, \cite[Theorem 4.7]{Chow2008} and \cite[Proposition 4.15]{Petersen-2007}.

In this paper, we mainly focus on the joint distribution of descents and inverse descents on the  hyperoctahedral group $\B_n$,
and several certain subsets of $\B_n$.
Particularly,
we give combinatorial interpretations for several recursive formulas of
the number of signed permutation of length $n$ with given numbers of descents and inverse descents.
For any $\pi\in\B_n$, denote by $\pi^{-1}$ the inverse signed permutation of $\pi$.
We call indices $i$ in the set
\[
\iDes^B(\pi)=\{i\in[0,n-1]\mid \pi^{-1}_i>\pi^{-1}_{i+1}\}
\]
\emph{inverse $\overline{B}$-descents}, and indices $i$ in the set
\[
\iDes_B(\pi)=\{i\in[0,n-1]\mid \pi^{-1}_i>_r\pi^{-1}_{i+1}\}
\]
\emph{inverse $\underline{B}$-descents}.
The inverse $\overline{B}$-descents and the inverse $\underline{B}$-descents are both regarded as
\emph{inverse descents} in the signed permutation $\pi$.
Moreover, they are also simply referred as \emph{idescents} (see \cite{Rawlings-1984}).
Let
\[
\ides^B(\pi)=|\Des^B(\pi)|, \qquad \ides_B(\pi)=|\Des_B(\pi)|.
\]

For $0\leq i,j\leq n$, let
\[
\overline{\mathcal{B}}_{n,i,j}=\{\pi\in\B_n\mid \mbox{ $\des^B(\pi)=i$ and $\ides^B(\pi)=j$}\}
\]
and $\overline{b}_{n,i,j}=|\overline{\mathcal{B}}_{n,i,j}|$.
The generating function of $\overline{b}_{n,i,j}$ is called as
the type $B$ two-side Eulerian polynomial (see \cite[Section 4.1]{Visontai-2013}):
\begin{equation}\label{genB(s,t)}
  \overline{B}_n(s,t)=\sum_{\pi\in\B_n}s^{\des^B(\pi)}t^{\ides^B(\pi)}=
  \sum_{0\leq i,j\leq n}\overline{b}_{n,i,j}s^it^j.
\end{equation}
We may think of $\overline{B}_n(s,t)$ as the generating function for
the joint distribution of $\des^B$ and $\ides^B$.
The classical two-side Eulerian polynomials $A_n(s,t)$ are defined on the symmetric group $\mathfrak{S}_n$, where one can see \cite[Page 167]{Petersen-2013} for details.
Thus the polynomial $\overline{B}_n(s,t)$ is a generalization of $A_n(s,t)$ from the Coxeter group of type $A$ to type $B$.
Starting with an identity of binomial coefficients given by Petersen to  prove \cite[Equation (9)]{Petersen-2013},
Visontai \cite{Visontai-2013} derived a differential equation about
$\overline{B}_n(s,t)$.

\begin{thm}[\cite{Visontai-2013}, Theorem 15]
For $n\geq 2$,
\begin{equation}\label{pdeofBn(st)}
\begin{aligned}
n \overline{B}_{n}(s, t)=&\left(2 n^{2} s t-n s t+n\right) \overline{B}_{n-1}(s, t) \\
&+(2 n s t(1-s)+s(1-s)(1-t)) \frac{\partial}{\partial s} \overline{B}_{n-1}(s, t) \\
&+(2 n s t(1-t)+t(1-s)(1-t)) \frac{\partial}{\partial t} \overline{B}_{n-1}(s, t) \\
&+2 s t(1-s)(1-t) \frac{\partial^{2}}{\partial s \partial t} \overline{B}_{n-1}(s, t)
\end{aligned}
\end{equation}
with initial value $\overline{B}_1(s,t)=1+st$.
\end{thm}

Comparing the coefficients of $\overline{B}_n(s,t)$ in the two sides of \eqref{pdeofBn(st)}, we obtain a recurrence of $\overline{b}_{n,i,j}$.
\begin{thm}\label{thmof^obnij}
For $n\geq 2$ and $0\leq i,j\leq n$,
we have
\begin{equation}\label{rec^obnij}
\begin{aligned}
n \overline{b}_{n, i, j}
=&(n+i+j+2ij) \overline{b}_{n-1, i, j}\\
&+((2n-1)i-(2i+1)(j-1)) \overline{b}_{n-1, i, j-1}\\
&+((2n-1)j-(2j+1)(i-1)) \overline{b}_{n-1, i-1, j}\\
&+((2n^2-n)+2(i-1)(j-1)+(1-2n)(i+j-2)) \overline{b}_{n-1, i-1, j-1},
\end{aligned}
\end{equation}
where $\overline{b}_{1,0,0}=\overline{b}_{1,1,1}=1$, $\overline{b}_{1,0,1}=\overline{b}_{1,1,0}=0$ and $\overline{b}_{n,i,j}=0$ if $i< 0$ or $j<0$.
\end{thm}

Motivated by the work in \cite{Li-Liu-Yang},
by generalizing the concept of permutation grids \cite[Section 1.5]{Stanley-2012} to signed permutations,
we give a combinatorial interpretation of \eqref{rec^obnij}.
Moreover, we also investigate the joint distribution of the statistics $\des_B$ and $\ides_B$. 
For $0\leq i,j\leq n$, let
\[
\underline{\mathcal{B}}_{n,i,j}=\{\pi\in\B_n\mid \mbox{ $\des_B(\pi)=i$ and $\ides_B(\pi)=j$}\},
\]
and $\underline{b}_{n,i,j}=|\underline{\mathcal{B}}_{n,i,j}|$.

\begin{thm}\label{thmof_bnij}
For $n\geq 2$ and $0\leq i,j\leq n$,
we have
\begin{equation}\label{rec_bnij}
\begin{aligned}
n \underline{b}_{n, i, j}
=&(n+i+j+2ij) \underline{b}_{n-1, i, j}\\
&+((2n-1)i-(2i+1)(j-1)) \underline{b}_{n-1, i, j-1}\\
&+((2n-1)j-(2j+1)(i-1)) \underline{b}_{n-1, i-1, j}\\
&+((2n^2-n)+2(i-1)(j-1)+(1-2n)(i+j-2)) \underline{b}_{n-1, i-1, j-1},
\end{aligned}
\end{equation}
where $\underline{b}_{1,0,0}=\underline{b}_{1,1,1}=1$, $\underline{b}_{1,0,1}=\underline{b}_{1,1,0}=0$ and $\underline{b}_{n,i,j}=0$ if $i<0$ or $j< 0$.
\end{thm}

Hence  by Theorems \ref{thmof^obnij} and \ref{thmof_bnij},
we find that $\underline{b}_{n,i,j}$ shares exactly the same recursive relation and initial values with $\overline{b}_{n,i,j}$ through combinatorial approaches,
which proves that $(\des^\ob, \ides^\ob)$ and $(\des^\ub,\ides^\ub)$ are equally distributed over the set $\mathfrak{B}_n$.The type $B$ two-side Eulerian polynomial \eqref{genB(s,t)} can be unified to
the two-side $\B_n$-Eulerian polynomial:
\[
B_n(s,t)=\sum_{\pi\in\B_n}s^{\des^B(\pi)}t^{\ides^B(\pi)}=\sum_{\pi\in\B_n}s^{\des_B(\pi)}t^{\ides_B(\pi)},
\]
which is a generalization of $\B_n$-Eulerian polynomial \eqref{BnEulerian}.
And it would be of interest to find a bijection between
$\overline{b}_{n, i, j}$ and $\underline{b}_{n, i, j}$.

Given any positive integer $n$, for a signed permutation $\pi\in\mathfrak{B}_n$,
we call $\pi$ a \emph{signed involution}, or \emph{involution} for short,
if $\pi^{-1}=\pi$.
Let $\mathcal{I}_{n}^B$ be the set of all involutions in $\mathfrak{B}_{n}$,
and $\mathcal{J}_{n}^B$ be the set of all fixed-point free signed involutions in $\mathfrak{B}_{n}$,
where a \emph{fixed-point free involution} $\pi\in\B_n$ satisfies
$\pi(i)\neq i$ and $\pi(i)\neq -i$ for each $i\in[n]$.
Let
\begin{equation*}
\mathcal{I}_{n,k}^\ob=\{\pi\in\mathcal{I}_n^B\mid \des^B(\pi)=k\},\quad
\mathcal{I}_{n,k}^\ub=\{\pi\in\mathcal{I}_n^B\mid \des_B(\pi)=k\},
\end{equation*}
\begin{equation*}
\mathcal{J}_{n,k}^\ob=\{\pi\in\mathcal{J}_n^B\mid \des^B(\pi)=k\},\quad
\mathcal{J}_{n,k}^\ub=\{\pi\in\mathcal{J}_n^B\mid \des_B(\pi)=k\},
\end{equation*}
and
\[
I_{n,k}^\ob=|\mathcal{I}_{n,k}^\ob|,\quad
I_{n,k}^\ub=|\mathcal{I}_{n,k}^\ub|,\quad
J_{n,k}^\ob=|\mathcal{J}_{n,k}^\ob|,\quad
J_{n,k}^\ub=|\mathcal{J}_{n,k}^\ub|.
\]
Since the idescents coincide with descents in signed involutions,
with the experience in proving Theorems \ref{thmof^obnij} and \ref{thmof_bnij},
we establish a bijection on permutation grids of signed involutions,
that proves the following linear recurrence formula of $I^\ob_{n,k}$.

\begin{thm}\label{thmIobnk}
For $n \geq 3$ and $k \geq 0$, we have
\begin{equation}\label{recIobnk}
\begin{aligned}
  n I_{n, k}^\ob=&(2k+1) I_{n-1, k}^\ob+(2n-2k+1) I_{n-1, k-1}^\ob+\left(n-1+2k(k+1)\right) I_{n-2, k}^\ob\\
  &+(2(n-1)+4(n-k-1)(k-1)) I_{n-2, k-1}^\ob\\
  &+\left((2n-3)(n-1)+2(k-2)(k-2n+1)\right) I_{n-2, k-2}^\ob,
\end{aligned}
\end{equation}
where $I_{1,0}^\ob=1, I_{1,1}^\ob=1, I_{2,0}^\ob=1, I_{2,1}^\ob=4, I_{2,2}^\ob=1$ and $I_{n, k}^\ob=0$ for $k<0$.
\end{thm}

Notice that in terms of a specialization of Poirier's signed quasisymmetric functions \cite[Section 3]{Poirier-1998},
Moustakas \cite{Moustakas-2019}
computed the generating function of $I_n^\ub(t)$,
where
\begin{equation}\label{genInubt}
  I_n^\ub(t)=\sum_{\pi\in\mathcal{I}^B_n}t^{\des_B(\pi)}
  =\sum_{k=0}^n I_{n,k}^\ub t^k.
\end{equation}

\begin{thm}[\cite{Moustakas-2019}, Theorem 1.1]
  We have
  \begin{equation}\label{genIub}
    \sum_{n=0}^{\infty}I_n^\ub(t)\frac{x^n}{(1-t)^{n+1}}
    =\sum_{m=0}^{\infty}\frac{t^m}{(1-x)^{2m+1}(1-x^2)^{m^2}},
  \end{equation}
  where $I_0^\ub(t):=1$.
\end{thm}

Based on the generating function \eqref{genIub},
Moustakas \cite[Proposition 5]{Moustakas-2019} algebraically derived a recurrence of
$I_{n,k}^\ub$, which is used to prove the symmetry and unimodality of
the sequence $I_{n,0}^\ub,I_{n,1}^\ub,\ldots,I_{n,n}^\ub$.
Similar to the proof of Theorem \ref{thmIobnk},
we give a combinatorial interpretation of this recurrence.

\begin{thm}[\cite{Moustakas-2019}, Proposition 5]\label{thmIubnk}
For $n \geq 3$ and $k \geq 0$, we have
\begin{equation}\label{recIubnk}
\begin{aligned}
  n I_{n, k}^\ub=&(2k+1) I_{n-1, k}^\ub+(2n-2k+1) I_{n-1, k-1}^\ub+\left(n-1+2k(k+1)\right) I_{n-2, k}^\ub\\
  &+(2(n-1)+4(n-k-1)(k-1)) I_{n-2, k-1}^\ub\\
  &+\left((2n-3)(n-1)+2(k-2)(k-2n+1)\right) I_{n-2, k-2}^\ub,
\end{aligned}
\end{equation}
where $I_{1,0}^\ub=1,I_{1,1}^\ub=1, I_{2,0}^\ub=1, I_{2,1}^\ub=4,I_{2,2}^\ub=1$ and $I_{n, k}^\ub=0$ for $k<0$.
\end{thm}

As mentioned both in \cite[Section 2.1.5]{Athan2018} and \cite[Remark 3]{Moustakas-2019},
the generating function $I_n^\ub(t)$ defined in \eqref{genInubt}
is presumably but not obviously unaffected when $\des_B$ is replaced by $\des^B$, which has been confirmed for $n\leq 5$.
Indeed, combining Theorems \ref{thmIobnk} and \ref{thmIubnk} verifies
the above fact for any positive integer $n$,
which implies that the $\ob$-descents and $\ub$-descents are
equally distributed on $\mathcal{I}^B_n$.
Thus, $I_n^\ub(t)$ can be actually unified as
\begin{equation*}
  I_n^B(t)=\sum_{\pi\in\mathcal{I}^B_n}t^{\des_B(\pi)}
  =\sum_{\pi\in\mathcal{I}^B_n}t^{\des^B(\pi)}.
\end{equation*}


Another part of our work is to deduce the recursive relations of  $J_{2 n, k}^\ob$ and  $J_{2 n, k}^\ub$ respectively via pure combinatorial approaches,
as sated by Theorems \ref{thmJobnk}  and \ref{thmJubnk} later.

\begin{thm}\label{thmJobnk}
For $n \geq 2$ and $k \geq 0$, we have
\begin{equation}\label{recJobnk}
\begin{aligned}
n J_{2 n, k}^\ob=& {(k^2+k+n-1) J_{2 n-2, k}^\ob+2((k-1)(2 n-k-1)+n) J_{2 n-2, k-1}^\ob } \\[3pt]
&+((2 n-k)(2 n-k+1)+(n-1)) J_{2 n-2, k-2}^\ob,
\end{aligned}
\end{equation}
where $J_{2,0}^\ob=0, J_{2,1}^\ob=2, J_{2,2}^\ob=0$ and $J_{2 n, k}^\ob=0$ if $k<0$.
\end{thm}

For any positive integer $n$, let $J^\ub_n(t)$ be the generating function of $J^\ub_{n,k}$:
\[
J_n^\ub(t)=\sum_{\pi\in\mathcal{J}^B_n}t^{\des_B(\pi)}
  =\sum_{k=0}^n J_{n,k}^\ub t^k.
\]
Motivated by the work of Moustakas \cite{Moustakas-2019} on the sequence $I^\ub_{n,k}$, Cao and Liu \cite[Theorem 3.1]{CaoLiu-2021} compute the generating function of $J^\ub_n(t)$ also by signed quasisymmetric functions.

\begin{thm}[\cite{CaoLiu-2021}, Theorem 3.1]
  We have
  \begin{equation}\label{genJub}
    \sum_{n=0}^{\infty}J_n^\ub(t)\frac{x^n}{(1-t)^{n+1}}
    =\sum_{m=0}^{\infty}\frac{x^m}{(1-x^2)^{m^2}},
  \end{equation}
  where $J_0^\ub(t):=1$.
\end{thm}

By comparing the coefficients in the generating function \eqref{genJub}, they \cite{CaoLiu-2021} obtained the following recurrence relation of the sequence $J^\ub_{n,k}$.

\begin{thm}[\cite{CaoLiu-2021}, Theorem 3.2]\label{thmJubnk}
For $n \geq 2$ and $k \geq 0$, we have
\begin{equation}\label{recJubnk}
\begin{aligned}
n J^\ub_{2 n, k}=& {(k^2+n-1) J^\ub_{2 n-2, k}+(2(k-1)(2n-k)+1) J^\ub_{2 n-2, k-1} } \\[3pt]
&+((k-2)(k-4n)+4n^2-3n) J^\ub_{2 n-2, k-2},
\end{aligned}
\end{equation}
where $J^\ub_{2,0}=0, J^\ub_{2,1}=1, J^\ub_{2,2}=1$ and $J^\ub_{2 n, k}=0$ if $k<0$.
\end{thm}

Using the recursive formula \eqref{recJubnk}, Cao and Liu \cite{CaoLiu-2021} verified the symmetric, unimodal and $\gamma$-positive properties of the sequence $J_{n,0}^\ub,J_{n,1}^\ub,\ldots,J_{n,n}^\ub$.
By modifying the construction in the proof of Theorem \ref{thmJobnk},
we give Theorem \ref{thmJubnk} a combinatorial proof.
Note that other than the identically distributed behavior of
the $\ub$-descents and the $\ob$-descents over $\mathcal{I}^B_n$,
the distributions of these two statistics are not the same over the set  $\mathcal{J}^B_n$,
which can be directly seen from comparing \eqref{recJobnk} to \eqref{recJubnk}.

\begin{remark}
Since the symmetric group $\mathfrak{S}_n$
is a subgroup of the hyperoctahedral group $\B_n$,
we generalize  the constructive methods in \cite{Li-Liu-Yang},
which give combinatorial proofs for recurrence formula related to the descents and inverse descents on $\mathfrak{S}_n$ (see \cite[Eq. (10)]{Petersen-2013}),
and recurrences of descents in involutions and fixed-point free involutions in $\mathfrak{S}_n$ respectively (see \cite[Theorems 2.2 and 2.1]{Guo-Zeng-2006}).

\end{remark}

The rest of this paper is organized as follows.
Section \ref{Sec-grids} is devoted to introduce  signed permutation grids and the related key operations in the construction of subsequent  proofs.
In Section \ref{Sec-GeoDes}, we investigate geometric properties of $\overline{B}$-descents and inverse $\overline{B}$-descents under natural order in signed permutation grids,
especially derive enumerative formulas for certain internal structures of grids.
After this, we combinatorially prove Theorems \ref{thmof^obnij} and \ref{thmIobnk}, \ref{thmJobnk} in Sections \ref{Sec-recBnij} and \ref{Sec-involup}, respectively.
By following the techniques of handling objects in natural order,
in Section \ref{Sec-geounder},
we study the properties of $\underline{B}$-descents and inverse $\underline{B}$-descents in grids under the $r$-order,
which is served as the preparation for combinatorial proofs of Theorem \ref{thmof_bnij}, \ref{thmIubnk} and
\ref{thmJubnk} given in Section \ref{Sec-recunderall}.

\section{Signed permutation grids and inserting operations}\label{Sec-grids}

In this section, we introduce the definition of signed permutation grids and
relevant notations that are frequently used in our proofs.

\subsection{Signed permutation grids}
The geometric representations of permutations \cite[Section 1.5]{Stanley-2012} such as matrices, grids, trees
provide us with another intuitive perspective to study the properties of permutations.
We first extend the definition of permutation grids to signed permutation grids.

\begin{definition}\label{defgrids}
  For any signed permutation $\pi=\pi_1\pi_2\ldots \pi_n\in\B_n$,
  the \emph{signed permutation grid} $P_\pi$ of $\pi$ is defined as
  an $n\times n$ grid with the $|\pi_i|$-th (from the left to the right) square in the $i$-th (from the top to the bottom) row
  filled in if $\pi_i>0$, and filled by slashes if $\pi_i<0$.
\end{definition}

For accuracy,
we refer to $\ij\in[n]\times [n]$ and $(i,j)\in[n+1]\times [n+1]$ as
the square located in the $i$-th row and the $j$-th column,
and the grid point intersected by the $i$-th horizontal line and the $j$-th vertical line, respectively.
For $1\leq i \leq n$, we define the sign of the square
by calling $\langle i, -\pi_i\rangle$ as \emph{negative square},
and $\langle i, \pi_i\rangle$ as \emph{positive square}.
For example, the signed permutation grid of $\pi=\overline{3}16\overline{5}24$ is given by
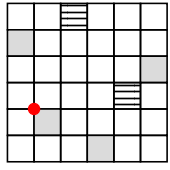
\begin{center}
\begin{tikzpicture}[scale = 0.7]
\def\hezi{-- +(5mm,0mm) -- +(5mm,5mm) -- +(0mm,5mm) -- cycle [line width=0.6pt]}
\def\judy{-- +(5mm,0mm) -- +(5mm,5mm) -- +(0mm,5mm) -- cycle [line width=0.6pt,fill=gainsboro]}
\def\nbx{-- +(5mm,0mm) -- +(5mm,5mm) -- +(0mm,5mm) -- cycle [line width=0.6pt,pattern={Lines[angle=45,distance=2.5pt, line width=0.6pt]}, pattern color=black]}
\tikzstyle{rdot}=[circle,fill=red,draw=red,inner sep=1.5]
\draw (0mm,0mm)\hezi;
\draw (5mm,0mm)\hezi;
\draw (10mm,0mm)\nbx;
\draw (15mm,0mm)\hezi;
\draw (20mm,0mm)\hezi;
\draw (25mm,0mm)\hezi;
\draw (0mm,-5mm)\judy;
\draw (5mm,-5mm)\hezi;
\draw (10mm,-5mm)\hezi;
\draw (15mm,-5mm)\hezi;
\draw (20mm,-5mm)\hezi;
\draw (25mm,-5mm)\hezi;
\draw (0mm,-10mm)\hezi;
\draw (5mm,-10mm)\hezi;
\draw (10mm,-10mm)\hezi;
\draw (15mm,-10mm)\hezi;
\draw (20mm,-10mm)\hezi;
\draw (25mm,-10mm)\judy;

\draw (0mm,-15mm)\hezi;
\draw (5mm,-15mm)\hezi;
\draw (10mm,-15mm)\hezi;
\draw (15mm,-15mm)\hezi;
\draw (20mm,-15mm)\nbx;
\draw (25mm,-15mm)\hezi;
\draw (0mm,-20mm)\hezi;
\draw (5mm,-20mm)\judy;
\draw (10mm,-20mm)\hezi;
\draw (15mm,-20mm)\hezi;
\draw (20mm,-20mm)\hezi;
\draw (25mm,-20mm)\hezi;
\draw (0mm,-25mm)\hezi;
\draw (5mm,-25mm)\hezi;
\draw (10mm,-25mm)\hezi;
\draw (15mm,-25mm)\judy;
\draw (20mm,-25mm)\hezi;
\draw (25mm,-25mm)\hezi;
\node[rdot] at (5mm, -15mm){};
\end{tikzpicture}
\end{center}
with a negative square $\langle4,5\rangle$ and
a positive square $\langle3,6\rangle$,
and the grid point $(5,2)$ addressed in red.

For any singed permutation,
transposing its grid yields
the grid of its inverse.
Since a signed involution $\pi$ satisfies  $\pi=\pi^{-1}$,
the grid $P_\pi$ is symmetrical about the main diagonal.
Moreover, the grid of a fixed-point free signed involution $\sigma$
cannot have filled squares on its main diagonal since $\sigma_i\neq \pm i$.
The  grids of involution $\pi=\overline{4}25\overline{1}3$ and fixed-point free involution $\sigma=\overline{5}326\overline{1}4$ are listed below.
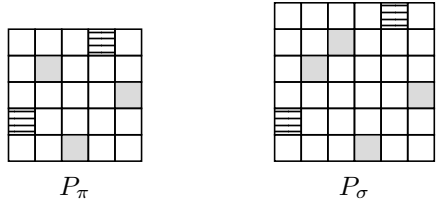
\begin{center}
\begin{tikzpicture}[scale = 0.7]
\def\hezi{-- +(5mm,0mm) -- +(5mm,5mm) -- +(0mm,5mm) -- cycle [line width=0.6pt]}
\def\judy{-- +(5mm,0mm) -- +(5mm,5mm) -- +(0mm,5mm) -- cycle [line width=0.6pt,fill=gainsboro]}
\def\nbx{-- +(5mm,0mm) -- +(5mm,5mm) -- +(0mm,5mm) -- cycle [line width=0.6pt,pattern={Lines[angle=45,distance=2.5pt, line width=0.6pt]}, pattern color=black]}
\tikzstyle{rdot}=[circle,fill=red,draw=red,inner sep=1.5]

\draw (0mm,-5mm)\hezi;
\draw (5mm,-5mm)\hezi;
\draw (10mm,-5mm)\hezi;
\draw (15mm,-5mm)\nbx;
\draw (20mm,-5mm)\hezi;
\draw (0mm,-10mm)\hezi;
\draw (5mm,-10mm)\judy;
\draw (10mm,-10mm)\hezi;
\draw (15mm,-10mm)\hezi;
\draw (20mm,-10mm)\hezi;
\draw (0mm,-15mm)\hezi;
\draw (5mm,-15mm)\hezi;
\draw (10mm,-15mm)\hezi;
\draw (15mm,-15mm)\hezi;
\draw (20mm,-15mm)\judy;
\draw (0mm,-20mm)\nbx;
\draw (5mm,-20mm)\hezi;
\draw (10mm,-20mm)\hezi;
\draw (15mm,-20mm)\hezi;
\draw (20mm,-20mm)\hezi;
\draw (0mm,-25mm)\hezi;
\draw (5mm,-25mm)\hezi;
\draw (10mm,-25mm)\judy;
\draw (15mm,-25mm)\hezi;
\draw (20mm,-25mm)\hezi;
\node at (12.5mm,-30mm){$P_{\pi}$};


\begin{scope}[shift={(50mm,0mm)}]
\draw (0mm,0mm)\hezi;
\draw (5mm,0mm)\hezi;
\draw (10mm,0mm)\hezi;
\draw (15mm,0mm)\hezi;
\draw (20mm,0mm)\nbx;
\draw (25mm,0mm)\hezi;
\draw (0mm,-5mm)\hezi;
\draw (5mm,-5mm)\hezi;
\draw (10mm,-5mm)\judy;
\draw (15mm,-5mm)\hezi;
\draw (20mm,-5mm)\hezi;
\draw (25mm,-5mm)\hezi;
\draw (0mm,-10mm)\hezi;
\draw (5mm,-10mm)\judy;
\draw (10mm,-10mm)\hezi;
\draw (15mm,-10mm)\hezi;
\draw (20mm,-10mm)\hezi;
\draw (25mm,-10mm)\hezi;

\draw (0mm,-15mm)\hezi;
\draw (5mm,-15mm)\hezi;
\draw (10mm,-15mm)\hezi;
\draw (15mm,-15mm)\hezi;
\draw (20mm,-15mm)\hezi;
\draw (25mm,-15mm)\judy;
\draw (0mm,-20mm)\nbx;
\draw (5mm,-20mm)\hezi;
\draw (10mm,-20mm)\hezi;
\draw (15mm,-20mm)\hezi;
\draw (20mm,-20mm)\hezi;
\draw (25mm,-20mm)\hezi;
\draw (0mm,-25mm)\hezi;
\draw (5mm,-25mm)\hezi;
\draw (10mm,-25mm)\hezi;
\draw (15mm,-25mm)\judy;
\draw (20mm,-25mm)\hezi;
\draw (25mm,-25mm)\hezi;
\node at (15mm,-30mm){$P_{\sigma}$};
\end{scope}
\end{tikzpicture}
\end{center}
By noticing the restrictions of symmetric form and no filled squares on the main diagonal,
we have $\mathcal{J}_{2k+1}^B=\emptyset$ with any nonnegative integer $k$.

\subsection{Inserting and deleting operations}
To establish the connection between $\B_n$ and $\B_{n+1}$,
we introduce the following inserting and deleting operations on signed permutation grids.

Let $\pi\in\mathfrak{B}_n$ and $i,j\in[n+1]$,
then the \emph{inserting operations} $\varphi_{(i,j)}$ and $\overline{\varphi}_{(i,j)}$ are defined by
respectively inserting a positive square and a negative square at the grid point $(i,j)$
of the grid $P_\pi$,
while preserving the relative positions of the original filled squares.
%
%
For example, if $\pi=\overline{3}16\overline{5}24$, then the grids of $\varphi_{(3,5)}(\pi)=\overline{3}157\overline{6}24$ and $\overline{\varphi}_{(4,3)}(\pi)=\overline{4}17\overline{3}\overline{6}25$
are given as follows.
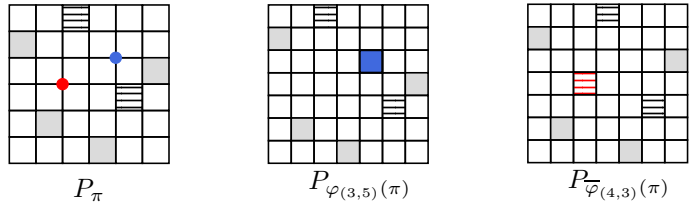
\begin{center}
\begin{tikzpicture}[scale = 0.7]
\def\hezi{-- +(5mm,0mm) -- +(5mm,5mm) -- +(0mm,5mm) -- cycle [line width=0.6pt]}
\def\judy{-- +(5mm,0mm) -- +(5mm,5mm) -- +(0mm,5mm) -- cycle [line width=0.6pt,fill=gainsboro]}
\def\nbx{-- +(5mm,0mm) -- +(5mm,5mm) -- +(0mm,5mm) -- cycle [line width=0.6pt,pattern={Lines[angle=45,distance=2.5pt, line width=0.6pt]}, pattern color=black]}
\tikzstyle{rdot}=[circle,fill=red,draw=red,inner sep=1.5]
\tikzstyle{gdot}=[circle,fill=royalblue,draw=royalblue,inner sep=1.5]

\draw (0mm,0mm)\hezi;
\draw (5mm,0mm)\hezi;
\draw (10mm,0mm)\nbx;
\draw (15mm,0mm)\hezi;
\draw (20mm,0mm)\hezi;
\draw (25mm,0mm)\hezi;
\draw (0mm,-5mm)\judy;
\draw (5mm,-5mm)\hezi;
\draw (10mm,-5mm)\hezi;
\draw (15mm,-5mm)\hezi;
\draw (20mm,-5mm)\hezi;
\draw (25mm,-5mm)\hezi;
\draw (0mm,-10mm)\hezi;
\draw (5mm,-10mm)\hezi;
\draw (10mm,-10mm)\hezi;
\draw (15mm,-10mm)\hezi;
\draw (20mm,-10mm)\hezi;
\draw (25mm,-10mm)\judy;
\draw (0mm,-15mm)\hezi;
\draw (5mm,-15mm)\hezi;
\draw (10mm,-15mm)\hezi;
\draw (15mm,-15mm)\hezi;
\draw (20mm,-15mm)\nbx;
\draw (25mm,-15mm)\hezi;
\draw (0mm,-20mm)\hezi;
\draw (5mm,-20mm)\judy;
\draw (10mm,-20mm)\hezi;
\draw (15mm,-20mm)\hezi;
\draw (20mm,-20mm)\hezi;
\draw (25mm,-20mm)\hezi;
\draw (0mm,-25mm)\hezi;
\draw (5mm,-25mm)\hezi;
\draw (10mm,-25mm)\hezi;
\draw (15mm,-25mm)\judy;
\draw (20mm,-25mm)\hezi;
\draw (25mm,-25mm)\hezi;
\node[gdot] at (20mm, -5mm){};
\node[rdot] at (10mm, -10mm){};
\node at (15mm,-30mm){$P_{\pi}$};
\end{tikzpicture}
\hspace{3em}
\begin{tikzpicture}[scale = 0.6]
\def\hezi{-- +(5mm,0mm) -- +(5mm,5mm) -- +(0mm,5mm) -- cycle [line width=0.6pt]}
\def\judy{-- +(5mm,0mm) -- +(5mm,5mm) -- +(0mm,5mm) -- cycle [line width=0.6pt,fill=gainsboro]}
\def\nbx{-- +(5mm,0mm) -- +(5mm,5mm) -- +(0mm,5mm) -- cycle [line width=0.6pt,pattern={Lines[angle=45,distance=2.5pt, line width=0.6pt]}, pattern color=black]}
\tikzstyle{rdot}=[circle,fill=purple,draw=purple,inner sep=1.5]
\tikzstyle{gdot}=[circle,fill=royalblue,draw=royalblue,inner sep=1.5]

\draw (0mm,0mm)\hezi;
\draw (5mm,0mm)\hezi;
\draw (10mm,0mm)\nbx;
\draw (15mm,0mm)\hezi;
\draw (20mm,0mm)\hezi;
\draw (25mm,0mm)\hezi;
\draw (30mm,0mm)\hezi;
\draw (0mm,-5mm)\judy;
\draw (5mm,-5mm)\hezi;
\draw (10mm,-5mm)\hezi;
\draw (15mm,-5mm)\hezi;
\draw (20mm,-5mm)\hezi;
\draw (25mm,-5mm)\hezi;
\draw (30mm,-5mm)\hezi;
\draw (0mm,-10mm)\hezi;
\draw (5mm,-10mm)\hezi;
\draw (10mm,-10mm)\hezi;
\draw (15mm,-10mm)\hezi;
\draw (20mm,-10mm)-- +(5mm,0mm) -- +(5mm,5mm) -- +(0mm,5mm) -- cycle [line width=0.8pt,fill=royalblue];
\draw (25mm,-10mm)\hezi;
\draw (30mm,-10mm)\hezi;
\draw (0mm,-15mm)\hezi;
\draw (5mm,-15mm)\hezi;
\draw (10mm,-15mm)\hezi;
\draw (15mm,-15mm)\hezi;
\draw (20mm,-15mm)\hezi;
\draw (25mm,-15mm)\hezi;
\draw (30mm,-15mm)\judy;
\draw (0mm,-20mm)\hezi;
\draw (5mm,-20mm)\hezi;
\draw (10mm,-20mm)\hezi;
\draw (15mm,-20mm)\hezi;
\draw (20mm,-20mm)\hezi;
\draw (25mm,-20mm)\nbx;
\draw (30mm,-20mm)\hezi;
\draw (0mm,-25mm)\hezi;
\draw (5mm,-25mm)\judy;
\draw (10mm,-25mm)\hezi;
\draw (15mm,-25mm)\hezi;
\draw (20mm,-25mm)\hezi;
\draw (25mm,-25mm)\hezi;
\draw (30mm,-25mm)\hezi;
\draw (0mm,-30mm)\hezi;
\draw (5mm,-30mm)\hezi;
\draw (10mm,-30mm)\hezi;
\draw (15mm,-30mm)\judy;
\draw (20mm,-30mm)\hezi;
\draw (25mm,-30mm)\hezi;
\draw (30mm,-30mm)\hezi;

\node at (20mm,-35mm){$P_{\varphi_{(3,5)}(\pi)}$};
\end{tikzpicture}
\hspace{3em}
\begin{tikzpicture}[scale = 0.6]
\def\hezi{-- +(5mm,0mm) -- +(5mm,5mm) -- +(0mm,5mm) -- cycle [line width=0.6pt]}
\def\judy{-- +(5mm,0mm) -- +(5mm,5mm) -- +(0mm,5mm) -- cycle [line width=0.6pt,fill=gainsboro]}
\def\nbx{-- +(5mm,0mm) -- +(5mm,5mm) -- +(0mm,5mm) -- cycle [line width=0.6pt,pattern={Lines[angle=45,distance=2.5pt, line width=0.6pt]}, pattern color=black]}
\tikzstyle{rdot}=[circle,fill=purple,draw=purple,inner sep=1.5]
\tikzstyle{gdot}=[circle,fill=royalblue,draw=royalblue,inner sep=1.5]

\draw (0mm,0mm)\hezi;
\draw (5mm,0mm)\hezi;
\draw (10mm,0mm)\hezi;
\draw (15mm,0mm)\nbx;
\draw (20mm,0mm)\hezi;
\draw (25mm,0mm)\hezi;
\draw (30mm,0mm)\hezi;
\draw (0mm,-5mm)\judy;
\draw (5mm,-5mm)\hezi;
\draw (10mm,-5mm)\hezi;
\draw (15mm,-5mm)\hezi;
\draw (20mm,-5mm)\hezi;
\draw (25mm,-5mm)\hezi;
\draw (30mm,-5mm)\hezi;
\draw (0mm,-10mm)\hezi;
\draw (5mm,-10mm)\hezi;
\draw (10mm,-10mm)\hezi;
\draw (15mm,-10mm)\hezi;
\draw (20mm,-10mm)\hezi;
\draw (25mm,-10mm)\hezi;
\draw (30mm,-10mm)\judy;
\draw (0mm,-15mm)\hezi;
\draw (5mm,-15mm)\hezi;
\draw (10mm,-15mm)-- +(5mm,0mm) -- +(5mm,5mm) -- +(0mm,5mm) -- cycle [line width=0.6pt,pattern={Lines[angle=45,distance=2.5pt, line width=0.6pt]}, pattern color=red];
\draw (15mm,-15mm)\hezi;
\draw (20mm,-15mm)\hezi;
\draw (25mm,-15mm)\hezi;
\draw (30mm,-15mm)\hezi;
\draw (0mm,-20mm)\hezi;
\draw (5mm,-20mm)\hezi;
\draw (10mm,-20mm)\hezi;
\draw (15mm,-20mm)\hezi;
\draw (20mm,-20mm)\hezi;
\draw (25mm,-20mm)\nbx;
\draw (30mm,-20mm)\hezi;
\draw (0mm,-25mm)\hezi;
\draw (5mm,-25mm)\judy;
\draw (10mm,-25mm)\hezi;
\draw (15mm,-25mm)\hezi;
\draw (20mm,-25mm)\hezi;
\draw (25mm,-25mm)\hezi;
\draw (30mm,-25mm)\hezi;
\draw (0mm,-30mm)\hezi;
\draw (5mm,-30mm)\hezi;
\draw (10mm,-30mm)\hezi;
\draw (15mm,-30mm)\hezi;
\draw (20mm,-30mm)\judy;
\draw (25mm,-30mm)\hezi;
\draw (30mm,-30mm)\hezi;

\node at (20mm,-35mm){$P_{\overline{\varphi}_{(4,3)}(\pi)}$};
\end{tikzpicture}
\end{center}

Note that
inserting squares at some specific grid points in the grid of an involution can still maintain the symmetry of the grid.

\begin{prop}\label{propi+1i}
Let $\pi\in\mathcal{I}_n^B$, then for $1\leq i\leq n$,
we have $\varphi_{(i+1,i)}(\pi)\in \mathcal{I}_{n+1}^B$ if $\pi_i=i$,
and $\overline{\varphi}_{(i+1,i)}(\pi)\in \mathcal{I}_{n+1}^B$ if $\pi_i=\overline{i}$.
\end{prop}

On the other hand,
for any signed permutation $\sigma\in \mathfrak{B}_{n+1}$ with a filled square $\ij$ in its grid $P_\sigma$,
by deleting the $i$-th row and the $j$-th column of squares from
the grid,
we obtain the grid $P_\pi$ for a signed permutation $\pi\in \mathfrak{B}_n$.
This process can be described as
deleting the positive square $\ij$ if $\sigma_i=j$,
or the negative square $\ij$ if $\sigma_i=\overline{j}$ from the grid $P_\sigma$,
and denoted by the \emph{deleting operations} $\varphi^{-1}_{\ij}$ or $\overline{\varphi}^{-1}_{\ij}$, respectively.
These two deleting operations can be roughly viewed as the inverse of
the inserting operations $\varphi_{(i,j)}$ and $\overline{\varphi}_{(i,j)}$,
since $\pi=\varphi^{-1}_{\ij}(\sigma)$ (resp. $\pi=\overline{\varphi}^{-1}_{\ij}(\sigma)$)
if and only if $\sigma=\varphi_{(i,j)}(\pi)$ (resp.$\sigma=\overline{\varphi}^{-1}_{\ij}(\pi)$).

\subsection{Double inserting and deleting operations on signed involutions}
Based on the inserting operations $\varphi_{(i,j)}$ and $\overline{\varphi}_{(i,j)}$,
we can define new operations to insert two filled squares
into grids of signed involutions at the same time.

For $\pi\in\mathcal{I}_n^B$ and $1\leq i,j\leq n+1$ with $i\neq j$,
the \emph{double inserting operations} $\xi_{(i,j)}$ and $\overline{\xi}_{(i,j)}$ are defined by
\begin{equation}\label{xi(ij)}
\xi_{(i,j)}=\left\{
        \begin{array}{ll}
          \varphi_{(j+1,i)}\circ \varphi_{(i,j)}, & i<j, \\[3pt]
          \varphi_{(j,i+1)}\circ \varphi_{(i,j)}, & i>j,
        \end{array}
      \right.
\end{equation}
\begin{equation}\label{bxi(ij)}
\overline{\xi}_{(i,j)}=\left\{
        \begin{array}{ll}
          \overline{\varphi}_{(j+1,i)}\circ \overline{\varphi}_{(i,j)}, & i<j, \\[3pt]
          \overline{\varphi}_{(j,i+1)}\circ \overline{\varphi}_{(i,j)}, & i>j,
        \end{array}
      \right.
\end{equation}
where the composition $f \circ g$ obeys the rule $f\circ g(x)=f(g(x))$ for the indeterminate $x$.

From a geometric view,
the grid of $\xi_{(i,j)}(\pi)$ (resp. $\overline{\xi}_{(i,j)}(\pi)$) is obtained by
inserting two positive (reps. negative) squares at the grid points $(i,j)$ and $(j,i)$ of the grid $P_\pi$ at the same time.
The following figures show the grids of signed involution $\overline{1}32$ and
$\xi_{(3,2)}(132)=\overline{1}4523$,
$\overline{\xi}_{(1,3)}(132)=\overline{4}\overline{2}5\overline{1}3$.
\begin{center}
\begin{tikzpicture}[scale = 1]
\def\hezi{-- +(5mm,0mm) -- +(5mm,5mm) -- +(0mm,5mm) -- cycle [line width=0.6pt]}
\def\judy{-- +(5mm,0mm) -- +(5mm,5mm) -- +(0mm,5mm) -- cycle [line width=0.6pt,fill=gainsboro]}
\def\nbx{-- +(5mm,0mm) -- +(5mm,5mm) -- +(0mm,5mm) -- cycle [line width=0.6pt,pattern={Lines[angle=45,distance=2.5pt, line width=0.6pt]}, pattern color=black]}

\tikzstyle{rdot}=[circle,fill=red,draw=red,inner sep=2]
\tikzstyle{gdot}=[circle,fill=royalblue,draw=royalblue,inner sep=2]
\draw (0mm,-15mm)\nbx;
\draw (5mm,-15mm)\hezi;
\draw (10mm,-15mm)\hezi;
\draw (0mm,-20mm)\hezi;
\draw (5mm,-20mm)\hezi;
\draw (10mm,-20mm)\judy;
\draw (0mm,-25mm)\hezi;
\draw (5mm,-25mm)\judy;
\draw (10mm,-25mm)\hezi;
\node at (7.5mm,-30mm){$P_{132}$};
\node[rdot] at (10mm,-10mm){};
\node[gdot] at (5mm,-20mm){};

\begin{scope}[shift={(30mm,-7mm)},scale=.8]
\draw (0mm,-5mm)\nbx;
\draw (5mm,-5mm)\hezi;
\draw (10mm,-5mm)\hezi;
\draw (15mm,-5mm)\hezi;
\draw (20mm,-5mm)\hezi;
\draw (0mm,-10mm)\hezi;
\draw (5mm,-10mm)\hezi;
\draw (10mm,-10mm)\hezi;
\draw (15mm,-10mm)-- +(5mm,0mm) -- +(5mm,5mm) -- +(0mm,5mm) -- cycle [line width=0.6pt,fill=royalblue];
\draw (20mm,-10mm)\hezi;
\draw (0mm,-15mm)\hezi;
\draw (5mm,-15mm)\hezi;
\draw (10mm,-15mm)\hezi;
\draw (15mm,-15mm)\hezi;
\draw (20mm,-15mm)\judy;
\draw (0mm,-20mm)\hezi;
\draw (5mm,-20mm)-- +(5mm,0mm) -- +(5mm,5mm) -- +(0mm,5mm) -- cycle [line width=0.6pt,fill=royalblue];
\draw (10mm,-20mm)\hezi;
\draw (15mm,-20mm)\hezi;
\draw (20mm,-20mm)\hezi;
\draw (0mm,-25mm)\hezi;
\draw (5mm,-25mm)\hezi;
\draw (10mm,-25mm)\judy;
\draw (15mm,-25mm)\hezi;
\draw (20mm,-25mm)\hezi;
\node at (12.5mm,-30mm){$P_{\overline{1}4523}$};
\end{scope}


\begin{scope}[shift={(65mm,-7mm)},scale=.8]
\draw (0mm,-5mm)\hezi;
\draw (5mm,-5mm)\hezi;
\draw (10mm,-5mm)\hezi;
\draw (15mm,-5mm)-- +(5mm,0mm) -- +(5mm,5mm) -- +(0mm,5mm) -- cycle [line width=0.6pt,pattern={Lines[angle=45,distance=2.5pt, line width=0.6pt]}, pattern color=red];
\draw (20mm,-5mm)\hezi;
\draw (0mm,-10mm)\hezi;
\draw (5mm,-10mm)\nbx;
\draw (10mm,-10mm)\hezi;
\draw (15mm,-10mm)\hezi;
\draw (20mm,-10mm)\hezi;
\draw (0mm,-15mm)\hezi;
\draw (5mm,-15mm)\hezi;
\draw (10mm,-15mm)\hezi;
\draw (15mm,-15mm)\hezi;
\draw (20mm,-15mm)\judy;
\draw (0mm,-20mm)-- +(5mm,0mm) -- +(5mm,5mm) -- +(0mm,5mm) -- cycle [line width=0.6pt,pattern={Lines[angle=45,distance=2.5pt, line width=0.6pt]}, pattern color=red];
\draw (5mm,-20mm)\hezi;
\draw (10mm,-20mm)\hezi;
\draw (15mm,-20mm)\hezi;
\draw (20mm,-20mm)\hezi;
\draw (0mm,-25mm)\hezi;
\draw (5mm,-25mm)\hezi;
\draw (10mm,-25mm)\judy;
\draw (15mm,-25mm)\hezi;
\draw (20mm,-25mm)\hezi;
\node at (12.5mm,-30mm){$P_{\overline{4}\overline{2}5\overline{1}3}$};
\end{scope}
\end{tikzpicture}
\end{center}
Thus we see $\xi_{(i,j)}(\pi), \overline{\xi}_{(i,j)}(\pi)\in\mathcal{I}_{n+2}^B$,
and $\xi_{(i,j)}(\pi)=\xi_{(j,i)}(\pi)$,
$\overline{\xi}_{(i,j)}(\pi)=\overline{\xi}_{(j,i)}(\pi)$.

For $\pi\in\mathcal{I}_n^B$ and $1<i<n+1$, the \emph{double inserting operations} $\eta_i$ and $\overline{\eta}_i$ are defined as
\begin{equation}\label{eta_i}
  \eta_i=\varphi_{(i,i)}^2 \quad\mbox{and} \quad \overline{\eta}_i=\overline{\varphi}_{(i,i)}^2
\end{equation}
and the \emph{double inserting operation} $\eta'_i$ and $\overline{\eta'}_i$ are defined as
\begin{equation}\label{eta'_i}
  \eta'_i=\varphi_{(i+1,i)}\circ \varphi_{(i,i)} \quad\mbox{and} \quad \overline{\eta'}_i=\overline{\varphi}_{(i+1,i)}\circ \overline{\varphi}_{(i,i)}.
\end{equation}

In order to see the insertion behavior of these operations on the grid, we name some filled square patterns on the grid as follows.
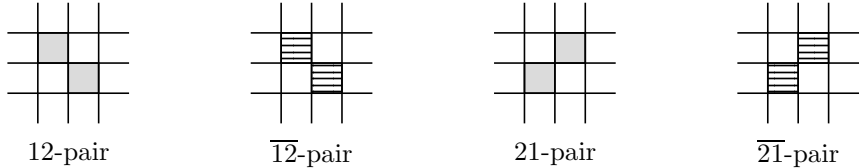
\begin{center}
\begin{tikzpicture}[scale = 0.8]
\def\hezi{-- +(5mm,0mm) -- +(5mm,5mm) -- +(0mm,5mm) -- cycle [line width=0.6pt]}
\def\judy{-- +(5mm,0mm) -- +(5mm,5mm) -- +(0mm,5mm) -- cycle [line width=0.6pt,fill=gainsboro]}
\def\nbx{-- +(5mm,0mm) -- +(5mm,5mm) -- +(0mm,5mm) -- cycle [line width=0.6pt,pattern={Lines[angle=45,distance=2.5pt, line width=0.6pt]}, pattern color=black]}
\tikzstyle{bdot}=[circle,fill=black,draw=black,inner sep=1.2]
\tikzstyle{gdot}=[circle,fill=royalblue,draw=royalblue,inner sep=1.5]
\draw (5mm,0mm)\judy;
\draw (10mm,-5mm)\judy;
\draw[line width=0.6pt] (0mm,5mm)--(20mm,5mm);
\draw[line width=0.6pt] (0mm,0mm)--(20mm,0mm);
\draw[line width=0.6pt] (0mm,-5mm)--(20mm,-5mm);
\draw[line width=0.6pt] (5mm,-10mm)--(5mm,10mm);
\draw[line width=0.6pt] (15mm,-10mm)--(15mm,10mm);
\draw[line width=0.6pt] (10mm,-10mm)--(10mm,10mm);
\node at (10mm,-15mm){$12$-pair};

\begin{scope}[shift={(40mm,0mm)}]
\draw (5mm,0mm)\nbx;
\draw (10mm,-5mm)\nbx;
\draw[line width=0.6pt] (0mm,5mm)--(20mm,5mm);
\draw[line width=0.6pt] (0mm,0mm)--(20mm,0mm);
\draw[line width=0.6pt] (0mm,-5mm)--(20mm,-5mm);
\draw[line width=0.6pt] (5mm,-10mm)--(5mm,10mm);
\draw[line width=0.6pt] (15mm,-10mm)--(15mm,10mm);
\draw[line width=0.6pt] (10mm,-10mm)--(10mm,10mm);
\node at (10mm,-15mm){$\overline{12}$-pair};
\end{scope}

\begin{scope}[shift={(80mm,0mm)}]
\draw (5mm,-5mm)\judy;
\draw (10mm,0mm)\judy;
\draw[line width=0.6pt] (0mm,5mm)--(20mm,5mm);
\draw[line width=0.6pt] (0mm,0mm)--(20mm,0mm);
\draw[line width=0.6pt] (0mm,-5mm)--(20mm,-5mm);
\draw[line width=0.6pt] (5mm,-10mm)--(5mm,10mm);
\draw[line width=0.6pt] (15mm,-10mm)--(15mm,10mm);
\draw[line width=0.6pt] (10mm,-10mm)--(10mm,10mm);
\node at (10mm,-15mm){$21$-pair};
\end{scope}

\begin{scope}[shift={(120mm,0mm)}]
\draw (5mm,-5mm)\nbx;
\draw (10mm,0mm)\nbx;
\draw[line width=0.6pt] (0mm,5mm)--(20mm,5mm);
\draw[line width=0.6pt] (0mm,0mm)--(20mm,0mm);
\draw[line width=0.6pt] (0mm,-5mm)--(20mm,-5mm);
\draw[line width=0.6pt] (5mm,-10mm)--(5mm,10mm);
\draw[line width=0.6pt] (15mm,-10mm)--(15mm,10mm);
\draw[line width=0.6pt] (10mm,-10mm)--(10mm,10mm);
\node at (10mm,-15mm){$\overline{21}$-pair};
\end{scope}
\end{tikzpicture}
\end{center}
Then the inserting operations $\eta_i$, $\overline{\eta}_i$
extend the grid point $(i,i)$ in $P_\pi$
to a $12$-pair, $\overline{1}\overline{2}$-pair, respectively,
and the inserting operations $\eta'_i$ and $\overline{\eta'}_i$
extend the grid point $(i,i)$ in $P_\pi$
to a $21$-pair and $\overline{2}\overline{1}$-pair, respectively.
The following example gives the signed permutation grids obtained
by employing $\overline{\eta}_2$ and $\eta'_4$ on the grid of $21\overline{3}$.
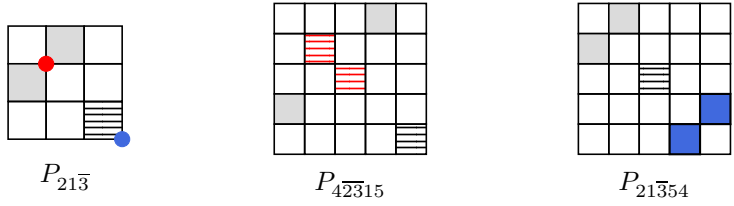
\begin{center}
\begin{tikzpicture}[scale = 1]
\def\hezi{-- +(5mm,0mm) -- +(5mm,5mm) -- +(0mm,5mm) -- cycle [line width=0.6pt]}
\def\judy{-- +(5mm,0mm) -- +(5mm,5mm) -- +(0mm,5mm) -- cycle [line width=0.6pt,fill=gainsboro]}
\def\nbx{-- +(5mm,0mm) -- +(5mm,5mm) -- +(0mm,5mm) -- cycle [line width=0.6pt,pattern={Lines[angle=45,distance=2.5pt, line width=0.6pt]}, pattern color=black]}
\tikzstyle{rdot}=[circle,fill=red,draw=red,inner sep=2]
\tikzstyle{gdot}=[circle,fill=royalblue,draw=royalblue,inner sep=2]

\draw (0mm,-15mm)\hezi;
\draw (5mm,-15mm)\judy;
\draw (10mm,-15mm)\hezi;
\draw (0mm,-20mm)\judy;
\draw (5mm,-20mm)\hezi;
\draw (10mm,-20mm)\hezi;
\draw (0mm,-25mm)\hezi;
\draw (5mm,-25mm)\hezi;
\draw (10mm,-25mm)\nbx;
\node at (7.5mm,-30mm){$P_{21\overline{3}}$};
\node[rdot] at (5mm,-15mm){};
\node[gdot] at (15mm,-25mm){};

\begin{scope}[shift={(35mm,-7mm)},scale=.8]
\draw (0mm,-5mm)\hezi;
\draw (5mm,-5mm)\hezi;
\draw (10mm,-5mm)\hezi;
\draw (15mm,-5mm)\judy;
\draw (20mm,-5mm)\hezi;
\draw (0mm,-10mm)\hezi;
\draw (5mm,-10mm)-- +(5mm,0mm) -- +(5mm,5mm) -- +(0mm,5mm) -- cycle [line width=0.6pt,pattern={Lines[angle=45,distance=2.5pt, line width=0.6pt]}, pattern color=red];
\draw (10mm,-10mm)\hezi;
\draw (15mm,-10mm)\hezi;
\draw (20mm,-10mm)\hezi;
\draw (0mm,-15mm)\hezi;
\draw (5mm,-15mm)\hezi;
\draw (10mm,-15mm)-- +(5mm,0mm) -- +(5mm,5mm) -- +(0mm,5mm) -- cycle [line width=0.6pt,pattern={Lines[angle=45,distance=2.5pt, line width=0.6pt]}, pattern color=red];
\draw (15mm,-15mm)\hezi;
\draw (20mm,-15mm)\hezi;
\draw (0mm,-20mm)\judy;
\draw (5mm,-20mm)\hezi;
\draw (10mm,-20mm)\hezi;
\draw (15mm,-20mm)\hezi;
\draw (20mm,-20mm)\hezi;
\draw (0mm,-25mm)\hezi;
\draw (5mm,-25mm)\hezi;
\draw (10mm,-25mm)\hezi;
\draw (15mm,-25mm)\hezi;
\draw (20mm,-25mm)\nbx;
\node at (12.5mm,-30mm){$P_{4\overline{2}\overline{3}15}$};
\end{scope}


\begin{scope}[shift={(75mm,-7mm)},scale=.8]
\draw (0mm,-5mm)\hezi;
\draw (5mm,-5mm)\judy;
\draw (10mm,-5mm)\hezi;
\draw (15mm,-5mm)\hezi;
\draw (20mm,-5mm)\hezi;
\draw (0mm,-10mm)\judy;
\draw (5mm,-10mm)\hezi;
\draw (10mm,-10mm)\hezi;
\draw (15mm,-10mm)\hezi;
\draw (20mm,-10mm)\hezi;
\draw (0mm,-15mm)\hezi;
\draw (5mm,-15mm)\hezi;
\draw (10mm,-15mm)\nbx;
\draw (15mm,-15mm)\hezi;
\draw (20mm,-15mm)\hezi;
\draw (0mm,-20mm)\hezi;
\draw (5mm,-20mm)\hezi;
\draw (10mm,-20mm)\hezi;
\draw (15mm,-20mm)\hezi;
\draw (20mm,-20mm)-- +(5mm,0mm) -- +(5mm,5mm) -- +(0mm,5mm) -- cycle [line width=0.8pt,fill=royalblue];
\draw (0mm,-25mm)\hezi;
\draw (5mm,-25mm)\hezi;
\draw (10mm,-25mm)\hezi;
\draw (15mm,-25mm)-- +(5mm,0mm) -- +(5mm,5mm) -- +(0mm,5mm) -- cycle [line width=0.8pt,fill=royalblue];
\draw (20mm,-25mm)\hezi;
\node at (12.5mm,-30mm){$P_{21\overline{3}54}$};

\end{scope}

\end{tikzpicture}
\end{center}
Clearly, we have $\eta_i(\pi),\eta'_i(\pi), \overline{\eta}_i(\pi),\overline{\eta'}_i(\pi)\in\mathcal{I}_{n+2}^B$ for any $\pi\in\mathcal{I}_n^B$ and $1\leq i\leq n+1$.

The corresponding deleting operations of the above double insertion operations can be constructed
by $\varphi^{-1}_{\ij}$ and $\overline{\varphi}^{-1}_{\ij}$.
Given $\sigma\in \mathcal{I}_{n+2}^B$, and $1\leq i,j\leq n+1$ with $i\neq j$,
the \emph{double deleting operation} $\xi^{-1}_{\ij}$ and $\overline{\xi}^{-1}_{\ij}$ are defined as
\begin{equation}\label{invxi(ij)}
\xi^{-1}_{\ij}=\left\{
        \begin{array}{ll}
          \varphi^{-1}_{\langle i,j\rangle}\circ \varphi^{-1}_{\langle j+1,i\rangle}, & i<j\text{ and }\sigma_i=j+1, \\[8pt]
          \varphi^{-1}_{\langle i,j\rangle}\circ \varphi^{-1}_{\langle j,i+1\rangle}, & i>j\text{ and }\sigma_{i+1}=j;
        \end{array}
      \right.
\end{equation}
and
\begin{equation}\label{invbxi(ij)}
\overline{\xi}^{-1}_{\ij}=\left\{
        \begin{array}{ll}
          \overline{\varphi}^{-1}_{\langle i,j\rangle}\circ \overline{\varphi}^{-1}_{\langle j+1,i\rangle}, & i<j\text{ and }\sigma_i=\overline{j+1}, \\[8pt]
          \overline{\varphi}^{-1}_{\langle i,j\rangle}\circ \overline{\varphi}^{-1}_{\langle j,i+1\rangle}, & i>j\text{ and }\sigma_{i+1}=\overline{j}.
        \end{array}
      \right.
\end{equation}
By the symmetry of the grid $P_\sigma$,
applying the operation $\xi^{-1}_{\ij}$ (resp. $\overline{\xi}^{-1}_{\langle i,j\rangle}$) with $i<j$
is to delete the positive (resp. negative) squares $\langle j+1,i\rangle$ and $\langle i,j+1\rangle$ from $P_\sigma$,
and applying the operation $\xi^{-1}_{\ij}$ (resp. $\overline{\xi}^{-1}_{\langle i,j\rangle}$) with $i>j$
is to delete the positive (resp. negative) squares $\langle j,i+1\rangle$ and $\langle i+1,j\rangle$ from $P_\sigma$.
This implies  $\xi^{-1}_{\ij}(\sigma)=\pi$ if $\xi_{(i,j)}(\pi)=\sigma$,
and $\overline{\xi}^{-1}_{\ij}(\sigma)=\pi$ if $\overline{\xi}_{(i,j)}(\pi)=\sigma$.
Note that $\xi^{-1}_{\ij}$ or $\overline{\xi}^{-1}_{\ij}$ cannot delete $21$-pair or $\overline{2}\overline{1}$-pair from the main diagonal of the grid $P_\sigma$.

For $\sigma\in\mathcal{I}_{n+2}^B$ and $1\leq i\leq n+1$, we define
\emph{the double deleting operations} $\eta^{-1}_i$, $\overline{\eta}^{-1}_i$,
$\eta'^{-1}_i$ and $\overline{\eta'}^{-1}_i$ by
\begin{equation}\label{inveta_i}
        \begin{array}{ll}
          \eta^{-1}_i=\left(\varphi^{-1}_{\langle i,i\rangle}\right)^2,
          &\text{if } \sigma_i=i\text{ and }\sigma_{i+1}=i+1, \\[8pt]
          \overline{\eta}^{-1}_i=\left(\overline{\varphi}^{-1}_{\langle i,i\rangle}\right)^2,
          &\text{if } \sigma_i=\overline{i}\text{ and }\sigma_{i+1}=\overline{i+1}, \\[8pt]
          \eta'^{-1}_i=\varphi^{-1}_{\langle i,i\rangle}\circ\varphi^{-1}_{\langle i+1,i\rangle},
          &\text{if } \sigma_i=i+1\text{ and }\sigma_{i+1}=i, \\[8pt]
          \overline{\eta'}^{-1}_i=\overline{\varphi}^{-1}_{\langle i,i\rangle}\circ\overline{\varphi}^{-1}_{\langle i+1,i\rangle},
          &\text{if } \sigma_i=\overline{i+1}\text{ and }\sigma_{i+1}=\overline{i}.
        \end{array}
\end{equation}
Hence from the view of the grid $P_\sigma$,
the operation $\eta^{-1}_i$ (resp. $\overline{\eta}^{-1}_i$, $\eta'^{-1}_i$, $\overline{\eta'}^{-1}_i$)
deletes the $12$-pair (resp. $\overline{1}\overline{2}$-pair, $21$-pair, $\overline{2}\overline{1}$-pair) centered at the grid point $(i+1, i+1)$.
It is clear that $\eta^{-1}_i(\sigma),\overline{\eta}^{-1}_i(\sigma), \eta'^{-1}_i(\sigma),\overline{\eta'}^{-1}_i(\sigma)\in\mathcal{I}_n^B$,
and $\pi=\eta^{-1}_i(\sigma)$ (reps. $\overline{\eta}^{-1}_i(\sigma), \eta'^{-1}_i(\sigma),\overline{\eta'}^{-1}_i(\sigma)$)
if $\sigma=\eta'_i(\pi)$ (reps. $\overline{\eta}_i(\pi), \eta'_i(\pi),\overline{\eta'}_i(\pi)$).

\section{Geometric properties of $\ob$-descents and inverse $\ob$-descents}\label{Sec-GeoDes}

In this section, we characterize the geometric properties of $\ob$-descents and $\ob$-idescents in signed permutation grids,
and their changes during the growth processes from $n\times n$ grids
to $(n+1)\times (n+1)$ grids.

\subsection{Descent types of square pairs and $d$-types of grid points}
According to the different forms of $\ob$-descents and $\ob$-idescents,
we group the pairs of filled squares in adjacent rows or adjacent columns
into eight different types.

Since the filled squares in the signed permutation grid are positive and negative ,
there are four different pairs of filled squares in adjacent rows as \emph{des-pairs} can form $\ob$-descents,
which are called \emph{des$^+_+$-pair}, \emph{des$^-_+$-pair}, \emph{des$^+_-$-pair} and \emph{des$^-_-$-pair} as follows.
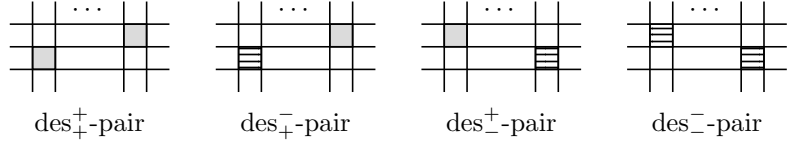
\begin{center}

\end{center}
For $\pi=\pi_1\pi_2\ldots\pi_n\in\mathfrak{B}_n$ with $\pi_0=0$,
the index $i\in [0,n-1]$ is a \emph{$\ob$-ascent} if $\pi_i<\pi_{i+1}$,
and a \emph{$\ob$-iascent} if $\pi^{-1}_i<\pi^{-1}_{i+1}$.
Thus the other four pairs of filled squares in adjacent rows as \emph{asc-pairs} form $\ob$-ascents,
which are called as \emph{asc$^+_+$-pair}, \emph{asc$^-_+$-pair}, \emph{asc$^+_-$-pair} and \emph{asc$^-_-$-pair}.
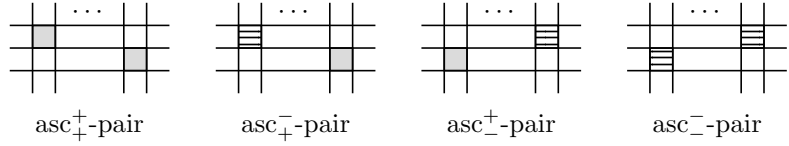
\begin{center}
%
\end{center}

It can be directly seen that
among the eight different pairs of filled squares in adjacent columns,
four of them as \emph{ides-pairs} reflect the $\ob$-idescents,
and the remanning four as \emph{iasc-pairs} reflect the $\ob$-iascents.
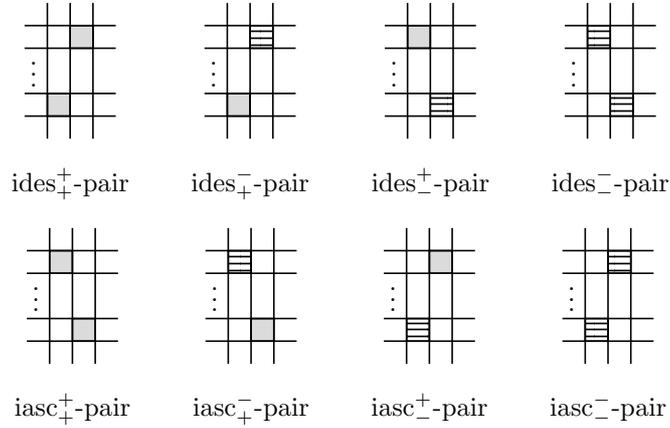
\begin{center}
%
\end{center}

Furthermore, since $\pi_0=0$ for any $\pi\in\mathfrak{B}_n$,
the filled squares in the first row or column of the grid $P_\pi$
are classified into the following four different types as \emph{des$^-$-square}, \emph{asc$^+$-square}, \emph{ides$^-$-square} and \emph{iasc$^+$-square}.
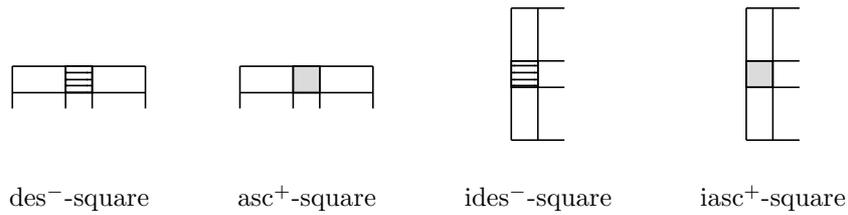
\begin{center}
%
\end{center}


\begin{definition}\label{defd+type}
Let $\pi\in\mathfrak{B}_n$ and $P_\pi$ be its grid, then for $1\leq i,j\leq n+1$,
the grid point $(i,j)$ is of \emph{$d^+$-type} $(p,q)$
and denoted by $d^+(i,j)=(p,q)$ if
$$(\des^B(\varphi_{(i,j)}(\pi)),\ides^B(\varphi_{(i,j)}(\pi)))-(\des^B(\pi),\ides^B(\pi))=(p,q),$$
where $p$ and $q$ are called as the \emph{$d^+_h$-type} and  the \emph{$d^+_v$-type},
and denoted by $d^+_h(i,j)=p$ and $d^+_v(i,j)=q$, respectively.
\end{definition}

\begin{definition}\label{defd-type}
Let $\pi\in\mathfrak{B}_n$ and $P_\pi$ be its grid, then for $1\leq i,j\leq n+1$,
the grid point $(i,j)$ is of \emph{$d^-$-type} $(p',q')$
and denoted by $d^-(i,j)=(p',q')$ if
$$(\des^B(\overline{\varphi}_{(i,j)}(\pi)),\ides^B(\overline{\varphi}_{(i,j)}(\pi)))
-(\des^B(\pi),\ides^B(\pi))=(p',q'),$$
where $p'$ and $q'$ are called as the \emph{$d^-_h$-type} and the \emph{$d^-_v$-type},
and denoted by $d^-_h(i,j)=p'$ and $d^-_v(i,j)=q'$, respectively.
\end{definition}

We call $d^+$-type and $d^-$-type,
$d^+_h$-type and $d^-_h$-type, $d^+_v$-type and $d^-_v$-type
collectively as \emph{$d$-type}, \emph{$d_h$-type}, \emph{$d_v$-type}, respectively,
and use the superscripts $+$ and $-$ to represent the signs of the corresponding type.
We see from the following propositions that the distribution of $d_h$-types is affected by des-pairs and asc-pairs,
and the distribution of $d_v$-types is affected by ides-pairs and iasc-pairs.

\begin{prop}\label{prop-dhdes}
The $d_h$-types of the grid points on the middle horizontal  line of
des-pairs and asc-pairs are indicated in the figure below,
where the numbers above are $d_h^+$-types,
and the numbers below are $d_h^-$-types.
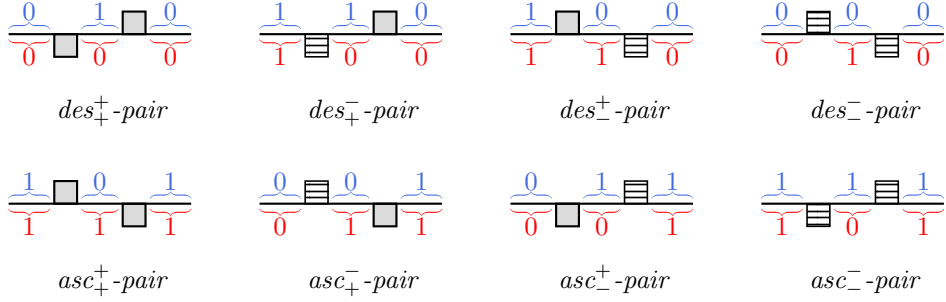
\begin{center}
\begin{tikzpicture}[scale =0.6]
\def\judy{-- +(5mm,0mm) -- +(5mm,5mm) -- +(0mm,5mm) -- cycle [line width=0.8pt, fill=gainsboro]}
\def\nbx{-- +(5mm,0mm) -- +(5mm,5mm) -- +(0mm,5mm) -- cycle [line width=0.6pt,pattern={Lines[angle=45,distance=2.5pt, line width=0.6pt]}, pattern color=black]}
\def\bbrace{[decorate,pen colour={royalblue},
    decoration = {calligraphic brace,mirror,raise=0.7mm,aspect=0.5}]}
\def\rbrace{[decorate,pen colour={red},
    decoration = {calligraphic brace,mirror,raise=0.7mm,aspect=0.5}]}


\draw(0mm,0mm)\judy;
\draw(15mm,5mm)\judy;
\draw[line width=0.8pt](-10mm,5mm)--(30mm,5mm);


\draw\bbrace  (-1mm,5mm) --  (-10mm,5mm);
\node at (-5.5mm,10mm){\royalblue{$0$}};

\draw\bbrace (14mm,5mm) --  (6mm,5mm);
\node at (10mm,10mm){\royalblue{$1$}};

\draw\bbrace (30mm,5mm) --  (21mm,5mm);
\node at (25.5mm,10mm){\royalblue{$0$}};


\draw\rbrace (-10mm,5mm) --  (-1mm,5mm);
\node at (-5.5mm,0mm){\red{$0$}};

\draw\rbrace (6mm,5mm) --  (14mm,5mm);
\node at (10mm,0mm){\red{$0$}};

\draw\rbrace (21mm,5mm) --  (30mm,5mm);
\node at (25.5mm,0mm){\red{$0$}};

\node at (12.5mm,-12mm){des$^+_+$-pair};


\begin{scope}[shift={(55mm,0mm)}]
\draw(0mm,0mm)\nbx;
\draw(15mm,5mm)\judy;
\draw[line width=0.8pt](-10mm,5mm)--(30mm,5mm);


\draw\bbrace (-1mm,5mm) --  (-10mm,5mm);
\node at (-5.5mm,10mm){\royalblue{$1$}};

\draw\bbrace (14mm,5mm) --  (6mm,5mm);
\node at (10mm,10mm){\royalblue{$1$}};

\draw\bbrace (30mm,5mm) --  (21mm,5mm);
\node at (25.5mm,10mm){\royalblue{$0$}};

\draw\rbrace (-10mm,5mm) --  (-1mm,5mm);
\node at (-5.5mm,0mm){\red{$1$}};

\draw\rbrace (6mm,5mm) --  (14mm,5mm);
\node at (10mm,0mm){\red{$0$}};

\draw\rbrace (21mm,5mm) --  (30mm,5mm);
\node at (25.5mm,0mm){\red{$0$}};
\node at (12.5mm,-12mm){des$^-_+$-pair};
\end{scope}


\begin{scope}[shift={(110mm,0mm)}]
\draw(0mm,5mm)\judy;
\draw(15mm,0mm)\nbx;
\draw[line width=0.8pt](-10mm,5mm)--(30mm,5mm);


\draw\bbrace (-1mm,5mm) --  (-10mm,5mm);
\node at (-5.5mm,10mm){\royalblue{$1$}};

\draw\bbrace (14mm,5mm) --  (6mm,5mm);
\node at (10mm,10mm){\royalblue{$0$}};

\draw\bbrace (30mm,5mm) --  (21mm,5mm);
\node at (25.5mm,10mm){\royalblue{$0$}};


\draw\rbrace (-10mm,5mm) --  (-1mm,5mm);
\node at (-5.5mm,0mm){\red{$1$}};

\draw\rbrace (6mm,5mm) --  (14mm,5mm);
\node at (10mm,0mm){\red{$1$}};

\draw\rbrace (21mm,5mm) --  (30mm,5mm);
\node at (25.5mm,0mm){\red{$0$}};
\node at (12.5mm,-12mm){des$^+_-$-pair};
\end{scope}


\begin{scope}[shift={(165mm,0mm)}]
\draw(0mm,5mm)\nbx;
\draw(15mm,0mm)\nbx;
\draw[line width=0.8pt](-10mm,5mm)--(30mm,5mm);


\draw\bbrace (-1mm,5mm) --  (-10mm,5mm);
\node at (-5.5mm,10mm){\royalblue{$0$}};

\draw\bbrace (14mm,5mm) --  (6mm,5mm);
\node at (10mm,10mm){\royalblue{$0$}};

\draw\bbrace (30mm,5mm) --  (21mm,5mm);
\node at (25.5mm,10mm){\royalblue{$0$}};


\draw\rbrace (-10mm,5mm) --  (-1mm,5mm);
\node at (-5.5mm,0mm){\red{$0$}};

\draw\rbrace (6mm,5mm) --  (14mm,5mm);
\node at (10mm,0mm){\red{$1$}};

\draw\rbrace (21mm,5mm) --  (30mm,5mm);
\node at (25.5mm,0mm){\red{$0$}};

\node at (12.5mm,-12mm){des$^-_-$-pair};
\end{scope}
\end{tikzpicture}
\end{center}
\begin{center}
\begin{tikzpicture}[scale =0.6]
\def\judy{-- +(5mm,0mm) -- +(5mm,5mm) -- +(0mm,5mm) -- cycle [line width=0.8pt, fill=gainsboro]}
\def\nbx{-- +(5mm,0mm) -- +(5mm,5mm) -- +(0mm,5mm) -- cycle [line width=0.6pt,pattern={Lines[angle=45,distance=2.5pt, line width=0.6pt]}, pattern color=black]}
\def\bbrace{[decorate,pen colour={royalblue},
    decoration = {calligraphic brace,mirror,raise=0.7mm,aspect=0.5}]}
\def\rbrace{[decorate,pen colour={red},
    decoration = {calligraphic brace,mirror,raise=0.7mm,aspect=0.5}]}


\draw(0mm,5mm)\judy;
\draw(15mm,0mm)\judy;
\draw[line width=0.8pt](-10mm,5mm)--(30mm,5mm);


\draw\bbrace  (-1mm,5mm) --  (-10mm,5mm);
\node at (-5.5mm,10mm){\royalblue{$1$}};

\draw\bbrace (14mm,5mm) --  (6mm,5mm);
\node at (10mm,10mm){\royalblue{$0$}};

\draw\bbrace (30mm,5mm) --  (21mm,5mm);
\node at (25.5mm,10mm){\royalblue{$1$}};


\draw\rbrace (-10mm,5mm) --  (-1mm,5mm);
\node at (-5.5mm,0mm){\red{$1$}};

\draw\rbrace (6mm,5mm) --  (14mm,5mm);
\node at (10mm,0mm){\red{$1$}};

\draw\rbrace (21mm,5mm) --  (30mm,5mm);
\node at (25.5mm,0mm){\red{$1$}};

\node at (12.5mm,-12mm){asc$^+_+$-pair};


\begin{scope}[shift={(55mm,0mm)}]
\draw(0mm,5mm)\nbx;
\draw(15mm,0mm)\judy;
\draw[line width=0.8pt](-10mm,5mm)--(30mm,5mm);


\draw\bbrace (-1mm,5mm) --  (-10mm,5mm);
\node at (-5.5mm,10mm){\royalblue{$0$}};

\draw\bbrace (14mm,5mm) --  (6mm,5mm);
\node at (10mm,10mm){\royalblue{$0$}};

\draw\bbrace (30mm,5mm) --  (21mm,5mm);
\node at (25.5mm,10mm){\royalblue{$1$}};

\draw\rbrace (-10mm,5mm) --  (-1mm,5mm);
\node at (-5.5mm,0mm){\red{$0$}};

\draw\rbrace (6mm,5mm) --  (14mm,5mm);
\node at (10mm,0mm){\red{$1$}};

\draw\rbrace (21mm,5mm) --  (30mm,5mm);
\node at (25.5mm,0mm){\red{$1$}};
\node at (12.5mm,-12mm){asc$^-_+$-pair};
\end{scope}


\begin{scope}[shift={(110mm,0mm)}]
\draw(0mm,0mm)\judy;
\draw(15mm,5mm)\nbx;
\draw[line width=0.8pt](-10mm,5mm)--(30mm,5mm);


\draw\bbrace (-1mm,5mm) --  (-10mm,5mm);
\node at (-5.5mm,10mm){\royalblue{$0$}};

\draw\bbrace (14mm,5mm) --  (6mm,5mm);
\node at (10mm,10mm){\royalblue{$1$}};

\draw\bbrace (30mm,5mm) --  (21mm,5mm);
\node at (25.5mm,10mm){\royalblue{$1$}};


\draw\rbrace (-10mm,5mm) --  (-1mm,5mm);
\node at (-5.5mm,0mm){\red{$0$}};

\draw\rbrace (6mm,5mm) --  (14mm,5mm);
\node at (10mm,0mm){\red{$0$}};

\draw\rbrace (21mm,5mm) --  (30mm,5mm);
\node at (25.5mm,0mm){\red{$1$}};
\node at (12.5mm,-12mm){asc$^+_-$-pair};
\end{scope}


\begin{scope}[shift={(165mm,0mm)}]
\draw(0mm,0mm)\nbx;
\draw(15mm,5mm)\nbx;
\draw[line width=0.8pt](-10mm,5mm)--(30mm,5mm);


\draw\bbrace (-1mm,5mm) --  (-10mm,5mm);
\node at (-5.5mm,10mm){\royalblue{$1$}};

\draw\bbrace (14mm,5mm) --  (6mm,5mm);
\node at (10mm,10mm){\royalblue{$1$}};

\draw\bbrace (30mm,5mm) --  (21mm,5mm);
\node at (25.5mm,10mm){\royalblue{$1$}};


\draw\rbrace (-10mm,5mm) --  (-1mm,5mm);
\node at (-5.5mm,0mm){\red{$1$}};

\draw\rbrace (6mm,5mm) --  (14mm,5mm);
\node at (10mm,0mm){\red{$0$}};

\draw\rbrace (21mm,5mm) --  (30mm,5mm);
\node at (25.5mm,0mm){\red{$1$}};

\node at (12.5mm,-12mm){asc$^-_-$-pair};
\end{scope}
\end{tikzpicture}
\end{center}
\end{prop}

\begin{prop}\label{prop-dviasc}
The $d_v$-types of the grid points on the middle vertical line of  ides-pairs and iasc-pairs are indicated in the figure below,
where the numbers on the left are $d_v^+$-types,
and the numbers on the right are $d_v^-$-types.
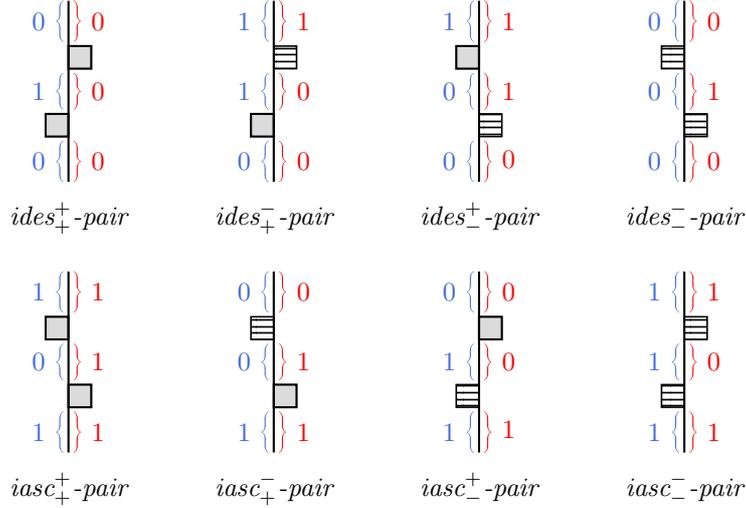
\begin{center}
\begin{tikzpicture}[scale =0.6]
\def\judy{-- +(5mm,0mm) -- +(5mm,5mm) -- +(0mm,5mm) -- cycle [line width=0.8pt, fill=gainsboro]}
\def\nbx{-- +(5mm,0mm) -- +(5mm,5mm) -- +(0mm,5mm) -- cycle [line width=0.6pt,pattern={Lines[angle=45,distance=2.5pt, line width=0.6pt]}, pattern color=black]}
\def\bbrace{[decorate,pen colour={royalblue}, decoration = {calligraphic brace,mirror,raise=0.7mm,aspect=0.5}]}
\def\rbrace{[decorate,pen colour={red}, decoration = {calligraphic brace,mirror,raise=0.7mm,aspect=0.5}]}


\draw(10mm,10mm)\judy;
\draw(5mm,-5mm)\judy;
\draw[line width=0.8pt](10mm,-15mm)--(10mm,25mm);


\draw\bbrace (10mm,25mm) --  (10mm,16mm);
\node at (3.5mm,20.5mm){\royalblue{$0$}};

\draw\bbrace (10mm,9mm) --  (10mm,1mm);
\node at (3.5mm,5mm){\royalblue{$1$}};

\draw\bbrace (10mm,-6mm) --  (10mm,-15mm);
\node at (3.5mm,-10.5mm){\royalblue{$0$}};

\draw\rbrace (10mm,16mm) --  (10mm,25mm);
\node at (16.5mm,20.5mm){\red{$0$}};

\draw\rbrace (10mm,1mm) --  (10mm,9mm);
\node at (16.5mm,5mm){\red{$0$}};

\draw\rbrace (10mm,-15mm) --  (10mm,-6mm);
\node at (16.5mm,-10.5mm){\red{$0$}};

\node at (10mm,-23mm){ides$^+_+$-pair};


\begin{scope}[shift={(45mm,0mm)}]

\draw(10mm,10mm)\nbx;
\draw(5mm,-5mm)\judy;
\draw[line width=0.8pt](10mm,-15mm)--(10mm,25mm);


\draw\bbrace (10mm,25mm) --  (10mm,16mm);
\node at (3.5mm,20.5mm){\royalblue{$1$}};

\draw\bbrace (10mm,9mm) --  (10mm,1mm);
\node at (3.5mm,5mm){\royalblue{$1$}};

\draw\bbrace  (10mm,-6mm) --  (10mm,-15mm);
\node at (3.5mm,-10.5mm){\royalblue{$0$}};

\draw\rbrace  (10mm,16mm) --  (10mm,25mm);
\node at (16.5mm,20.5mm){\red{$1$}};

\draw\rbrace  (10mm,1mm) --  (10mm,9mm);
\node at (16.5mm,5mm){\red{$0$}};

\draw\rbrace (10mm,-15mm) --  (10mm,-6mm);
\node at (16.5mm,-10.5mm){\red{$0$}};

\node at (10mm,-23mm){ides$^-_+$-pair};
\end{scope}


\begin{scope}[shift={(90mm,0mm)}]

\draw(5mm,10mm)\judy;
\draw(10mm,-5mm)\nbx;
\draw[line width=0.8pt](10mm,-15mm)--(10mm,25mm);


\draw\bbrace (10mm,25mm) --  (10mm,16mm);
\node at (3.5mm,20.5mm){\royalblue{$1$}};

\draw\bbrace (10mm,9mm) --  (10mm,1mm);
\node at (3.5mm,5mm){\royalblue{$0$}};

\draw\bbrace (10mm,-6mm) --  (10mm,-15mm);
\node at (3.5mm,-10.5mm){\royalblue{$0$}};


\draw\rbrace (10mm,16mm) --  (10mm,25mm);
\node at (16.5mm,20.5mm){\red{$1$}};

\draw\rbrace (10mm,1mm) --  (10mm,9mm);
\node at (16.5mm,5mm){\red{$1$}};

\draw\rbrace (10mm,-15mm) --  (10mm,-6mm);
\node at (16.5mm,-10mm){\red{$0$}};

\node at (10mm,-23mm){ides$^+_-$-pair};
\end{scope}


\begin{scope}[shift={(135mm,0mm)}]

\draw(5mm,10mm)\nbx;
\draw(10mm,-5mm)\nbx;
\draw[line width=0.8pt](10mm,-15mm)--(10mm,25mm);


\draw\bbrace (10mm,25mm) --  (10mm,16mm);
\node at (3.5mm,20.5mm){\royalblue{$0$}};

\draw\bbrace (10mm,9mm) --  (10mm,1mm);
\node at (3.5mm,5mm){\royalblue{$0$}};

\draw\bbrace (10mm,-6mm) --  (10mm,-15mm);
\node at (3.5mm,-10.5mm){\royalblue{$0$}};


\draw\rbrace (10mm,16mm) --  (10mm,25mm);
\node at (16.5mm,20.5mm){\red{$0$}};

\draw\rbrace (10mm,1mm) --  (10mm,9mm);
\node at (16.5mm,5mm){\red{$1$}};

\draw\rbrace (10mm,-15mm) --  (10mm,-6mm);
\node at (16.5mm,-10.5mm){\red{$0$}};

\node at (10mm,-23mm){ides$^-_-$-pair};
\end{scope}
\end{tikzpicture}
\end{center}
\begin{center}
\begin{tikzpicture}[scale =0.6]
\def\judy{-- +(5mm,0mm) -- +(5mm,5mm) -- +(0mm,5mm) -- cycle [line width=0.8pt, fill=gainsboro]}
\def\nbx{-- +(5mm,0mm) -- +(5mm,5mm) -- +(0mm,5mm) -- cycle [line width=0.6pt,pattern={Lines[angle=45,distance=2.5pt, line width=0.6pt]}, pattern color=black]}
\def\bbrace{[decorate,pen colour={royalblue}, decoration = {calligraphic brace,mirror,raise=0.7mm,aspect=0.5}]}
\def\rbrace{[decorate,pen colour={red}, decoration = {calligraphic brace,mirror,raise=0.7mm,aspect=0.5}]}


\draw(5mm,10mm)\judy;
\draw(10mm,-5mm)\judy;
\draw[line width=0.8pt](10mm,-15mm)--(10mm,25mm);


\draw\bbrace (10mm,25mm) --  (10mm,16mm);
\node at (3.5mm,20.5mm){\royalblue{$1$}};

\draw\bbrace (10mm,9mm) --  (10mm,1mm);
\node at (3.5mm,5mm){\royalblue{$0$}};

\draw\bbrace (10mm,-6mm) --  (10mm,-15mm);
\node at (3.5mm,-10.5mm){\royalblue{$1$}};

\draw\rbrace (10mm,16mm) --  (10mm,25mm);
\node at (16.5mm,20.5mm){\red{$1$}};

\draw\rbrace (10mm,1mm) --  (10mm,9mm);
\node at (16.5mm,5mm){\red{$1$}};

\draw\rbrace (10mm,-15mm) --  (10mm,-6mm);
\node at (16.5mm,-10.5mm){\red{$1$}};

\node at (10mm,-23mm){iasc$^+_+$-pair};


\begin{scope}[shift={(45mm,0mm)}]

\draw(5mm,10mm)\nbx;
\draw(10mm,-5mm)\judy;
\draw[line width=0.8pt](10mm,-15mm)--(10mm,25mm);


\draw\bbrace (10mm,25mm) --  (10mm,16mm);
\node at (3.5mm,20.5mm){\royalblue{$0$}};

\draw\bbrace (10mm,9mm) --  (10mm,1mm);
\node at (3.5mm,5mm){\royalblue{$0$}};

\draw\bbrace  (10mm,-6mm) --  (10mm,-15mm);
\node at (3.5mm,-10.5mm){\royalblue{$1$}};

\draw\rbrace  (10mm,16mm) --  (10mm,25mm);
\node at (16.5mm,20.5mm){\red{$0$}};

\draw\rbrace  (10mm,1mm) --  (10mm,9mm);
\node at (16.5mm,5mm){\red{$1$}};

\draw\rbrace (10mm,-15mm) --  (10mm,-6mm);
\node at (16.5mm,-10.5mm){\red{$1$}};

\node at (10mm,-23mm){iasc$^-_+$-pair};
\end{scope}


\begin{scope}[shift={(90mm,0mm)}]

\draw(10mm,10mm)\judy;
\draw(5mm,-5mm)\nbx;
\draw[line width=0.8pt](10mm,-15mm)--(10mm,25mm);


\draw\bbrace (10mm,25mm) --  (10mm,16mm);
\node at (3.5mm,20.5mm){\royalblue{$0$}};

\draw\bbrace (10mm,9mm) --  (10mm,1mm);
\node at (3.5mm,5mm){\royalblue{$1$}};

\draw\bbrace (10mm,-6mm) --  (10mm,-15mm);
\node at (3.5mm,-10.5mm){\royalblue{$1$}};


\draw\rbrace (10mm,16mm) --  (10mm,25mm);
\node at (16.5mm,20.5mm){\red{$0$}};

\draw\rbrace (10mm,1mm) --  (10mm,9mm);
\node at (16.5mm,5mm){\red{$0$}};

\draw\rbrace (10mm,-15mm) --  (10mm,-6mm);
\node at (16.5mm,-10mm){\red{$1$}};

\node at (10mm,-23mm){iasc$^+_-$-pair};
\end{scope}


\begin{scope}[shift={(135mm,0mm)}]

\draw(10mm,10mm)\nbx;
\draw(5mm,-5mm)\nbx;
\draw[line width=0.8pt](10mm,-15mm)--(10mm,25mm);


\draw\bbrace (10mm,25mm) --  (10mm,16mm);
\node at (3.5mm,20.5mm){\royalblue{$1$}};

\draw\bbrace (10mm,9mm) --  (10mm,1mm);
\node at (3.5mm,5mm){\royalblue{$1$}};

\draw\bbrace (10mm,-6mm) --  (10mm,-15mm);
\node at (3.5mm,-10.5mm){\royalblue{$1$}};


\draw\rbrace (10mm,16mm) --  (10mm,25mm);
\node at (16.5mm,20.5mm){\red{$1$}};

\draw\rbrace (10mm,1mm) --  (10mm,9mm);
\node at (16.5mm,5mm){\red{$0$}};

\draw\rbrace (10mm,-15mm) --  (10mm,-6mm);
\node at (16.5mm,-10.5mm){\red{$1$}};

\node at (10mm,-23mm){iasc$^-_-$-pair};
\end{scope}
\end{tikzpicture}
\end{center}
\end{prop}

The above two propositions can be directly checked by taking operations $\varphi_{(i,j)}$
and $\overline{\varphi}_{(i,j)}$ at corresponding grid points $(i,j)$.
Using the same analysis with $\pi_0=0$, we have
the next proposition that  describes the distributions of $d_h$-types and $d_v$-types of grid points on the grid border, and it can be

%
%
%

\begin{prop}\label{prop-dtype-rim}
The $d_h$-types of the grid points on the top and bottom boundaries of the grid
are indicated in the figure below,
where the numbers above are $d_h^+$-types,
and the numbers below are $d_h^-$-types.
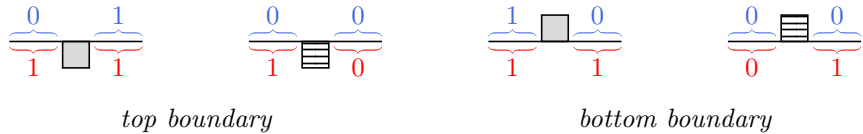
\begin{center}
\begin{tikzpicture}[scale =0.7]
\def\nbx{-- +(5mm,0mm) -- +(5mm,5mm) -- +(0mm,5mm) -- cycle [line width=0.6pt,pattern={Lines[angle=45,distance=2.5pt, line width=0.6pt]}, pattern color=black]}
\def\judy{-- +(5mm,0mm) -- +(5mm,5mm) -- +(0mm,5mm) -- cycle [line width=0.6pt,fill=gainsboro]}
\def\bbrace{[decorate,pen colour={royalblue}, decoration = {calligraphic brace,mirror,raise=0.7mm,aspect=0.5}]}
\def\rbrace{[decorate,pen colour={red}, decoration = {calligraphic brace,mirror,raise=0.7mm,aspect=0.5}]}

\draw (10mm,0mm)\judy;
\draw[line width=0.6pt] (0mm,5mm)--(25mm,5mm);

\draw\bbrace (9mm,5mm) --  (0mm,5mm);
\node at (4.5mm,10mm){\royalblue{$0$}};
\draw\bbrace (25mm,5mm) --  (16mm,5mm);
\node at (20.5mm,10mm){\royalblue{$1$}};

\draw\rbrace (0mm,5mm) --  (9mm,5mm);
\node at (4.5mm,0mm){\red{$1$}};
\draw\rbrace (16mm,5mm) --  (25mm,5mm);
\node at (20.5mm,0mm){\red{$1$}};

\draw (55mm,0mm)\nbx;
\draw[line width=0.6pt] (45mm,5mm)--(70mm,5mm);

\draw\bbrace (54mm,5mm) --  (45mm,5mm);
\node at (49.5mm,10mm){\royalblue{$0$}};
\draw\bbrace (70mm,5mm) --  (61mm,5mm);
\node at (65.5mm,10mm){\royalblue{$0$}};

\draw\rbrace (45mm,5mm) --  (54mm,5mm);
\node at (49.5mm,0mm){\red{$1$}};
\draw\rbrace (61mm,5mm) --  (70mm,5mm);
\node at (65.5mm,0mm){\red{$0$}};
\node at (35mm,-10mm){top boundary};
\begin{scope}[xshift=90mm]
\draw (10mm,5mm)\judy;
\draw[line width=0.6pt] (0mm,5mm)--(25mm,5mm);

\draw\bbrace (9mm,5mm) --  (0mm,5mm);
\node at (4.5mm,10mm){\royalblue{$1$}};
\draw\bbrace (25mm,5mm) --  (16mm,5mm);
\node at (20.5mm,10mm){\royalblue{$0$}};

\draw\rbrace (0mm,5mm) --  (9mm,5mm);
\node at (4.5mm,0mm){\red{$1$}};
\draw\rbrace (16mm,5mm) --  (25mm,5mm);
\node at (20.5mm,0mm){\red{$1$}};
\draw (55mm,5mm)\nbx;
\draw[line width=0.6pt] (45mm,5mm)--(70mm,5mm);

\draw\bbrace (54mm,5mm) --  (45mm,5mm);
\node at (49.5mm,10mm){\royalblue{$0$}};
\draw\bbrace (70mm,5mm) --  (61mm,5mm);
\node at (65.5mm,10mm){\royalblue{$0$}};

\draw\rbrace (45mm,5mm) --  (54mm,5mm);
\node at (49.5mm,0mm){\red{$0$}};
\draw\rbrace (61mm,5mm) --  (70mm,5mm);
\node at (65.5mm,0mm){\red{$1$}};
\node at (35mm,-10mm){bottom boundary};
\end{scope}
\end{tikzpicture}
\end{center}
And the $d_v$-types of the grid points on the left and right boundaries of the grid are indicated in the figure below,
where the numbers on the left are $d_v^+$-types,
and the numbers on the right are $d_v^-$-types.
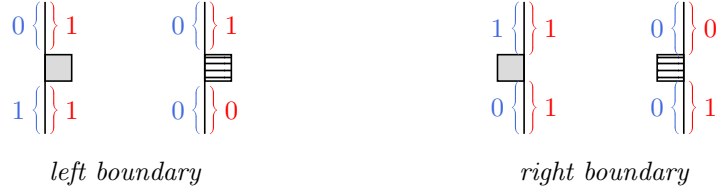
\begin{center}

\begin{tikzpicture}[scale =0.7]
\def\nbx{-- +(5mm,0mm) -- +(5mm,5mm) -- +(0mm,5mm) -- cycle [line width=0.6pt,pattern={Lines[angle=45,distance=2.5pt, line width=0.6pt]}, pattern color=black]}
\def\judy{-- +(5mm,0mm) -- +(5mm,5mm) -- +(0mm,5mm) -- cycle [line width=0.6pt,fill=gainsboro]}
\def\bbrace{[decorate,pen colour={royalblue}, decoration = {calligraphic brace,mirror,raise=0.7mm,aspect=0.5}]}
\def\rbrace{[decorate,pen colour={red}, decoration = {calligraphic brace,mirror,raise=0.7mm,aspect=0.5}]}

\draw (10mm,0mm)\judy;
\draw[line width=0.6pt] (10mm,-10mm)--(10mm,15mm);

\draw\bbrace (10mm,15mm) --  (10mm,6mm);
\node at (5mm,10.5mm){\royalblue{$0$}};
\draw\bbrace (10mm,-1mm) --  (10mm,-10mm);
\node at (5mm,-5.5mm){\royalblue{$1$}};

\draw\rbrace (10mm,6mm) --  (10mm,15mm);
\node at (15mm,10.5mm){\red{$1$}};
\draw\rbrace (10mm,-10mm) --  (10mm,-1mm);
\node at (15mm,-5.5mm){\red{$1$}};

\draw(40mm,0mm)\nbx;
\draw[line width=0.6pt] (40mm,-10mm)--(40mm,15mm);


\draw\bbrace (40mm,15mm) --  (40mm,6mm);
\node at (35mm,10.5mm){\royalblue{$0$}};
\draw\bbrace (40mm,-1mm) --  (40mm,-10mm);
\node at (35mm,-5.5mm){\royalblue{$0$}};

\draw\rbrace (40mm,6mm) --  (40mm,15mm);
\node at (45mm,10.5mm){\red{$1$}};
\draw\rbrace (40mm,-10mm) --  (40mm,-1mm);
\node at (45mm,-5.5mm){\red{$0$}};
\node at (25mm,-17.5mm){left boundary};
\begin{scope}[xshift=90mm]

\draw (5mm,0mm)\judy;
\draw[line width=0.6pt] (10mm,-10mm)--(10mm,15mm);

\draw\bbrace (10mm,15mm) --  (10mm,5mm);
\node at (5mm,10mm){\royalblue{$1$}};
\draw\bbrace (10mm,0mm) --  (10mm,-10mm);
\node at (5mm,-5mm){\royalblue{$0$}};

\draw\rbrace (10mm,5mm) --  (10mm,15mm);
\node at (15mm,10mm){\red{$1$}};
\draw\rbrace (10mm,-10mm) --  (10mm,0mm);
\node at (15mm,-5mm){\red{$1$}};

\draw(35mm,0mm)\nbx;
\draw[line width=0.6pt] (40mm,-10mm)--(40mm,15mm);


\draw\bbrace (40mm,15mm) --  (40mm,5mm);
\node at (35mm,10mm){\royalblue{$0$}};
\draw\bbrace (40mm,0mm) --  (40mm,-10mm);
\node at (35mm,-5mm){\royalblue{$0$}};

\draw\rbrace (40mm,5mm) --  (40mm,15mm);
\node at (45mm,10mm){\red{$0$}};
\draw\rbrace (40mm,-10mm) --  (40mm,0mm);
\node at (45mm,-5mm){\red{$1$}};
\node at (25mm,-17.5mm){right boundary};

\end{scope}
\end{tikzpicture}
\end{center}

\end{prop}
%
%
%

Combining Propositions \ref{prop-dhdes}--\ref{prop-dtype-rim},
we obtain all possible values of $d$-types,
and the $d$-types of the four corners of filled squares.

\begin{prop}\label{prop-possible}
Let $\pi\in\B_n$ and $P_\pi$ be its grid, then
all possible $d^+$-types or $d^-$-types of each grid point are
only $(0,0)$, $(1,0),(0,1)$ and $(1,1)$.
\end{prop}

\begin{prop}\label{prop-4corners}
Let $\pi\in\B_n$ and $P_\pi$ be its grid, then for $1\leq i,j\leq n$, we have
\[
d^+(i,j)=d^+(i+1,j+1)=(0,0)\quad\mbox{and}\quad d^+(i,j+1)=d^+(i+1,j)=(1,1),
\]
if $\ij$ is a positive square in $P_\pi$, and
\[
d^-(i,j)=d^-(i+1,j+1)=(1,1)\quad\mbox{and}\quad d^-(i,j+1)=d^-(i+1,j)=(0,0).
\]
if $\ij$ is a negative square in $P_\pi$, as shown below.
\begin{center}
\begin{tikzpicture}[scale=1.5]
\def\judy{-- +(5mm,0mm) -- +(5mm,5mm) -- +(0mm,5mm) -- cycle [line width=1pt,fill=gainsboro]}
\tikzstyle{gdot}=[circle,fill=royalblue,draw=royalblue,inner sep=2]

\draw(15mm,0mm)\judy;
\draw[line width=1pt,dashed](10mm,5mm)--(25mm,5mm);
\draw[line width=1pt,dashed](10mm,0mm)--(25mm,0mm);
\draw[line width=1pt,dashed] (15mm,-5mm)--(15mm,10mm);
\draw[line width=1pt,dashed] (20mm,-5mm)--(20mm,10mm);

\node[gdot] at (15mm,5mm){};
\node at (12mm,7mm){\royalblue{$(0,0)$}};

\node[gdot] at (20mm,0mm){};
\node at (23mm,-2mm){\royalblue{$(0,0)$}};

\node[gdot] at (15mm,0mm){};
\node at (12mm,-2mm){\royalblue{$(1,1)$}};

\node[gdot] at (20mm,5mm){};
\node at (23mm,7mm){\royalblue{$(1,1)$}};

\node at (17.5mm,-10mm){$d^+$-types on a positive square};

\end{tikzpicture}
\hspace{3em}
\begin{tikzpicture}[scale =1.5]
\def\nbx{-- +(5mm,0mm) -- +(5mm,5mm) -- +(0mm,5mm) -- cycle [line width=1pt,pattern={Lines[angle=45,distance=4pt, line width=1pt]}, pattern color=black]}
\tikzstyle{rdot}=[circle,fill=red,draw=red,inner sep=2]

\draw(15mm,0mm)\nbx;
\draw[line width=1pt,dashed](10mm,5mm)--(25mm,5mm);
\draw[line width=1pt,dashed](10mm,0mm)--(25mm,0mm);
\draw[line width=1pt,dashed] (15mm,-5mm)--(15mm,10mm);
\draw[line width=1pt,dashed] (20mm,-5mm)--(20mm,10mm);

\node[rdot] at (15mm,5mm){};
\node at (12mm,7mm){\red{$(1,1)$}};

\node[rdot] at (20mm,0mm){};
\node at (23mm,-2mm){\red{$(1,1)$}};

\node[rdot] at (15mm,0mm){};
\node at (12mm,-2mm){\red{$(0,0)$}};

\node[rdot] at (20mm,5mm){};
\node at (23mm,7mm){\red{$(0,0)$}};

\node at (17.5mm,-10mm){$d^-$-types on a negative square};

\end{tikzpicture}
\end{center}
\end{prop}

Additionally, we arrive at  more general conclusions that describe the circumstances under which the $d_h$-types or $d_v$-types on the same horizontal or vertical grid line change.

\begin{prop}\label{prop-dh=!}
Let $\pi\in\mathfrak{B}_n$ and $P_\pi$ be its grid, then for any two grid points $(i,j_1)$ and $(i,j_2)$ on the $i$-th horizontal grid line with $1\leq i\leq n+1$ and $1\leq j_1<j_2\leq n+1$, we have
\begin{enumerate}[(a)]
  \item $d_h^+(i,j_1)\neq d_h^+(i,j_1)$ if and only if there exists
  exactly one positive square $\langle i', j'\rangle$ satisfying $i-1\leq i'\leq i$ and $j_1\leq j< j_2$,
  \item $d_h^-(i,j_1)\neq d_h^-(i,j_1)$ if and only if there exists
  exactly one negative square $\langle i', j'\rangle$ satisfying $i-1\leq i'\leq i$ and $j_1\leq j< j_2$.
\end{enumerate}
\end{prop}

\begin{prop}\label{prop-dv=!}
Let $\pi\in\mathfrak{B}_n$ and $P_\pi$ be its grid, then for any two grid points $(i_1,j)$ and $(i_2,j)$ on the $j$-th vertical grid line with $1\leq j\leq n+1$ and $1\leq i_1<i_2\leq n+1$, we have
\begin{enumerate}[(a)]
  \item $d_v^+(i_1,j)\neq d_h^+(i_2,j)$ if and only if there exists
  exactly one positive square $\langle i', j'\rangle$ satisfying $j-1\leq j' \leq j$ and $i_1\leq i< i_2$,
  \item $d_v^-(i_1,j)\neq d_h^-(i_2,j)$ if and only if there exists
  exactly one negative square $\langle i', j'\rangle$ satisfying $j-1\leq j' \leq j$ and $i_1\leq i< i_2$.
\end{enumerate}
\end{prop}

\subsection{$p_h$-paths and $q_v$-paths for $p,q\in\{0,1\}$}

As a preparation for counting the number of grid points of certain $d$-types,
we connect grid points with the same $d_h$-type $p$ or $d_v$-type $q$
for $p,q\in\{0,1\}$ by specific rules
to build $p_h$-paths or $q_v$-paths.

\begin{definition}\label{def-path}
Let $\pi\in\mathfrak{B}_n$ and $P_\pi$ be its grid, then for $p,q\in\{0,1\}$,
\begin{enumerate}[(a)]
\item the \emph{$p_h^+$-paths} (reps. \emph{$p_h^-$-paths}) are established by
connecting grid points of $d_h^+$-types (resp. $d_h^-$-types) $p$  only along horizontal grid lines,
and then connecting these segments of length less than $n+1$ by diagonal lines
only if they touch the same filled square;
\item the \emph{$q_v^+$-paths} (reps. \emph{$q_v^-$-paths}) are established by
connecting grid points of $d_v^+$-types (resp. $d_v^-$-types) $q$ only along vertical grid lines,
and then connecting these segments of length less than $n+1$ by diagonal lines
only if they touch the same filled square.
\end{enumerate}
\end{definition}

We refer $p_h^+$-paths and $p_h^-$-paths collectively as $p_h$-paths,
and refer $q_v^+$-paths and $q_v^-$-paths collectively as $q_v$-paths.
By the construction above, we know that for any grid point, it must be on exactly
one $p_h^+$-path, one $p_h^-$-path, one $q_v^+$-path and one $q_v^-$-path.

\begin{example}
  The following figures present each one of $p_x^\ast$-path
  in the permutation grid  of $\pi=\bar{3}1652\bar{4}$,
  where $p\in\{0,1\}$, $x\in\{h,v\}$ and $\ast\in\{+,-\}$.
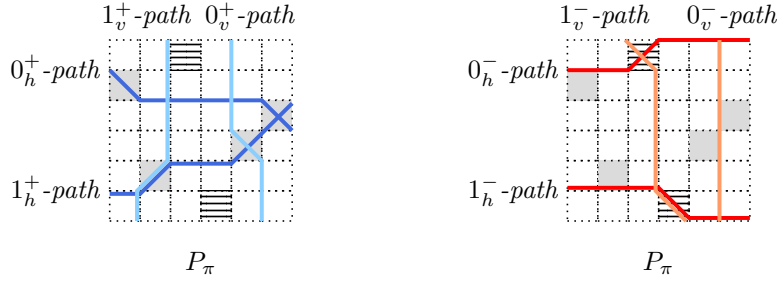
\begin{center}
\begin{tikzpicture}[scale =0.8]
\def\hezi{-- +(5mm,0mm) -- +(5mm,5mm) -- +(0mm,5mm) -- cycle [line width=0.6pt, dotted]}
\def\judy{-- +(5mm,0mm) -- +(5mm,5mm) -- +(0mm,5mm) -- cycle [line width=0.6pt,fill=gainsboro,dotted]}
\def\nbx{-- +(5mm,0mm) -- +(5mm,5mm) -- +(0mm,5mm) -- cycle [line width=0.6pt,pattern={Lines[angle=45,distance=2.5pt, line width=0.6pt]}, pattern color=black,dotted]}
\draw (0mm,0mm)\hezi;
\draw (5mm,0mm)\hezi;
\draw (10mm,0mm)\nbx;
\draw (15mm,0mm)\hezi;
\draw (20mm,0mm)\hezi;
\draw (25mm,0mm)\hezi;
\draw (0mm,-5mm)\judy;
\draw (5mm,-5mm)\hezi;
\draw (10mm,-5mm)\hezi;
\draw (15mm,-5mm)\hezi;
\draw (20mm,-5mm)\hezi;
\draw (25mm,-5mm)\hezi;
\draw (0mm,-10mm)\hezi;
\draw (5mm,-10mm)\hezi;
\draw (10mm,-10mm)\hezi;
\draw (15mm,-10mm)\hezi;
\draw (20mm,-10mm)\hezi;
\draw (25mm,-10mm)\judy;
\draw (0mm,-15mm)\hezi;
\draw (5mm,-15mm)\hezi;
\draw (10mm,-15mm)\hezi;
\draw (15mm,-15mm)\hezi;
\draw (20mm,-15mm)\judy;
\draw (25mm,-15mm)\hezi;
\draw (0mm,-20mm)\hezi;
\draw (5mm,-20mm)\judy;
\draw (10mm,-20mm)\hezi;
\draw (15mm,-20mm)\hezi;
\draw (20mm,-20mm)\hezi;
\draw (25mm,-20mm)\hezi;
\draw (0mm,-25mm)\hezi;
\draw (5mm,-25mm)\hezi;
\draw (10mm,-25mm)\hezi;
\draw (15mm,-25mm)\nbx;
\draw (20mm,-25mm)\hezi;
\draw (25mm,-25mm)\hezi;
\node at (-9mm,0mm){$0_h^+$-path};
\draw(0mm,0mm)--(5mm,-5mm)--(25mm,-5mm)--(30mm,-10mm)[line width=1.5pt,draw=royalblue];
\node at (-9mm,-20mm){$1_h^+$-path};
\draw[shift={(0,-0.05)}] (0mm,-20mm)--(5mm,-20mm)--(10mm,-15mm)--(20mm,-15mm)--(25mm,-10mm)--(30mm,-5mm)[line width=1.5pt,draw=royalblue];

\node at (23mm,9mm){$0_v^+$-path};
\draw(20mm,5mm)--(20mm,-10mm)--(25mm,-15mm)--(25mm,-25mm)[line width=1.5pt,draw=babyblue];
\node at (6mm,9mm){$1_v^+$-path};
\draw[shift={(-0.05,0)}] (10mm,5mm)--(10mm,-15mm)--(5mm,-20mm)--(5mm,-25mm)[line width=1.5pt,draw=babyblue];

\node at (15mm,-32mm){$P_\pi$};
\end{tikzpicture}
\hspace{5em}
\begin{tikzpicture}[scale =0.8]
\def\hezi{-- +(5mm,0mm) -- +(5mm,5mm) -- +(0mm,5mm) -- cycle [line width=0.6pt, dotted]}
\def\judy{-- +(5mm,0mm) -- +(5mm,5mm) -- +(0mm,5mm) -- cycle [line width=0.6pt,fill=gainsboro,dotted]}
\def\nbx{-- +(5mm,0mm) -- +(5mm,5mm) -- +(0mm,5mm) -- cycle [line width=0.6pt,pattern={Lines[angle=45,distance=2.5pt, line width=0.6pt]}, pattern color=black,dotted]}
\tikzstyle{cc}=[circle,draw=black,fill=yellow, line width=0.5pt, inner sep=1.5]
\draw (0mm,0mm)\hezi;
\draw (5mm,0mm)\hezi;
\draw (10mm,0mm)\nbx;
\draw (15mm,0mm)\hezi;
\draw (20mm,0mm)\hezi;
\draw (25mm,0mm)\hezi;
\draw (0mm,-5mm)\judy;
\draw (5mm,-5mm)\hezi;
\draw (10mm,-5mm)\hezi;
\draw (15mm,-5mm)\hezi;
\draw (20mm,-5mm)\hezi;
\draw (25mm,-5mm)\hezi;
\draw (0mm,-10mm)\hezi;
\draw (5mm,-10mm)\hezi;
\draw (10mm,-10mm)\hezi;
\draw (15mm,-10mm)\hezi;
\draw (20mm,-10mm)\hezi;
\draw (25mm,-10mm)\judy;
\draw (0mm,-15mm)\hezi;
\draw (5mm,-15mm)\hezi;
\draw (10mm,-15mm)\hezi;
\draw (15mm,-15mm)\hezi;
\draw (20mm,-15mm)\judy;
\draw (25mm,-15mm)\hezi;
\draw (0mm,-20mm)\hezi;
\draw (5mm,-20mm)\judy;
\draw (10mm,-20mm)\hezi;
\draw (15mm,-20mm)\hezi;
\draw (20mm,-20mm)\hezi;
\draw (25mm,-20mm)\hezi;
\draw (0mm,-25mm)\hezi;
\draw (5mm,-25mm)\hezi;
\draw (10mm,-25mm)\hezi;
\draw (15mm,-25mm)\nbx;
\draw (20mm,-25mm)\hezi;
\draw (25mm,-25mm)\hezi;
\node at (-9mm,0mm){$0_h^-$-path};
\draw(0mm,0mm)--(10mm,0mm)--(15mm,5mm)--(30mm,5mm)[line width=1.5pt,draw=red];

\node at (-9mm,-20mm){$1_h^-$-path};
\draw[shift={(0,0.05)}] (0mm,-20mm)--(15mm,-20mm)--(20mm,-25mm)--(30mm,-25mm)[line width=1.5pt,draw=red];
\node at (6mm,9mm){$1_v^-$-path};
\draw[shift={(-0.05,0)}] (10mm,5mm)--(15mm,0mm)--(15mm,-20mm)--(20mm,-25mm)[line width=1.5pt,draw=orange];

\node at (27mm,9mm){$0_v^-$-path};
\draw(25mm,5mm)--(25mm,-25mm)[line width=1.5pt,draw=orange];


\node at (15mm,-32mm){$P_\pi$};
\end{tikzpicture}
\end{center}

\end{example}

The above example completely present the walking trends of
these $p_h$-paths and the $q_v$-paths for $p,q\in\{0,1\}$.

\begin{thm}\label{lempath}
Let $\pi\in\mathfrak{B}_n$, then in the grid $P_\pi$,
\begin{enumerate}[(a)]
\item each $0_h^+$-path (resp. $1_h^+$-path) goes from the left boundary to the right boundary along the horizontal grid lines except for carrying out a southeast (resp. northeast) step when encountering a positive square;
\item each $0_v^+$-path (resp. $1_v^+$-path) goes from the top boundary to the bottom boundary along the vertical grid lines except for carrying out a southeast (resp. southwest) step when encountering a positive square;
\item each $0_h^-$-path (resp. $1_h^-$-path) goes from the left boundary to the right boundary along the horizontal grid lines except for carrying out a northeast (resp. southeast) step when encountering a negative square;
\item each $0_v^-$-path (resp. $1_v^-$-path) goes from the top boundary to the bottom boundary along the vertical grid lines except for carrying out a southwest (resp. southeast) step when encountering a negative square.
\end{enumerate}
\end{thm}

\pf
By Propositions \ref{prop-dh=!} and \ref{prop-dv=!},
we note that the $d_h$-types (resp. $d_v$-types) on the same horizontal (resp. vertical) grid line change
from $0$ to $1$ or from $1$ to $0$
only when it touches a filled square of the same sign.
Thus in constructing processes in Definition \ref{def-path}, after connecting the grid points of the same $d_h$-type or $d_v$-type
by horizontal or vertical grid lines,
those $d_h$-segments or $d_v$-segments of length less than $n+1$
must touch filled squares with the same sign as the $d_h$-type or $d_v$-type.

Hence, for $p,q\in\{0,1\}$, connecting the $p_h$-segments or $q_v$-segments by diagonal lines in the filled square of the same sign as $p_h$ or $q_v$
is equivalent to that
the $p_h$-path or the $q_v$-path taking a diagonal step when it meets that square,
and the direction of this diagonal step is completely determined by Proposition \ref{prop-4corners}.
\qed

Based on the above observations, for $p,q\in\{0,1\}$,
we utilize the number of $\overline{B}$-descents to count the number of $p_h$-paths,
and the number of $\overline{B}$-idescents to count the number of $q_v$-paths.

\begin{thm}\label{thmpaths}
Let $\pi\in\mathfrak{B}_n$, then in the  grid $P_\pi$,
\begin{enumerate}[(a)]
    \item the number of $0_h^+$-paths is  $\des^B(\pi)+1$,
    \item the number of $0_h^-$-paths is  $\des^B(\pi)$,
    \item the number of $1_h^+$-paths is  $n-\des^B(\pi)$,
    \item the number of $1_h^-$-paths is  $n-\des^B(\pi)+1$;
\end{enumerate}
and
\begin{enumerate}[(a)]
\addtocounter{enumi}{4}
    \item the number of $0_v^+$-paths is  $\ides^B(\pi)+1$,
    \item the number of $0_v^-$-paths is  $\ides^B(\pi)$,
    \item the number of $1_v^+$-paths is  $n-\ides^B(\pi)$,
    \item the number of $1_v^-$-paths is  $n-\ides^B(\pi)+1$.
\end{enumerate}
\end{thm}

\proof
By Theorem \ref{lempath}, for $p,q\in\{0,1\}$,
the number of $p_h$-paths equals the number of grid points $(i,1)$ of $d_h$-type $p$ for $1\leq i\leq n+1$,
and the number of $q_v$-paths equals  the number of grid points $(1,j)$ of $d_v$-type $q$ for $1\leq j\leq n+1$.

We first assume $\pi_i>0$ for all $1\leq i\leq n$.
Hence in the grid $P_\pi$, there are only four types of pairs of filled squares: des$^+_+$-pair, asc$^+_+$-pair, ides$^+_+$-pair and iasc$^+_+$-pair,
which are counted by $\des^B(\pi)$, $n-1-\des^B(\pi)$, $\ides^B(\pi)$ and $n-1-\ides^B(\pi)$, respectively.

It follows from Proposition \ref{prop-dhdes} that
for $2\leq i\leq n$, the numbers of grid points $(i,1)$ with $d_h^+(i,1)=0$, $d_h^-(i,1)=0$, $d_h^+(i,1)=1$ and $d_h^-(i,1)=1$ are $\des^B(\pi)$, $\des^B(\pi)$,
$n-1-\des^B(\pi)$ and $n-1-\des^B(\pi)$, respectively.
Since $\pi_1>0$ and $\pi_n>0$, by Proposition \ref{prop-dtype-rim},
we have $d_h^+(1,1)=0$, $d_h^-(1,1)=1$, $d_h^+(n+1,1)=1$ and $d_h^-(n+1,1)=1$.
Therefore, for $1\leq i\leq n+1$,
the number of grid points $(i,1)$ with $d_h^+(i,1)=0$, $d_h^-(i,1)=0$, $d_h^+(i,1)=1$ and $d_h^-(i,1)=1$ are $\des^B(\pi)+1$, $\des^B(\pi)$,
$n-\des^B(\pi)$ and $n+1-\des^B(\pi)$, respectively.

On the other side, by Proposition \ref{prop-dviasc}, for $2\leq j\leq n$,
the number of grid points $(1,j)$ with $d_v^+(i,1)=0$, $d_v^-(i,1)=0$, $d_v^+(i,1)=1$ and $d_v^-(i,1)=1$ are given by $\ides^B(\pi)$, $\ides^B(\pi)$,
$n-1-\ides^B(\pi)$ and $n-1-\ides^B(\pi)$, respectively.
Note that $d_v^+(1,1)=0$, $d_v^-(1,1)=1$, $d_v^+(1,n+1)=1$ and $d_v^-(1,n+1)=1$ by Proposition \ref{prop-dtype-rim}.
Thus, for $1\leq j\leq n+1$, the number of grid points $(1,j)$ with $d_v^+(1,j)=0$, $d_v^-(1,j)=0$, $d_v^+(1,j)=1$ and $d_v^-(1,j)=1$ are $\ides^B(\pi)+1$, $\ides^B(\pi)$,
$n-\ides^B(\pi)$ and $n+1-\ides^B(\pi)$, respectively.
So the theorem holds
provided all filled squares in the grid $P_\pi$ are positive.


To complete the proof, we are about to show that the quantitative relationships between $p_h$-paths and $\ob$-descents, and between $q_v$-paths and $\ob$-idescents, remain unchanged after any positive square is replaced by a negative square in $P_\pi$.

Let $\ij$ with $1\leq i,j\leq n$  be the replaced square and $\pi'$ be the signed permutation after the replacement.
In particular, when the positive square $\ij$ is replaced,
the effect is limited to the neighborhood formed by itself and the squares on its adjacent rows or columns,
which corresponds to the 16 cases listed below.
We only offer a detailed proof for Case 1,
while the validity of the remaining cases can be also checked.

The symbols in the following proof are explained as:
the filled squares in $P_\pi$ and $P_{\pi'}$ indicate the square $\langle i,j\rangle$
before and after the substitution;
the blank squares in $P_\pi$ and $P_{\pi'}$ indicate the filled squares with undefined signs
 that remain unchanged after the substitution;
the numbers $a/b$ indicates the $d_h^+$-type (resp. $d_v^+$-type) $a$
and $d_h^-$-type (resp. $d_v^-$-type) $b$ of the corresponding grid point on the first column (resp. row);
$\# p_x^\ast$ and $\# {p_x^\ast}'$ denote the numbers of $p_x^\ast$-paths in $P_{\pi}$ and $P_{\pi'}$, respectively, where $p\in\{0,1\}$, $x\in\{h,v\}$ and $\ast\in\{+,-\}$.

\noindent{\bf Case 1}: $\bm{|\pi_{i-1}|<\pi_i<|\pi_{i+1}|}$ {\bf for} $\bm{2\leq i\leq n-1}$.
By Proposition \ref{prop-dhdes}, although the signs of $\pi_{i-1}$ and $\pi_{i+1}$ are not known,
the squares $\langle i-1, |\pi_{i-1}| \rangle$ and $\langle i,\pi_i\rangle$ in $P_\pi$
always changes from an asc-pair to a des-pair in $P_{\pi'}$ but preserve the $d_h^+$-type and $d_h^-$-type of the grid point $(i,1)$,
and the descent type of the pair $\langle i,\pi_i\rangle$ and $\langle i+1, |\pi_{i+1}| \rangle$ remains the same but the $d_h^+$-type and $d_h^-$-type of $(i+1,1)$
both change from $1$ to $0$.
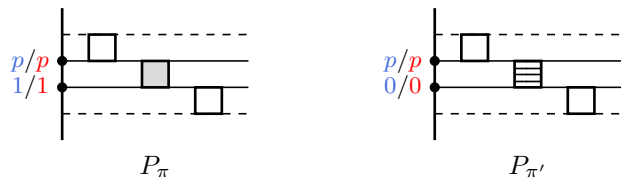
\begin{center}

\end{center}
\begin{center}
\renewcommand\arraystretch{1.2}
\begin{tabular}{m{0.27\textwidth}<{\centering}|m{0.27\textwidth}<{\centering}|m{0.27\textwidth}<{\centering}}
    \hlineB{2.5}
    square pairs & $\pi_{i-1}>0\slash \pi_{i+1}>0$ & $\pi_{i-1}<0\slash\pi_{i+1}<0$\\ \hline
    $\langle i-1, |\pi_{i-1}| \rangle$, $\langle i,\pi_i\rangle$ & asc$^+_+$-pair $\rightarrow$ des$^+_-$-pair & asc$^-_+$-pair $\rightarrow$ des$^-_-$-pair\\ \hline
    $\langle i,\pi_i\rangle$, $\langle i+1, |\pi_{i+1}| \rangle$ & asc$^+_+$-pair $\rightarrow$ asc$^-_+$-pair & des$^+_-$-pair $\rightarrow$ des$^-_-$-pair\\
    \hlineB{2.5}
  \end{tabular}
\end{center}
Thus, we obtain the desired change in the quantitative relationships,
which are $\des^B(\pi')=\des^B(\pi)+1$, and
$$\# {0_h^+}'=\# {0_h^+}+1, \,\,\# {0_h^-}'=\# {0_h^-}+1,\,\,
\# {1_h^+}'=\# {1_h^+}-1, \,\, \# {1_h^-}'=\# {1_h^-}-1.$$

\noindent{\bf Case 2}: $\bm{|\pi_{i+1}|<\pi_i<|\pi_{i-1}|}$ {\bf for} $\bm{2\leq i\leq n-1}$.
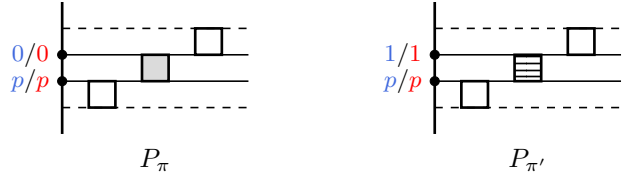
\begin{center}

\end{center}
\begin{center}
\renewcommand\arraystretch{1.2}
\begin{tabular}{m{0.27\textwidth}<{\centering}|m{0.27\textwidth}<{\centering}|m{0.27\textwidth}<{\centering}}
    \hlineB{2.5}
    square pairs & $\pi_{i-1}>0\slash \pi_{i+1}>0$ & $\pi_{i-1}<0\slash\pi_{i+1}<0$\\ \hline
    $\langle i-1, |\pi_{i-1}| \rangle$, $\langle i,\pi_i\rangle$ &
    des$^+_+$-pair $\rightarrow$ des$^-_+$-pair &
    asc$^+_-$-pair $\rightarrow$ asc$^-_-$-pair\\ \hline
    $\langle i,\pi_i\rangle$, $\langle i+1, |\pi_{i+1}| \rangle$ &
    des$^+_+$-pair $\rightarrow$ asc$^+_-$-pair &
    des$^-_+$-pair $\rightarrow$ asc$^-_-$-pair\\
    \hlineB{2.5}
  \end{tabular}
\end{center}
After replacing the positive square $\langle i,\pi_i\rangle$ by a negative square,
we have $\des^B(\pi')=\des^B(\pi)-1$, and
$$\# {0_h^+}'=\# {0_h^+}-1,\,\,\# {0_h^-}'=\# {0_h^-}-1, \,\,
\# {1_h^+}'=\# {1_h^+}+1,\,\, \# {1_h^-}'=\# {1_h^-}+1.$$

\noindent{\bf Case 3}: $\bm{|\pi_{i-1}|<\pi_i}$ {\bf and} $\bm{|\pi_{i+1}|<\pi_i}$ {\bf for} $\bm{2\leq i\leq n-1}$.
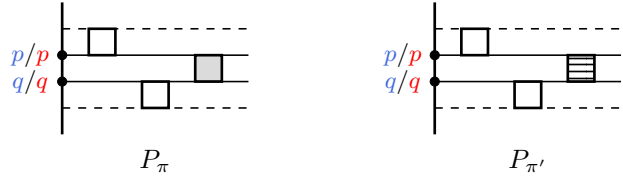
\begin{center}

\end{center}
\begin{center}
\renewcommand\arraystretch{1.2}
\begin{tabular}{m{0.27\textwidth}<{\centering}|m{0.27\textwidth}<{\centering}|m{0.27\textwidth}<{\centering}}
    \hlineB{2.5}
    square pairs & $\pi_{i-1}>0\slash \pi_{i+1}>0$ & $\pi_{i-1}<0\slash\pi_{i+1}<0$\\ \hline
    $\langle i-1, |\pi_{i-1}| \rangle$, $\langle i,\pi_i\rangle$ &
    asc$^+_+$-pair $\rightarrow$ des$^+_-$-pair &
    asc$^-_+$-pair $\rightarrow$ des$^-_-$-pair\\ \hline
    $\langle i,\pi_i\rangle$, $\langle i+1, |\pi_{i+1}| \rangle$ &
    des$^+_+$-pair $\rightarrow$ asc$^+_-$-pair &
    des$^-_+$-pair $\rightarrow$ asc$^-_-$-pair\\
    \hlineB{2.5}
  \end{tabular}
\end{center}
After replacing the positive square $\langle i,\pi_i\rangle$ by a negative square,
we have $\des^B(\pi')=\des^B(\pi)$, and
$$\# {0_h^+}'=\# {0_h^+},\,\,\# {0_h^-}'=\# {0_h^-}, \,\,
\# {1_h^+}'=\# {1_h^+},\,\, \# {1_h^-}'=\# {1_h^-}.$$

\noindent{\bf Case 4}: $\bm{\pi_i<|\pi_{i-1}|}$ {\bf and} $\bm{\pi_i<|\pi_{i+1}|}$ {\bf for} $\bm{2\leq i\leq n-1}$.
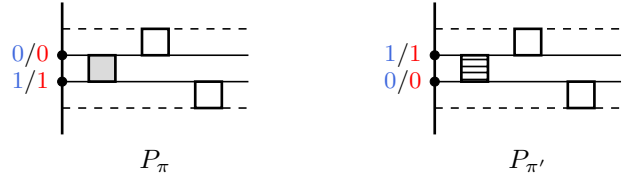
\begin{center}

\end{center}
\begin{center}
\renewcommand\arraystretch{1.2}
\begin{tabular}{m{0.27\textwidth}<{\centering}|m{0.27\textwidth}<{\centering}|m{0.27\textwidth}<{\centering}}
    \hlineB{2.5}
    square pairs & $\pi_{i-1}>0\slash \pi_{i+1}>0$ & $\pi_{i-1}<0\slash\pi_{i+1}<0$\\ \hline
    $\langle i-1, |\pi_{i-1}| \rangle$, $\langle i,\pi_i\rangle$ &
    des$^+_+$-pair $\rightarrow$ des$^-_+$-pair &
    asc$^+_-$-pair $\rightarrow$ asc$^-_-$-pair\\ \hline
    $\langle i,\pi_i\rangle$, $\langle i+1, |\pi_{i+1}| \rangle$ &
    asc$^+_+$-pair $\rightarrow$ asc$^-_+$-pair &
    des$^+_-$-pair $\rightarrow$ des$^-_-$-pair\\
    \hlineB{2.5}
  \end{tabular}
\end{center}
After replacing the positive square $\langle i,\pi_i\rangle$ by a negative square,
we have $\des^B(\pi')=\des^B(\pi)$, and
$$\# {0_h^+}'=\# {0_h^+},\,\,\# {0_h^-}'=\# {0_h^-}, \,\,
\# {1_h^+}'=\# {1_h^+},\,\, \# {1_h^-}'=\# {1_h^-}.$$

\noindent{\bf Case 5: $\bm{|\pi_{i+1}|<\pi_i}$ for $\bm{i=1}$.}
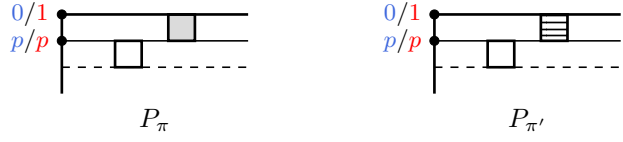
\begin{center}

\end{center}
\begin{center}
\renewcommand\arraystretch{1.2}
\begin{tabular}{m{0.27\textwidth}<{\centering}|m{0.27\textwidth}<{\centering}|m{0.27\textwidth}<{\centering}}
    \hlineB{2.5}
    square (pair) & $\pi_{2}>0$ & $\pi_{2}<0$\\ \hline
    $\langle 1,\pi_1\rangle$ &
    \multicolumn{2}{c}{asc$^+$-square $\rightarrow$ des$^-$-square} \\ \hline
    $\langle 1,\pi_1\rangle$, $\langle 2, |\pi_{2}| \rangle$ &
    des$^+_+$-pair $\rightarrow$ asc$^+_-$-pair &
    des$^-_+$-pair $\rightarrow$ asc$^-_-$-pair\\
    \hlineB{2.5}
  \end{tabular}
\end{center}
After replacing the positive square $\langle i,\pi_i\rangle$ by a negative square,
we have $\des^B(\pi')=\des^B(\pi)$, and
$$\# {0_h^+}'=\# {0_h^+},\,\,\# {0_h^-}'=\# {0_h^-}, \,\,
\# {1_h^+}'=\# {1_h^+},\,\, \# {1_h^-}'=\# {1_h^-}.$$

\noindent{\bf Case 6: $\bm{|\pi_i|<\pi_{i+1}}$ for $\bm{i=1}$.}
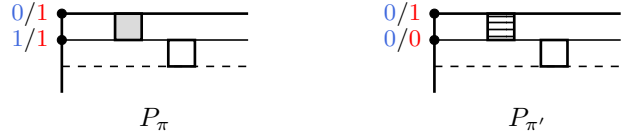
\begin{center}

\end{center}
\begin{center}
\renewcommand\arraystretch{1.2}
\begin{tabular}{m{0.27\textwidth}<{\centering}|m{0.27\textwidth}<{\centering}|m{0.27\textwidth}<{\centering}}
    \hlineB{2.5}
    square (pair) & $\pi_{2}>0$ & $\pi_{2}<0$\\ \hline
    $\langle 1,\pi_1\rangle$ &
    \multicolumn{2}{c}{asc$^+$-square $\rightarrow$ des$^-$-square} \\ \hline
    $\langle 1,\pi_1\rangle$, $\langle 2, |\pi_{2}| \rangle$ &
    asc$^+_+$-pair $\rightarrow$ asc$^-_+$-pair &
    des$^+_-$-pair $\rightarrow$ des$^-_-$-pair\\
    \hlineB{2.5}
  \end{tabular}
\end{center}
After replacing the positive square $\langle i,\pi_i\rangle$ by a negative square,
we have $\des^B(\pi')=\des^B(\pi)+1$, and
$$\# {0_h^+}'=\# {0_h^+}+1,\,\,\# {0_h^-}'=\# {0_h^-}+1, \,\,
\# {1_h^+}'=\# {1_h^+}-1,\,\, \# {1_h^-}'=\# {1_h^-}-1.$$

\noindent{\bf Case 7: $\bm{|\pi_{i-1}|<\pi_i}$ for $\bm{i=n}$.}
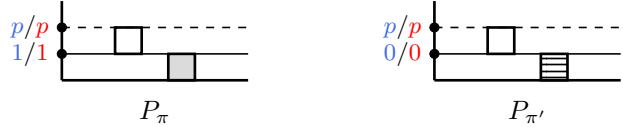
\begin{center}

\end{center}
\begin{center}
\renewcommand\arraystretch{1.2}
\begin{tabular}{m{0.27\textwidth}<{\centering}|m{0.27\textwidth}<{\centering}|m{0.27\textwidth}<{\centering}}
    \hlineB{2.5}
    square pair & $\pi_{n-1}>0$ & $\pi_{n-1}<0$\\ \hline
    $\langle n-1,|\pi_{n-1}|\rangle$, $\langle n, \pi_{n} \rangle$ &
    asc$^+_+$-pair $\rightarrow$ des$^+_-$-pair &
    asc$^-_+$-pair $\rightarrow$ des$^-_-$-pair\\
    \hlineB{2.5}
  \end{tabular}
\end{center}
After replacing the positive square $\langle i,\pi_i\rangle$ by a negative square,
we have $\des^B(\pi')=\des^B(\pi)+1$, and
$$\# {0_h^+}'=\# {0_h^+}+1,\,\,\# {0_h^-}'=\# {0_h^-}+1, \,\,
\# {1_h^+}'=\# {1_h^+}-1,\,\, \# {1_h^-}'=\# {1_h^-}-1.$$

\noindent{\bf Case 8: $\bm{\pi_i<|\pi_{i-1}|}$ for $\bm{i=n}$.}
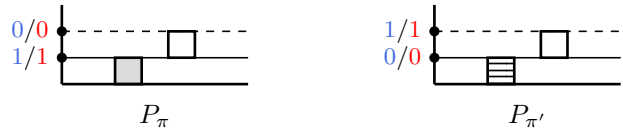
\begin{center}

\end{center}
\begin{center}
\renewcommand\arraystretch{1.2}
\begin{tabular}{m{0.27\textwidth}<{\centering}|m{0.27\textwidth}<{\centering}|m{0.27\textwidth}<{\centering}}
    \hlineB{2.5}
    square pair & $\pi_{n-1}>0$ & $\pi_{n-1}<0$\\ \hline
    $\langle n-1,|\pi_{n-1}|\rangle$, $\langle n, \pi_{n} \rangle$ &
    des$^+_+$-pair $\rightarrow$ des$^-_+$-pair &
    asc$^+_-$-pair $\rightarrow$ asc$^-_-$-pair\\
    \hlineB{2.5}
  \end{tabular}
\end{center}
After replacing the positive square $\langle i,\pi_i\rangle$ by a negative square,
we have $\des^B(\pi')=\des^B(\pi)$, and
$$\# {0_h^+}'=\# {0_h^+},\,\,\# {0_h^-}'=\# {0_h^-}, \,\,
\# {1_h^+}'=\# {1_h^+},\,\, \# {1_h^-}'=\# {1_h^-}.$$

\noindent{\bf Case 9: $\bm{|\pi_{j-1}^{-1}|<\pi_j^{-1}<|\pi_{j+1}^{-1}|}$ for $\bm{2\leq j\leq n-1}$.}
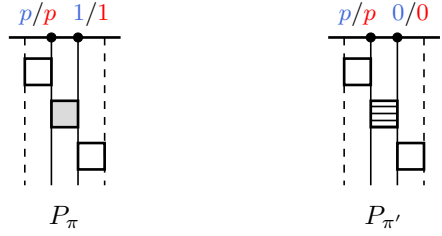
\begin{center}

\end{center}
\begin{center}
\renewcommand\arraystretch{1.2}
\begin{tabular}{m{0.27\textwidth}<{\centering}|m{0.28\textwidth}<{\centering}|m{0.28\textwidth}<{\centering}}
    \hlineB{2.5}
    square pairs & $\pi_{j-1}^{-1}>0\slash \pi_{j+1}^{-1}>0$ & $\pi_{j-1}^{-1}<0\slash\pi_{j+1}^{-1}<0$\\ \hline
    $\langle |\pi_{j-1}^{-1}|, j-1 \rangle$, $\langle \pi_j^{-1},j\rangle$ &
    iasc$^+_+$-pair $\rightarrow$ ides$^+_-$-pair &
    iasc$^-_+$-pair $\rightarrow$ ides$^-_-$-pair\\ \hline
    $\langle \pi_j^{-1},j\rangle$, $\langle |\pi_{j+1}^{-1}|, j+1 \rangle$ &
    iasc$^+_+$-pair $\rightarrow$ iasc$^-_+$-pair &
    ides$^+_-$-pair $\rightarrow$ ides$^-_-$-pair\\
    \hlineB{2.5}
  \end{tabular}
\end{center}
After replacing the positive square $\langle \pi_j^{-1},j\rangle$ by a negative square,
we have $\ides^B(\pi')=\ides^B(\pi)+1$, and
$$\# {0_v^+}'=\# {0_v^+}+1,\,\,\# {0_v^-}'=\# {0_v^-}+1, \,\,
\# {1_v^+}'=\# {1_v^+}-1,\,\, \# {1_v^-}'=\# {1_v^-}-1.$$

\noindent{\bf Case 10: $\bm{|\pi_{j+1}^{-1}|<\pi_j^{-1}<|\pi_{j-1}^{-1}|}$ for $\bm{2\leq j\leq n-1}$.}
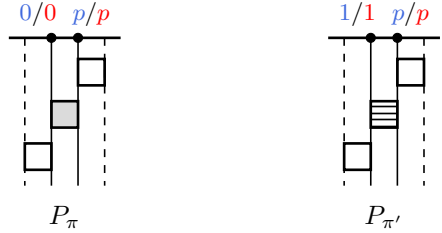
\begin{center}

\end{center}
\begin{center}
\renewcommand\arraystretch{1.2}
\begin{tabular}{m{0.27\textwidth}<{\centering}|m{0.28\textwidth}<{\centering}|m{0.28\textwidth}<{\centering}}
    \hlineB{2.5}
    square pairs & $\pi_{j-1}^{-1}>0\slash \pi_{j+1}^{-1}>0$ & $\pi_{j-1}^{-1}<0\slash\pi_{j+1}^{-1}<0$\\ \hline
    $\langle |\pi_{j-1}^{-1}|, j-1 \rangle$, $\langle \pi_j^{-1},j\rangle$ &
    ides$^+_+$-pair $\rightarrow$ ides$^-_+$-pair &
    iasc$^+_-$-pair $\rightarrow$ iasc$^-_-$-pair\\ \hline
    $\langle \pi_j^{-1},j\rangle$, $\langle |\pi_{j+1}^{-1}|, j+1 \rangle$ &
    ides$^+_+$-pair $\rightarrow$ iasc$^+_-$-pair &
    ides$^-_+$-pair $\rightarrow$ iasc$^-_-$-pair\\
    \hlineB{2.5}
  \end{tabular}
\end{center}
After replacing the positive square $\langle \pi_j^{-1},j\rangle$ by a negative square,
we have $\ides^B(\pi')=\ides^B(\pi)-1$, and
$$\# {0_v^+}'=\# {0_v^+}-1,\,\,\# {0_v^-}'=\# {0_v^-}-1, \,\,
\# {1_v^+}'=\# {1_v^+}+1,\,\, \# {1_v^-}'=\# {1_v^-}+1.$$

\noindent{\bf Case 11: $\bm{|\pi_{j+1}^{-1}|<\pi_j^{-1}}$ and $\bm{|\pi_{j-1}^{-1}|<\pi_j^{-1}}$ for $\bm{2\leq j\leq n-1}$.}
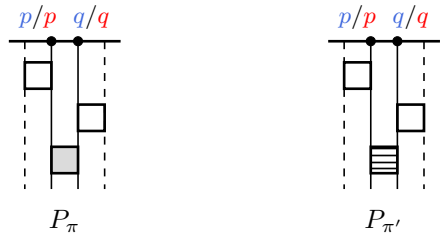
\begin{center}

\end{center}
\begin{center}
\renewcommand\arraystretch{1.2}
\begin{tabular}{m{0.27\textwidth}<{\centering}|m{0.28\textwidth}<{\centering}|m{0.28\textwidth}<{\centering}}
    \hlineB{2.5}
    square pairs & $\pi_{j-1}^{-1}>0\slash \pi_{j+1}^{-1}>0$ & $\pi_{j-1}^{-1}<0\slash\pi_{j+1}^{-1}<0$\\ \hline
    $\langle |\pi_{j-1}^{-1}|, j-1 \rangle$, $\langle \pi_j^{-1},j\rangle$ &
    iasc$^+_+$-pair $\rightarrow$ ides$^+_-$-pair &
    iasc$^-_+$-pair $\rightarrow$ ides$^-_-$-pair\\ \hline
    $\langle \pi_j^{-1},j\rangle$, $\langle |\pi_{j+1}^{-1}|, j+1 \rangle$ &
    ides$^+_+$-pair $\rightarrow$ iasc$^+_-$-pair &
    ides$^-_+$-pair $\rightarrow$ iasc$^-_-$-pair\\
    \hlineB{2.5}
  \end{tabular}
\end{center}
After replacing the positive square $\langle \pi_j^{-1},j\rangle$ by a negative square,
we have $\ides^B(\pi')=\ides^B(\pi)$, and
$$\# {0_v^+}'=\# {0_v^+},\,\,\# {0_v^-}'=\# {0_v^-}, \,\,
\# {1_v^+}'=\# {1_v^+},\,\, \# {1_v^-}'=\# {1_v^-}.$$

\noindent{\bf Case 12: $\bm{\pi_j^{-1}<|\pi_{j+1}^{-1}|}$ and $\bm{\pi_j^{-1}<|\pi_{j-1}^{-1}|}$ for $\bm{2\leq j\leq n-1}$.}
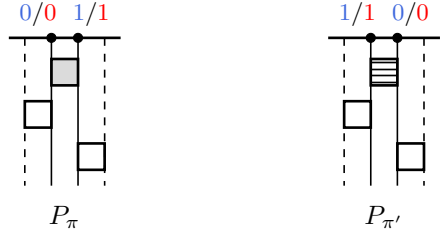
\begin{center}

\end{center}
\begin{center}
\renewcommand\arraystretch{1.2}
\begin{tabular}{m{0.27\textwidth}<{\centering}|m{0.28\textwidth}<{\centering}|m{0.28\textwidth}<{\centering}}
    \hlineB{2.5}
    square pairs & $\pi_{j-1}^{-1}>0\slash \pi_{j+1}^{-1}>0$ & $\pi_{j-1}^{-1}<0\slash\pi_{j+1}^{-1}<0$\\ \hline
    $\langle |\pi_{j-1}^{-1}|, j-1 \rangle$, $\langle \pi_j^{-1},j\rangle$ &
    ides$^+_+$-pair $\rightarrow$ ides$^-_+$-pair &
    iasc$^+_-$-pair $\rightarrow$ iasc$^-_-$-pair\\ \hline
    $\langle \pi_j^{-1},j\rangle$, $\langle |\pi_{j+1}^{-1}|, j+1 \rangle$ &
    iasc$^+_+$-pair $\rightarrow$ iasc$^-_+$-pair &
    ides$^+_-$-pair $\rightarrow$ ides$^-_-$-pair\\
    \hlineB{2.5}
  \end{tabular}
\end{center}
After replacing the positive square $\langle \pi_j^{-1},j\rangle$ by a negative square,
we have $\ides^B(\pi')=\ides^B(\pi)$, and
$$\# {0_v^+}'=\# {0_v^+},\,\,\# {0_v^-}'=\# {0_v^-}, \,\,
\# {1_v^+}'=\# {1_v^+},\,\, \# {1_v^-}'=\# {1_v^-}.$$

\noindent{\bf Case 13: $\bm{\pi_j^{-1}<|\pi_{j+1}^{-1}|}$ for $\bm{j=1}$.}
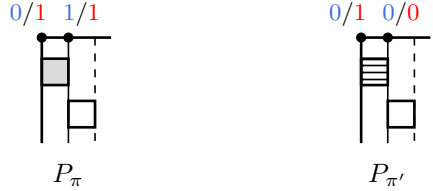
\begin{center}

\end{center}
\begin{center}
\renewcommand\arraystretch{1.2}
\begin{tabular}{m{0.27\textwidth}<{\centering}|m{0.28\textwidth}<{\centering}|m{0.28\textwidth}<{\centering}}
    \hlineB{2.5}
    square (pair) & $\pi_{2}^{-1}>0$ & $\pi_{2}^{-1}<0$\\ \hline
    $\langle \pi_1^{-1},1\rangle$ &
    \multicolumn{2}{c}{iasc$^+$-square $\rightarrow$ ides$^-$-square} \\ \hline
    $\langle \pi_1^{-1},1\rangle$, $\langle |\pi_2^{-1}|,2 \rangle$ &
    iasc$^+_+$-pair $\rightarrow$ iasc$^-_+$-pair &
    ides$^+_-$-pair $\rightarrow$ ides$^-_-$-pair\\
    \hlineB{2.5}
  \end{tabular}
\end{center}
After replacing the positive square $\langle \pi_j^{-1},j\rangle$ by a negative square,
we have $\ides^B(\pi')=\ides^B(\pi)+1$, and
$$\# {0_v^+}'=\# {0_v^+}+1,\,\,\# {0_v^-}'=\# {0_v^-}+1, \,\,
\# {1_v^+}'=\# {1_v^+}-1,\,\, \# {1_v^-}'=\# {1_v^-}-1.$$

\noindent{\bf Case 14: $\bm{|\pi_{j+1}^{-1}|<\pi_j^{-1}}$ for $\bm{j=1}$.}
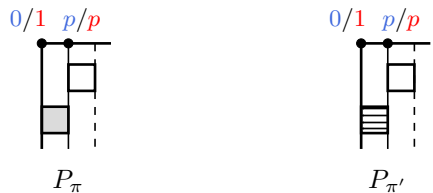
\begin{center}

\end{center}
\begin{center}
\renewcommand\arraystretch{1.2}
\begin{tabular}{m{0.27\textwidth}<{\centering}|m{0.28\textwidth}<{\centering}|m{0.28\textwidth}<{\centering}}
    \hlineB{2.5}
    square (pair) & $\pi_{2}^{-1}>0$ & $\pi_{2}^{-1}<0$\\ \hline
    $\langle \pi_1^{-1},1\rangle$ &
    \multicolumn{2}{c}{iasc$^+$-square $\rightarrow$ ides$^-$-square} \\ \hline
    $\langle \pi_1^{-1},1\rangle$, $\langle |\pi_2^{-1}|,2 \rangle$ &
    ides$^+_+$-pair $\rightarrow$ iasc$^+_-$-pair &
    ides$^-_+$-pair $\rightarrow$ iasc$^-_-$-pair\\
    \hlineB{2.5}
  \end{tabular}
\end{center}
After replacing the positive square $\langle \pi_j^{-1},j\rangle$ by a negative square,
we have $\ides^B(\pi')=\ides^B(\pi)$, and
$$\# {0_v^+}'=\# {0_v^+},\,\,\# {0_v^-}'=\# {0_v^-}, \,\,
\# {1_v^+}'=\# {1_v^+},\,\, \# {1_v^-}'=\# {1_v^-}.$$

\noindent{\bf Case 15: $\bm{\pi_j^{-1}<|\pi_{j-1}^{-1}|}$ for $\bm{j=n}$.}
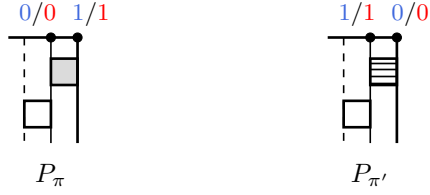
\begin{center}

\end{center}
\begin{center}
\renewcommand\arraystretch{1.2}
\begin{tabular}{m{0.27\textwidth}<{\centering}|m{0.28\textwidth}<{\centering}|m{0.28\textwidth}<{\centering}}
    \hlineB{2.5}
    square pair & $\pi_{n-1}^{-1}>0$ & $\pi_{n-1}^{-1}<0$\\ \hline
    $\langle |\pi_{n-1}^{-1}|,n-1 \rangle$, $\langle \pi_n^{-1},n\rangle$ &
    ides$^+_+$-pair $\rightarrow$ ides$^-_+$-pair &
    iasc$^+_-$-pair $\rightarrow$ iasc$^-_-$-pair\\
    \hlineB{2.5}
  \end{tabular}
\end{center}
After replacing the positive square $\langle \pi_j^{-1},j\rangle$ by a negative square,
we have $\ides^B(\pi')=\ides^B(\pi)$, and
$$\# {0_v^+}'=\# {0_v^+},\,\,\# {0_v^-}'=\# {0_v^-}, \,\,
\# {1_v^+}'=\# {1_v^+},\,\, \# {1_v^-}'=\# {1_v^-}.$$

\noindent{\bf Case 16: $\bm{|\pi_{j-1}^{-1}|<\pi_j^{-1}}$ for $\bm{j=n}$.}
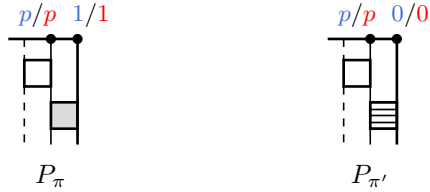
\begin{center}

\end{center}
\begin{center}
\renewcommand\arraystretch{1.2}
\begin{tabular}{m{0.27\textwidth}<{\centering}|m{0.28\textwidth}<{\centering}|m{0.28\textwidth}<{\centering}}
    \hlineB{2.5}
    square pair & $\pi_{n-1}^{-1}>0$ & $\pi_{n-1}^{-1}<0$\\ \hline
    $\langle |\pi_{n-1}^{-1}|,n-1 \rangle$, $\langle \pi_n^{-1},n\rangle$ &
    iasc$^+_+$-pair $\rightarrow$ ides$^+_-$-pair &
    iasc$^-_+$-pair $\rightarrow$ ides$^-_-$-pair\\
    \hlineB{2.5}
  \end{tabular}
\end{center}
After replacing the positive square $\langle \pi_j^{-1},j\rangle$ by a negative square,
we have $\ides^B(\pi')=\ides^B(\pi)+1$, and
$$\# {0_v^+}'=\# {0_v^+}+1,\,\,\# {0_v^-}'=\# {0_v^-}+1, \,\,
\# {1_v^+}'=\# {1_v^+}-1,\,\, \# {1_v^-}'=\# {1_v^-}-1.$$

Therefore, by induction, we complete the proof.
\qed

\subsection{Enumerations of grid points of certain $d$-types}
For  $\pi\in\mathfrak{B}_n$,
let $n(\pi)$ be the number of negative elements in $\pi_1\pi_2\cdots\pi_n$.
With the help of Theorems \ref{lempath} and \ref{thmpaths},
we obtain the number of grid points of  given $d^+$-type or $d^-$-type in the grid $P_\pi$ by $\des^B(\pi)$, $\ides^B(\pi)$ and $n(\pi)$.

\begin{thm}\label{thm-enum-dtype}
Let $\pi\in\B_n$, then in the grid $P_\pi$, there are
\begin{enumerate}[(a)]
    \item $(\des^B(\pi)+1)(\ides^B(\pi)+1)-n(\pi)+n$ grid points of $d^+$-type$ (0,0)$, \label{aaaa}
    \item $(\ides^B(\pi)+1)(n-\des^B(\pi))+n(\pi)-n$ grid points of $d^+$-type $(1,0)$,\label{bbbb}
    \item $(\des^B(\pi)+1)(n-\ides^B(\pi))+n(\pi)-n$ grid points of $d^+$-type $(0,1)$,\label{cccc}
    \item $(n-\des^B(\pi))(n-\ides^B(\pi))-n(\pi)+n$ grid points of $d^+$-type $(1,1)$;\label{dddd}
\end{enumerate}
and
\begin{enumerate}[(a)]\addtocounter{enumi}{4}
    \item $\des^B(\pi)\ides^B(\pi)+n(\pi)$ grid points of $d^-$-type$ (0,0)$, \label{eeee}
    \item $\ides^B(\pi)(n-\des^B(\pi)+1)-n(\pi)$ grid points of $d^-$-type $(1,0)$,\label{ffff}
    \item $\des^B(\pi)(n-\ides^B(\pi)+1)-n(\pi)$ grid points of $d^-$-type $(0,1)$,\label{gggg}
    \item $(n-\des^B(\pi)+1)(n-\ides^B(\pi)+1)+n(\pi)$ grid points of $d^-$-type $(1,1)$.\label{hhhh}
\end{enumerate}
\end{thm}

\pf
By Definition \ref{def-path}, for $p,q\in\{0,1\}$, a grid point is of $d^+$-type (resp. $d^-$-type) $(p,q)$ if and only if
it is an intersection of a $p_h^+$-path (resp. $p_h^-$-path) and a $q_v^+$-path (resp. $q_v^-$-path). Note that the numbers of positive and negative squares in the grid $P_\pi$ are $n-n(\pi)$ and $n(\pi)$, respectively.
Referring to Example \ref{examdtype}, we will have a clearer picture of the  proof that follows.

Since there are $(\des^B(\pi)+1)$ $0_h^+$-paths and $(\ides^B(\pi)+1)$ $0_v^+$-paths in $P_\pi$ by Theorem \ref{thmpaths},
the $0_h^+$-paths and $0_v^+$-paths meet each other $(\des^B(\pi)+1)(\ides^B(\pi)+1)$ times.
By noticing that both $0_h^+$-paths and $0_v^+$-paths take a southeast step when they encounter a positive square by Theorem \ref{lempath},
we see that the number of intersections formed by $0_h^+$-paths and $0_v^+$-paths at grid points is actually $n-n(\pi)$, the number of positive squares, more than the number of times they meet.
Thus we prove the statement \eqref{aaaa}.
In terms of the counting formulas in Theorem \ref{thmpaths},
and similar analyses but replacing $0_h^+$-paths with $1_h^+$-paths, $0_h^-$-paths and $1_h^-$-paths respectively,
and replacing $0_v^+$-paths with $1_v^+$-paths, $0_v^-$-paths and $1_v^-$-paths correspondingly,
we can prove statements \eqref{dddd}, \eqref{eeee} and \eqref{hhhh}.

To see the number of grid points of $d^+$-type $(1,0)$,
we notice that if a $1_h^+$-path and a $0_v^+$-path meet in a positive square,
their intersection is in the interior of the square instead of on the grid point,
because by Theorem \ref{lempath}, the $1_h^+$-path and the $0_v^+$-path take the northeast and southeast steps, respectively, inside the square.
Thus, the number of intersections on grid points formed by $1_h^+$-paths and $0_v^+$-paths is $n-n(\pi)$ less than the number of times they meet.
By Theorem \ref{thmpaths}, such intersections are counted by $(\ides^B(\pi)+1)(n-\des^B(\pi))+n(\pi)-n$,
which verifies statement \eqref{bbbb}.
By replacing the $1_h^+$-paths in the above proof process with $0_h^+$-paths, $1_h^-$-paths, $0_h^-$-paths,
and the $0_v^+$-paths with $1_v^+$-paths, $0_v^-$-paths, $1_v^-$-paths, respectively,
and using the related quantities in Theorem \ref{thmpaths},
we can  prove statements \eqref{cccc}, \eqref{ffff}, \eqref{gggg}, respectively.
\qed

\begin{exam}\label{examdtype}
Let $\pi=2\bar{4}3\bar{1}5\in\B_5$ with $\des^B(\pi)=2$, $\ides^B(\pi)=2$ and $n(\pi)=2$,
the following diagrams present the $d^+$-type and $d^-$-type of each grid point in the grid $P_\pi$ by Theorem \ref{thm-enum-dtype},
where $\#d^\ast(p,q)$ denote the number of grid points of $d^\ast$-type $(p,q)$
for $\ast\in\{+,-\}$ and $p,q\in\{0,1\}$.
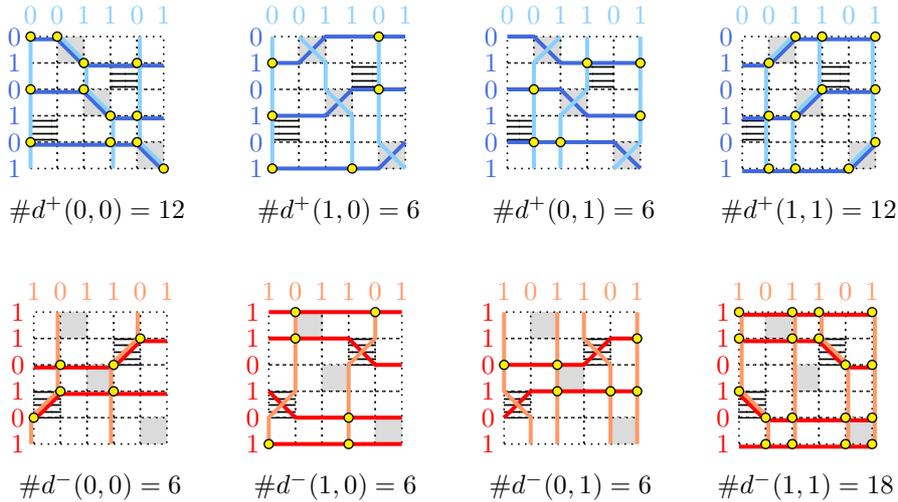
\begin{center}

\end{center}
\end{exam}

\section{The recurrence of $\ob$-descents and $\ob$-inverse descents on $\B_n$}\label{Sec-recBnij}

To give a combinatorial proof of the four-term recurrence \eqref{rec^obnij}
 of $\overline{b}_{n,i,j}$ given in Theorem \ref{thmof^obnij},
we first introduce some sets of ordered pairs.
Let
$$
\mathcal{C}_{n,i,j}=\{(\sigma, k) \mid \pi\in\overline{\mathcal{B}}_{n, i, j} \text{ and }1\leq k\leq n \},
$$
and recall that
$$
\overline{\mathcal{B}}_{n,i,j}=\{\sigma\in\B_n\mid \mbox{ $\des^B(\sigma)=i$ and $\ides^B(\sigma)=j$}\}.
$$
Thus we have $|\mathcal{C}_{n,i,j}|=n\overline{b}_{n,i,j}$,
which is equal to the left side of \eqref{rec^obnij}.
Let
\begin{equation}\label{defDpq}
\mathcal{D}_{n,i,j}^{(p,q)^\ast}=\left\{(\pi,(r,s))\mid \text{ $\pi\in\overline{\mathcal{B}}_{n,i,j}$ and $d^\ast(r,s)=(p,q)$ in $P_\pi$}\right\}
\end{equation}
with $\ast\in\{+,-\}$ and $p,q\in\{0,1\}$.
Clearly, for any given $1\leq i, j\leq n$,
the eight sets $\mathcal{D}_{n,i,j}^{(0,0)^+}$, $\mathcal{D}_{n,i,j}^{(0,0)^-}$,
$\mathcal{D}_{n,i,j-1}^{(0,1)^+}$, $\mathcal{D}_{n,i,j-1}^{(0,1)^-}$,
$\mathcal{D}_{n,i-1,j}^{(1,0)^+}$, $\mathcal{D}_{n,i-1,j}^{(1,0)^-}$,
$\mathcal{D}_{n,i-1,j-1}^{(1,1)^+}$ and $\mathcal{D}_{n,i-1,j-1}^{(1,1)^-}$
are pairwise disjoint,
and by Theorem \ref{thm-enum-dtype}, we see
\begin{equation*}
\begin{aligned}
&\left|\mathcal{D}_{n-1,i,j}^{(0,0)^+}\uplus\mathcal{D}_{n-1,i,j}^{(0,0)^-}\right|
=(n+i+j+2ij) \overline{b}_{n-1, i, j},\\[3pt]
&\left|\mathcal{D}_{n-1,i,j-1}^{(0,1)^+}\uplus\mathcal{D}_{n-1,i,j-1}^{(0,1)^-}\right|
=((2n-1)i-(2i+1)(j-1)) \overline{b}_{n-1, i, j-1},\\[3pt]
&\left|\mathcal{D}_{n-1,i-1,j}^{(1,0)^+}\uplus\mathcal{D}_{n-1,i-1,j}^{(1,0)^-}\right|
=((2n-1)j-(2j+1)(i-1)) \overline{b}_{n-1, i-1, j},\\[3pt]
&\left|\mathcal{D}_{n-1,i-1,j-1}^{(1,1)^+}\uplus\mathcal{D}_{n-1,i-1,j-1}^{(1,1)^-}\right|
=((2n^2-n)+2(i-1)(j-1)+(1-2n)(i+j-2)) \overline{b}_{n-1, i-1, j-1},
\end{aligned}
\end{equation*}
where $\biguplus$ indicates the union of pairwise disjoint sets.
Hence the cardinality of the set $\mathcal{D}_{n-1,i,j}$ defined as
\begin{equation*}
  \mathcal{D}_{n-1,i,j}:=\biguplus_{p,q\in\{0,1\}}\mathcal{D}_{n-1,i-p,j-q}^{(p,q)^+}\uplus\mathcal{D}_{n-1,i-p,j-q}^{(p,q)^-}
\end{equation*}
is equal to the right side of \eqref{rec^obnij}.

{\noindent{\emph{Combinatorial Proof of Theorem \ref{thmof^obnij}.}\hskip 2pt}}
We shall establish a bijection between $\mathcal{D}_{n-1,i,j}$ and $\mathcal{C}_{n,i,j}$ by employing the inserting operation $\varphi_{(i,j)}$ or $\overline{\varphi}_{(i,j)}$.
In particular, define
\begin{equation*}
\Psi((\pi,(r,s)))=(\sigma, r)
\end{equation*}
for any pair $(\pi,(r,s))\in{\mathcal{D}}_{n-1,i,j}$, where
\begin{equation}\label{defpsi}
\sigma=\left\{
      \begin{array}{ll}
        \varphi_{(r,s)}(\pi), & \hbox{if $(\pi,(r,s))\in\mathcal{D}_{n-1,i-p,j-q}^{(p,q)^+}$ ,} \\[6pt]
        \overline{\varphi}_{(r,s)}(\pi), & \hbox{if $(\pi,(r,s))\in\mathcal{D}_{n-1,i-p,j-q}^{(p,q)^-}$ ,}
      \end{array}
    \right.
\end{equation}
for $p,q\in\{0,1\}$.

By the construction as given in \eqref{defDpq},
for
$$((\pi,(r,s)))\in\mathcal{D}_{n-1,i-p,j-q}^{(p,q)^+}
\uplus\mathcal{D}_{n-1,i-p,j-q}^{(p,q)^-} ,$$
we have $\pi\in\overline{\mathcal{B}}_{n-1,i-p,j-q}$, and
the grid point $(r,s)$ is of $d^+$-type or $d^-$-type $(p,q)$,
which implies
\begin{equation*}
(\des(\varphi_{(r,s)}(\pi)),\ides(\varphi_{(r,s)}(\pi)))-(\des(\pi),\ides(\pi))=(p,q),
\end{equation*}
or
\begin{equation*}
(\des(\overline{\varphi}_{(r,s)}(\pi)),\ides(\overline{\varphi}_{(r,s)}(\pi)))-(\des(\pi),\ides(\pi))=(p,q)
\end{equation*}
by Definitions \ref{defd+type} and \ref{defd-type} of $d$-types.
Thus we obtain $\sigma\in \mathcal{B}_{n,i,j}$,
and $(\sigma,r)\in\mathcal{C}_{n,i,j}$.

On the other hand,
for any pair $(\sigma ,r)\in\mathcal{C}_{n,i,j}$,
notice that the signed permutation $\sigma$ has $i$ descents and $j$ idescents.
Then deleting the square $\langle r, \sigma_r\rangle$ in the grid $P_\sigma$
by the operation $\varphi^{-1}_{\langle r, \sigma_r\rangle}$ if $\sigma_r>0$ and $\overline{\varphi}^{-1}_{\langle r, \sigma_r\rangle}$ if $\sigma_r<0$,
we obtain a signed permutation $\pi$ that must belong the set $\B_{n,i-p,j-p}$ for some
$p,q\in\{0,1\}$ by Proposition \ref{prop-possible}.
Moreover,
the grid point $(r,\sigma_r)$ in the grid $P_\pi$ must have $d^+$-type $(p,q)$
if $\sigma_r>0$ or $d^-$-type $(p,q)$ if $\sigma_r<0$.
Thus, we deduce
$$(\pi,(r,\sigma_r))\in\mathcal{D}_{n-1,i-p,j-q}^{(p,q)^+}
\uplus\mathcal{D}_{n-1,i-p,j-q}^{(p,q)^-} \, ,$$
and the inverse of $\Psi$ is given as
\begin{equation*}
\Psi^{-1}((\sigma, r))=(\pi,(r,\sigma_r))
\end{equation*}
for any pair $(\sigma, r)\in{\mathcal{C}}_{n,i,j}$, where
\begin{equation}\label{defpsi^-1}
\pi=\left\{
      \begin{array}{ll}
        \varphi^{-1}_{\langle r, \sigma_r\rangle}(\sigma), & \hbox{if $\sigma_r>0$,} \\[6pt]
        \overline{\varphi}^{-1}_{\langle r, \sigma_r\rangle}(\sigma), & \hbox{if $\sigma_r<0$.}
      \end{array}
    \right.
\end{equation}

Therefore, with the combination of \eqref{defpsi} and \eqref{defpsi^-1},
we get the desired bijection
\[
\Psi:\,{\mathcal{D}}_{n-1,i,j}\,\leftrightarrow\,\mathcal{C}_{n,i,j},
\]
which completes the proof of Theorem \ref{thmof^obnij}.
\qed

\begin{exam}
For $\sigma=3\overline{2}\overline{5}1\overline{4}\in\overline{B}_{5,3,2}$,
as shown below, for $1\leq r\leq 5$, we have
\begin{equation*}
\begin{aligned}
&\Psi^{-1}(\sigma,1)=(\overline{24}1\overline{3},(1,3))\in\mathcal{D}_{4,3,2}^{(0,0)^+},\quad
\Psi^{-1}(\sigma,2)=(2\overline{4}1\overline{3},(2,2))\in\mathcal{D}_{4,2,2}^{(1,0)^-},\\[3pt]
&\Psi^{-1}(\sigma,3)=(3\overline{2}1\overline{4},(3,5))\in\mathcal{D}_{4,2,2}^{(1,0)^-},\quad
\Psi^{-1}(\sigma,4)=(2\overline{143},(4,1))\in\mathcal{D}_{4,2,1}^{(1,1)^+},\\[3pt]
&\Psi^{-1}(\sigma,5)=(2\overline{34}1,(5,4))\in\mathcal{D}_{4,2,2}^{(1,0)^-}.
\end{aligned}
\end{equation*}

\begin{center}
\begin{tikzpicture}[scale = 0.6]
\def\hezi{-- +(5mm,0mm) -- +(5mm,5mm) -- +(0mm,5mm) -- cycle [line width=0.6pt]}
\def\judy{-- +(5mm,0mm) -- +(5mm,5mm) -- +(0mm,5mm) -- cycle [line width=0.6pt,fill=gainsboro]}
\def\nbx{-- +(5mm,0mm) -- +(5mm,5mm) -- +(0mm,5mm) -- cycle [line width=0.6pt,pattern={Lines[angle=45,distance=2.5pt, line width=0.6pt]}, pattern color=black,dotted]}
\tikzstyle{rdot}=[circle,fill=red,draw=red,inner sep=1.8]
\tikzstyle{cc}=[circle,draw=black,fill=yellow, line width=0.5pt, inner sep=1.5]

\begin{scope}[yshift=-3mm]
\draw[line width=0.6pt] (-5mm,34mm)--(185mm,34mm);
\node at (12.5mm, 40mm) {$\sigma\in\mathcal{B}_{5,3,2}$};
\node at (50mm, 40mm) {$r=1$};
\node at (80mm, 40mm) {$r=2$};
\node at (110mm, 40mm) {$r=3$};
\node at (140mm, 40mm) {$r=4$};
\node at (170mm, 40mm) {$r=5$};
\draw[line width=1pt] (-5mm,46mm)--(185mm,46mm);
\draw[line width=1pt] (-5mm,-11mm)--(185mm,-11mm);
\draw[line width=0.6pt] (32mm,46mm)--(32mm,-11mm);
\end{scope}

\draw[line width=0.8pt] (0mm,0mm)--(0mm,25mm);
\draw[line width=0.8pt] (5mm,0mm)--(5mm,25mm);
\draw[line width=0.8pt] (10mm,0mm)--(10mm,25mm);
\draw[line width=0.8pt] (15mm,0mm)--(15mm,25mm);
\draw[line width=0.8pt] (20mm,0mm)--(20mm,25mm);
\draw[line width=0.8pt] (25mm,0mm)--(25mm,25mm);

\draw[line width=0.8pt] (0mm,0mm)--(25mm,0mm);
\draw[line width=0.8pt] (0mm,5mm)--(25mm,5mm);
\draw[line width=0.8pt] (0mm,10mm)--(25mm,10mm);
\draw[line width=0.8pt] (0mm,15mm)--(25mm,15mm);
\draw[line width=0.8pt] (0mm,20mm)--(25mm,20mm);
\draw[line width=0.8pt] (0mm,25mm)--(25mm,25mm);

\draw (10mm,20mm)\judy;\draw (5mm,15mm)\nbx;\draw (20mm,10mm)\nbx;
\draw (0mm,5mm)\judy;\draw (15mm,0mm)\nbx;

\node at (12.5mm, -6mm) {$3\overline{2}\overline{5}1\overline{4}$};

\begin{scope}[shift={(40mm,0mm)}]
\draw[line width=0.8pt] (0mm,0mm)--(0mm,20mm);
\draw[line width=0.8pt] (5mm,0mm)--(5mm,20mm);
\draw[line width=0.8pt] (10mm,0mm)--(10mm,20mm);
\draw[line width=0.8pt] (15mm,0mm)--(15mm,20mm);
\draw[line width=0.8pt] (20mm,0mm)--(20mm,20mm);

\draw[line width=0.8pt] (0mm,0mm)--(20mm,0mm);
\draw[line width=0.8pt] (0mm,5mm)--(20mm,5mm);
\draw[line width=0.8pt] (0mm,10mm)--(20mm,10mm);
\draw[line width=0.8pt] (0mm,15mm)--(20mm,15mm);
\draw[line width=0.8pt] (0mm,20mm)--(20mm,20mm);

\draw (5mm,15mm)\nbx;\draw (15mm,10mm)\nbx;
\draw (0mm,5mm)\judy;\draw (10mm,0mm)\nbx;

\node at (-2mm,20mm){\royalblue{$0$}};
\node at (10mm,24mm){\babyblue{$0$}};

\draw (0mm,20mm)--(20mm,20mm)[line width=1.5pt,draw=royalblue];
\draw (10mm,20mm)--(10mm,0mm)[line width=1.5pt,draw=babyblue];

\node[cc] at (10mm,20mm){};

\node at (10mm, -6mm) {$\overline{24}1\overline{3}$};

\end{scope}

\begin{scope}[shift={(70mm,0mm)}]
\draw[line width=0.8pt] (0mm,0mm)--(0mm,20mm);
\draw[line width=0.8pt] (5mm,0mm)--(5mm,20mm);
\draw[line width=0.8pt] (10mm,0mm)--(10mm,20mm);
\draw[line width=0.8pt] (15mm,0mm)--(15mm,20mm);
\draw[line width=0.8pt] (20mm,0mm)--(20mm,20mm);

\draw[line width=0.8pt] (0mm,0mm)--(20mm,0mm);
\draw[line width=0.8pt] (0mm,5mm)--(20mm,5mm);
\draw[line width=0.8pt] (0mm,10mm)--(20mm,10mm);
\draw[line width=0.8pt] (0mm,15mm)--(20mm,15mm);
\draw[line width=0.8pt] (0mm,20mm)--(20mm,20mm);

\draw (5mm,15mm)\judy;\draw (15mm,10mm)\nbx;
\draw (0mm,5mm)\judy;\draw (10mm,0mm)\nbx;

\node at (-2mm,15mm){\red{$1$}};
\node at (5mm,24mm){\orange{$0$}};

\draw (0mm,15mm)--(15mm,15mm)--(20mm,10mm)[line width=1.5pt,draw=red];
\draw (5mm,20mm)--(5mm,0mm)[line width=1.5pt,draw=orange];

\node[cc] at (5mm,15mm){};

\node at (10mm, -6mm) {$2\overline{4}1\overline{3}$};
\end{scope}

\begin{scope}[shift={(100mm,0mm)}]
\draw[line width=0.8pt] (0mm,0mm)--(0mm,20mm);
\draw[line width=0.8pt] (5mm,0mm)--(5mm,20mm);
\draw[line width=0.8pt] (10mm,0mm)--(10mm,20mm);
\draw[line width=0.8pt] (15mm,0mm)--(15mm,20mm);
\draw[line width=0.8pt] (20mm,0mm)--(20mm,20mm);

\draw[line width=0.8pt] (0mm,0mm)--(20mm,0mm);
\draw[line width=0.8pt] (0mm,5mm)--(20mm,5mm);
\draw[line width=0.8pt] (0mm,10mm)--(20mm,10mm);
\draw[line width=0.8pt] (0mm,15mm)--(20mm,15mm);
\draw[line width=0.8pt] (0mm,20mm)--(20mm,20mm);

\draw (10mm,15mm)\judy;\draw (5mm,10mm)\nbx;
\draw (0mm,5mm)\judy;\draw (15mm,0mm)\nbx;

\node at (-2mm,15mm){\red{$1$}};
\node at (20mm,24mm){\orange{$0$}};

\draw (0mm,15mm)--(5mm,15mm)--(10mm,10mm)--(20mm,10mm)[line width=1.5pt,draw=red];
\draw (20mm,20mm)--(20mm,5mm)--(15mm,0mm)[line width=1.5pt,draw=orange];

\node[cc] at (20mm,10mm){};

\node at (10mm, -6mm) {$3\overline{2}1\overline{4}$};
\end{scope}

\begin{scope}[shift={(130mm,0mm)}]
\draw[line width=0.8pt] (0mm,0mm)--(0mm,20mm);
\draw[line width=0.8pt] (5mm,0mm)--(5mm,20mm);
\draw[line width=0.8pt] (10mm,0mm)--(10mm,20mm);
\draw[line width=0.8pt] (15mm,0mm)--(15mm,20mm);
\draw[line width=0.8pt] (20mm,0mm)--(20mm,20mm);

\draw[line width=0.8pt] (0mm,0mm)--(20mm,0mm);
\draw[line width=0.8pt] (0mm,5mm)--(20mm,5mm);
\draw[line width=0.8pt] (0mm,10mm)--(20mm,10mm);
\draw[line width=0.8pt] (0mm,15mm)--(20mm,15mm);
\draw[line width=0.8pt] (0mm,20mm)--(20mm,20mm);

\draw (5mm,15mm)\judy;\draw (0mm,10mm)\nbx;
\draw (15mm,5mm)\nbx;\draw (10mm,0mm)\nbx;

\node at (-2mm,5mm){\royalblue{$1$}};
\node at (0mm,24mm){\babyblue{$1$}};

\draw (0mm,5mm)--(20mm,5mm)[line width=1.5pt,draw=royalblue];
\draw (0mm,20mm)--(0mm,0mm)[line width=1.5pt,draw=babyblue];

\node[cc] at (0mm,5mm){};

\node at (10mm, -6mm) {$2\overline{143}$};
\end{scope}

\begin{scope}[shift={(160mm,0mm)}]
\draw[line width=0.8pt] (0mm,0mm)--(0mm,20mm);
\draw[line width=0.8pt] (5mm,0mm)--(5mm,20mm);
\draw[line width=0.8pt] (10mm,0mm)--(10mm,20mm);
\draw[line width=0.8pt] (15mm,0mm)--(15mm,20mm);
\draw[line width=0.8pt] (20mm,0mm)--(20mm,20mm);

\draw[line width=0.8pt] (0mm,0mm)--(20mm,0mm);
\draw[line width=0.8pt] (0mm,5mm)--(20mm,5mm);
\draw[line width=0.8pt] (0mm,10mm)--(20mm,10mm);
\draw[line width=0.8pt] (0mm,15mm)--(20mm,15mm);
\draw[line width=0.8pt] (0mm,20mm)--(20mm,20mm);

\draw (10mm,15mm)\judy;\draw (5mm,10mm)\nbx;
\draw (15mm,5mm)\nbx;\draw (0mm,0mm)\judy;

\node at (-2mm,0mm){\red{$1$}};
\node at (20mm,24mm){\orange{$0$}};

\draw (0mm,0mm)--(20mm,0mm)[line width=1.5pt,draw=red];
\draw (20mm,20mm)--(20mm,10mm)--(15mm,5mm)--(15mm,0mm)[line width=1.5pt,draw=orange];

\node[cc] at (15mm,0mm){};

\node at (10mm, -6mm) {$2\overline{34}1$};
\end{scope}
\end{tikzpicture}
\end{center}
\end{exam}

\section{Recurrences of $I_{n,k}^\ob$ and $J_{n,k}^\ob$ on signed involutions}\label{Sec-involup}

By applying double inserting and double deleting operations,
we give combinatorial proofs to the recurrences of $I_{n,k}^\ob$ and $J_{n,k}^\ob$ in Theorems \ref{thmIobnk} and \ref{thmJobnk}.

\subsection{Properties of $d$-types in signed involution grids}
We show that the symmetry of signed involution grids implies the symmetry of $d$-types of its grid points.
Recall that $\mathcal{I}_n^B$ is the set all involutions in $\B_n$.

\begin{prop}\label{lemij=ji}
  Let $\pi\in\mathcal{I}_n^B$ and $P_\pi$ be its grid,
then for $1\leq i,j \leq n+1$,  we have
  \begin{enumerate}[(a)]
    \item $d^+(i,j)=d^+(j,i)$ and $d^-(i,j)=d^-(j,i)$; \label{alemij=ji}
    \item $d^+(i,i)=(0,0)$ or $(1,1)$ and $d^-(i,i)=(0,0)$ or $(1,1)$. \label{blemij=ji}
  \end{enumerate}
\end{prop}

\pf
The statement \eqref{alemij=ji} is directly from
the symmetry of the grid and the symmetry of walk trends of
$p_h^+$-paths (resp. $q_v^+$-paths) and $p_h^-$-paths ($q_v^-$-paths) in Theorem \ref{lempath}, where $p,q\in\{0,1\}$.
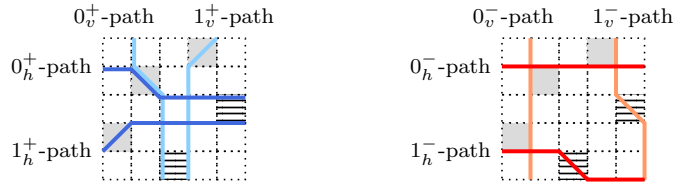
\begin{center}
\begin{tikzpicture}[scale = 0.75]
\def\hezi{-- +(5mm,0mm) -- +(5mm,5mm) -- +(0mm,5mm) -- cycle [line width=0.6pt, dotted]}
\def\judy{-- +(5mm,0mm) -- +(5mm,5mm) -- +(0mm,5mm) -- cycle [line width=0.6pt,fill=gainsboro,dotted]}
\def\nbx{-- +(5mm,0mm) -- +(5mm,5mm) -- +(0mm,5mm) -- cycle [line width=0.6pt,pattern={Lines[angle=45,distance=2.5pt, line width=0.6pt]}, pattern color=black,dotted]}
\tikzstyle{cc}=[circle,draw=black,fill=yellow, line width=0.5pt, inner sep=1.5]

\draw (0mm,0mm)\hezi;
\draw (5mm,0mm)\hezi;
\draw (10mm,0mm)\hezi;
\draw (15mm,0mm)\judy;
\draw (20mm,0mm)\hezi;
\draw (0mm,-5mm)\hezi;
\draw (5mm,-5mm)\judy;
\draw (10mm,-5mm)\hezi;
\draw (15mm,-5mm)\hezi;
\draw (20mm,-5mm)\hezi;
\draw (0mm,-10mm)\hezi;
\draw (5mm,-10mm)\hezi;
\draw (10mm,-10mm)\hezi;
\draw (15mm,-10mm)\hezi;
\draw (20mm,-10mm)\nbx;
\draw (0mm,-15mm)\judy;
\draw (5mm,-15mm)\hezi;
\draw (10mm,-15mm)\hezi;
\draw (15mm,-15mm)\hezi;
\draw (20mm,-15mm)\hezi;
\draw (0mm,-20mm)\hezi;
\draw (5mm,-20mm)\hezi;
\draw (10mm,-20mm)\nbx;
\draw (15mm,-20mm)\hezi;
\draw (20mm,-20mm)\hezi;

\draw (20mm,5mm)--(15mm,0mm)--(15mm,-20mm)[line width=1.5pt,draw=babyblue];

\draw[shift={(0.05,0)}] (5mm,5mm)--(5mm,0mm)--(10mm,-5mm)--(10mm,-15mm)--(10mm,-20mm)[line width=1.5pt,draw=babyblue];

\draw[shift={(0,-0.05)}] (0mm,0mm)--(5mm,0mm)--(10mm,-5mm)--(20mm,-5mm)--(25mm,-5mm)[line width=1.5pt,draw=royalblue];

\draw (0mm,-15mm)--(5mm,-10mm)--(20mm,-10mm)--(25mm,-10mm)[line width=1.5pt,draw=royalblue];


\node at (2mm,9mm){\footnotesize $0_v^+$-path};
\node at (-9mm,0mm){\footnotesize $0_h^+$-path};

\node at (23mm,9mm){\footnotesize $1_v^+$-path};
\node at (-9mm,-15mm){\footnotesize $1_h^+$-path};

\begin{scope}[xshift=70mm]

\draw (0mm,0mm)\hezi;
\draw (5mm,0mm)\hezi;
\draw (10mm,0mm)\hezi;
\draw (15mm,0mm)\judy;
\draw (20mm,0mm)\hezi;
\draw (0mm,-5mm)\hezi;
\draw (5mm,-5mm)\judy;
\draw (10mm,-5mm)\hezi;
\draw (15mm,-5mm)\hezi;
\draw (20mm,-5mm)\hezi;
\draw (0mm,-10mm)\hezi;
\draw (5mm,-10mm)\hezi;
\draw (10mm,-10mm)\hezi;
\draw (15mm,-10mm)\hezi;
\draw (20mm,-10mm)\nbx;
\draw (0mm,-15mm)\judy;
\draw (5mm,-15mm)\hezi;
\draw (10mm,-15mm)\hezi;
\draw (15mm,-15mm)\hezi;
\draw (20mm,-15mm)\hezi;
\draw (0mm,-20mm)\hezi;
\draw (5mm,-20mm)\hezi;
\draw (10mm,-20mm)\nbx;
\draw (15mm,-20mm)\hezi;
\draw (20mm,-20mm)\hezi;

\draw (20mm,5mm)--(20mm,-5mm)--(25mm,-10mm)--(25mm,-20mm)[line width=1.5pt,draw=orange];
\draw(5mm,5mm)--(5mm,-20mm)[line width=1.5pt,draw=orange];

\draw(0mm,0mm)--(25mm,0mm)[line width=1.5pt,draw=red];
\draw (0mm,-15mm)--(10mm,-15mm)--(15mm,-20mm)--(25mm,-20mm)[line width=1.5pt,draw=red];

\node at (2mm,9mm){\footnotesize $0_v^-$-path};
\node at (-9mm,0mm){\footnotesize $0_h^-$-path};

\node at (23mm,9mm){\footnotesize $1_v^-$-path};
\node at (-9mm,-15mm){\footnotesize $1_h^-$-path};
\end{scope}
\end{tikzpicture}
\end{center}
And the statement \eqref{blemij=ji} is directly from \eqref{alemij=ji} by setting $i=j$.
\qed


\begin{prop}\label{lemi+1i}
Let $\pi\in\mathcal{I}_n^B$, then for $1\leq i\leq n$, if $\pi_i=i$, we have
$$
\des^B(\varphi_{(i+1,i)}(\pi))=\des^B(\pi)+1,
$$
and if $\pi_i=\overline{i}$, we have
$$
\des^B(\overline{\varphi}_{(i+1,i)}(\pi)) =\des^B(\pi).
$$
\end{prop}

\proof
Since the grid point $(i+1,i)$ is the southwest corner of the filled square $\langle i,i\rangle$ if $|\pi_i|=i$,
then by Proposition \ref{prop-4corners},
we deduce that
$\des^B(\varphi_{(i+1,i)}(\pi))=\des^B(\pi)+1$ if $\pi_i=i$,
and $\des^B(\overline{\varphi}_{(i+1,i)}(\pi)) =\des^B(\pi)$ if $\pi_i=\overline{i}$.
\qed

\begin{prop}\label{lemxi}
Let $\pi\in\mathcal{I}_n^B$ and $P_\pi$ be its grid, then for $1\leq i\neq j\leq n+1$ and $p,q\in\{0,1\}$, if $d^+(i,j)=(p,q)$, we have
$$\des^B(\xi_{(i,j)}(\pi))=\des^B(\pi)+p+q,$$
and if $d^-(i,j)=(p,q)$, we have
$$\des^B(\overline{\xi}_{(i,j)}(\pi))=\des^B(\pi)+p+q.$$
\end{prop}

\proof
We only proof the case of $d^+(i,j)=(p,q)$.
By the symmetry of $P_\pi$,
the change of idescents after the insertion at the grid point $(j,i)$ can completely reflect
the change of desents after the insertion at the grid point $(i,j)$.
Set $d^+(i,j)=(p,q)$, then by Proposition \ref{lemij=ji}, we have $d^+(j,i)=(p,q)$.
Thus, the operation $\xi_{(i,j)}$ defined in \eqref{xi(ij)}
inserting two positive squares at the points
$(i,j)$ and $(j,i)$ respectively
yields $p$ new descent at the $i$-th row and $q$ new descent at the $j$-th row,
which leads to $\des(\xi_{(i,j)}(\pi))=\des(\pi)+p+q$.
\qed

\begin{prop}\label{lemetaeta'}
Let $\pi\in\mathcal{I}_n^B$ and $P_\pi$ be its grid, then for $1\leq i \leq n+1$ and $p\in\{0,1\}$, if $d^+(i,i)=(p,p)$, we have
$$\des^B(\eta_i(\pi))=\des^B(\pi)+p\quad\mbox{and}\quad
\des^B(\eta'_i(\pi))=\des^B(\pi)+p+1,
$$
and if $d^-(i,i)=(p,p)$, we have
$$\des^B(\overline{\eta}_i(\pi))=\des^B(\pi)+p+1\quad\mbox{and}\quad
\des^B(\overline{\eta'}_i(\pi))=\des^B(\pi)+p.
$$

\end{prop}

\proof
Recall the definitions of $\eta_i, \overline{\eta}_i, \eta_i', \overline{\eta'}_i$ as given in \eqref{eta_i} and \eqref{eta'_i}.
The insertion processes of these four opeartions can be taken as
inserting a positive or negative square at the grid point $(i,i)$, which increases $\des^B(\pi)$ by $p$ depending on $d^+(i,i)=(p,p)$ or $d^-(i,i)=(p,p)$,
then extending the inserted square to a $12$-pair, a $\overline{12}$-pair, a $21$-pair or a $\overline{21}$-pair correspondingly.
Thus the possible additional increase of $\des^B(\pi)$
can only come from the interior of the pair.
Since the $12$-pair and the $\overline{21}$-pair do not contain a $\ob$-descent,
while both the $21$-pair and the $\overline{12}$-pair contain a $\ob$-descent each,
we see that all the formulas in the proposition hold.
\qed

\subsection{Combinatorial proof of the recurrence of $I_{n,k}^\ob$}\label{subsecInkob}
Recall that for $1\leq k\leq n$, $\mathcal{I}_{n,k}^\ob$ is the set of involutions in $\mathcal{I}^B_n$ with $k$ $\ob$-descents, and $|\mathcal{I}_{n,k}^\ob|=I_{n,k}^\ob$.
We now proceed to prove the five-term recurrence \eqref{recIobnk} of $I_{n,k}^\ob$ in Theorem \ref{thmIobnk}.

{\noindent{\emph{Combinatorial Proof of Theorem \ref{thmIobnk}.}\hskip 2pt}}
Let
\[
\mathcal{G}_{n,k}=\{(\sigma,i)\mid \pi\in \mathcal{I}_{n,k}^\ob\text{ and }1\leq i\leq n \},
\]
then we have $|\mathcal{G}_{n,k}|=nI_{n,k}^B$
that equals the left side of \eqref{recIobnk}.
Our next goal is to construct several pairwise disjoint sets
such that their union is counted by the sum of the five terms on the right side of \eqref{recIobnk}.
We specify that for any two grid points, even if their positions in the grid are the same, they are still treated as different
when the signs of the $d$-types we consider are different.

Let
\begin{equation*}
  \begin{aligned}
\mathcal{E}^{(1)^+}_{n,k}=\{(\pi,(i,j))\mid &\text{ $\pi\in\mathcal{I}_{n,k}^\ob$ and $(i,j)$ is the grid point where}\\
    &\text{ a $0_h^+$-path first touches the main diagonal of $P_\pi$}\}
  \end{aligned}
\end{equation*}
and
\begin{equation*}
  \begin{aligned}
\mathcal{E}^{(1)^-}_{n,k}=\{(\pi,(i,j))\mid &\text{ $\pi\in\mathcal{I}_{n,k}^\ob$ and $(i,j)$ is the grid point where a $0_h^-$-path first}\\
    &\text{ touches the main diagonal of $P_\pi$ \emph{or} a negative square on it}\}.
  \end{aligned}
\end{equation*}
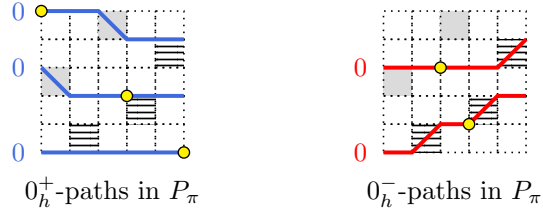
\begin{center}
\begin{tikzpicture}[scale = 0.75]
\def\hezi{-- +(5mm,0mm) -- +(5mm,5mm) -- +(0mm,5mm) -- cycle [line width=0.6pt, dotted]}
\def\judy{-- +(5mm,0mm) -- +(5mm,5mm) -- +(0mm,5mm) -- cycle [line width=0.6pt,fill=gainsboro,dotted]}
\def\nbx{-- +(5mm,0mm) -- +(5mm,5mm) -- +(0mm,5mm) -- cycle [line width=0.6pt,pattern={Lines[angle=45,distance=2.5pt, line width=0.6pt]}, pattern color=black,dotted]}
\tikzstyle{cc}=[circle,draw=black,fill=yellow, line width=0.5pt, inner sep=1.5]

\draw (0mm,0mm)\hezi;
\draw (5mm,0mm)\hezi;
\draw (10mm,0mm)\judy;
\draw (15mm,0mm)\hezi;
\draw (20mm,0mm)\hezi;
\draw (0mm,-5mm)\hezi;
\draw (5mm,-5mm)\hezi;
\draw (10mm,-5mm)\hezi;
\draw (15mm,-5mm)\hezi;
\draw (20mm,-5mm)\nbx;
\draw (0mm,-10mm)\judy;
\draw (5mm,-10mm)\hezi;
\draw (10mm,-10mm)\hezi;
\draw (15mm,-10mm)\hezi;
\draw (20mm,-10mm)\hezi;
\draw (0mm,-15mm)\hezi;
\draw (5mm,-15mm)\hezi;
\draw (10mm,-15mm)\hezi;
\draw (15mm,-15mm)\nbx;
\draw (20mm,-15mm)\hezi;
\draw (0mm,-20mm)\hezi;
\draw (5mm,-20mm)\nbx;
\draw (10mm,-20mm)\hezi;
\draw (15mm,-20mm)\hezi;
\draw (20mm,-20mm)\hezi;

\draw (0mm,5mm)--(10mm,5mm)--(15mm,0mm)--(25mm,0mm)[line width=1.5pt,draw=royalblue];
\node[cc] at (0mm,5mm){};
\node at (-4mm,5mm){\royalblue{$0$}};

\draw (0mm,-5mm)--(5mm,-10mm)--(25mm,-10mm)[line width=1.5pt,draw=royalblue];
\node[cc] at (15mm,-10mm){};
\node at (-4mm,-5mm){\royalblue{$0$}};

\draw (0mm,-20mm)--(25mm,-20mm)[line width=1.5pt,draw=royalblue];
\node[cc] at (25mm,-20mm){};
\node at (-4mm,-20mm){\royalblue{$0$}};

\node at (12.5mm,-27mm){$0_h^+$-paths in $P_\pi$};

\begin{scope}[xshift=60mm]
\draw (0mm,0mm)\hezi;
\draw (5mm,0mm)\hezi;
\draw (10mm,0mm)\judy;
\draw (15mm,0mm)\hezi;
\draw (20mm,0mm)\hezi;
\draw (0mm,-5mm)\hezi;
\draw (5mm,-5mm)\hezi;
\draw (10mm,-5mm)\hezi;
\draw (15mm,-5mm)\hezi;
\draw (20mm,-5mm)\nbx;
\draw (0mm,-10mm)\judy;
\draw (5mm,-10mm)\hezi;
\draw (10mm,-10mm)\hezi;
\draw (15mm,-10mm)\hezi;
\draw (20mm,-10mm)\hezi;
\draw (0mm,-15mm)\hezi;
\draw (5mm,-15mm)\hezi;
\draw (10mm,-15mm)\hezi;
\draw (15mm,-15mm)\nbx;
\draw (20mm,-15mm)\hezi;
\draw (0mm,-20mm)\hezi;
\draw (5mm,-20mm)\nbx;
\draw (10mm,-20mm)\hezi;
\draw (15mm,-20mm)\hezi;
\draw (20mm,-20mm)\hezi;

\draw (0mm,-5mm)--(20mm,-5mm)--(25mm,0mm)[line width=1.5pt,draw=red];
\node[cc] at (10mm,-5mm){};
\node at (-4mm,-5mm){\red{$0$}};

\draw (0mm,-20mm)--(5mm,-20mm)--(10mm,-15mm)--(15mm,-15mm)--(20mm,-10mm)--(25mm,-10mm)[line width=1.5pt,draw=red];
\node[cc] at (15mm,-15mm){};
\node at (-4mm,-20mm){\red{$0$}};

\node at (12.5mm,-27mm){$0_h^-$-paths in $P_\pi$};
\end{scope}
\end{tikzpicture}
\end{center}
Since there are $(k+1)$ $0_h^+$-paths and $k$ $0_h^-$-paths in $P_\pi$ for  $\pi\in \mathcal{I}_{n-1,k}^B$ from Theorem \ref{thmpaths},
we have $$\left|\mathcal{E}^{(1)^+}_{n-1,k}\uplus\mathcal{E}^{(1)^-}_{n-1,k}\right|
=(k+1)I_{n-1,k}^\ob+kI_{n-1,k}^\ob=(2k+1)I_{n-1,k}^\ob,$$ which is the \textbf{first} term on the right side of \eqref{recIobnk}.
For any $(\pi,(i,j))\in\mathcal{E}^{(1)^+}_{n-1,k}\uplus \mathcal{E}^{(1)^-}_{n-1,k}$,
set
\begin{equation*}
  \Psi((\pi,(i,j)))=(\sigma,i),
\end{equation*}
where
\begin{equation}\label{eqpfIn1st}
\sigma=\left\{
      \begin{array}{ll}
        \varphi_{(i,j)}(\pi), & \hbox{if $(\pi,(i,j))\in\mathcal{E}^{(1)^+}_{n-1,k}$ ,} \\[6pt]
        \overline{\varphi}_{(i,j)}(\pi), & \hbox{if $(\pi,(i,j))\in\mathcal{E}^{(1)^-}_{n-1,k}$.}
      \end{array}
    \right.
\end{equation}
Note that the grid point $(i,j)$ must have $d^+$-type or $d^-$-type $(0,0)$ by Propositions \ref{prop-4corners} and \ref{lemij=ji},
which leads to $\des^B(\sigma)=\des^B(\pi)=k$.
Moreover $\sigma\in\mathcal{I}^B_n$ by Proposition \ref{propi+1i}.
Hence  we have $(\sigma,i)\in\mathcal{G}_{n,k}$.

Let
\begin{equation*}
  \begin{aligned}
\mathcal{E}^{(2)^+}_{n,k}=\{(\pi,(i,j))\mid& \text{ $\pi\in\mathcal{I}_{n,k}^\ob$ and $(i,j)$ is grid point where a $1_h^+$-path}\\
&\text{ touches the main diagonal of $P_\pi$ \emph{or} a positive square on it}\},
  \end{aligned}
\end{equation*}
and
\begin{equation*}
  \begin{aligned}
\mathcal{E}^{(2)^-}_{n,k}=\{(\pi,(i,j))\mid& \text{ $\pi\in\mathcal{I}_{n,k}^\ob$ and $(i,j)$ is the grid point where}\\
    &\text{ a $1_h^-$-path first touches the main diagonal of $P_\pi$}\}.
  \end{aligned}
\end{equation*}
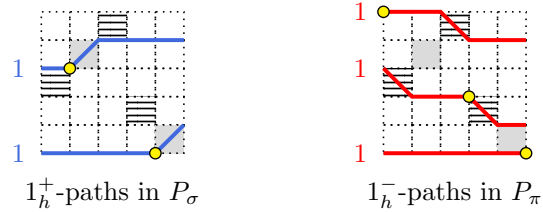
\begin{center}
\begin{tikzpicture}[scale = 0.75]
\def\hezi{-- +(5mm,0mm) -- +(5mm,5mm) -- +(0mm,5mm) -- cycle [line width=0.6pt, dotted]}
\def\judy{-- +(5mm,0mm) -- +(5mm,5mm) -- +(0mm,5mm) -- cycle [line width=0.6pt,fill=gainsboro,dotted]}
\def\nbx{-- +(5mm,0mm) -- +(5mm,5mm) -- +(0mm,5mm) -- cycle [line width=0.6pt,pattern={Lines[angle=45,distance=2.5pt, line width=0.6pt]}, pattern color=black,dotted]}
\tikzstyle{cc}=[circle,draw=black,fill=yellow, line width=0.5pt, inner sep=1.5]

\draw (0mm,0mm)\hezi;
\draw (5mm,0mm)\hezi;
\draw (10mm,0mm)\nbx;
\draw (15mm,0mm)\hezi;
\draw (20mm,0mm)\hezi;
\draw (0mm,-5mm)\hezi;
\draw (5mm,-5mm)\judy;
\draw (10mm,-5mm)\hezi;
\draw (15mm,-5mm)\hezi;
\draw (20mm,-5mm)\hezi;
\draw (0mm,-10mm)\nbx;
\draw (5mm,-10mm)\hezi;
\draw (10mm,-10mm)\hezi;
\draw (15mm,-10mm)\hezi;
\draw (20mm,-10mm)\hezi;
\draw (0mm,-15mm)\hezi;
\draw (5mm,-15mm)\hezi;
\draw (10mm,-15mm)\hezi;
\draw (15mm,-15mm)\nbx;
\draw (20mm,-15mm)\hezi;
\draw (0mm,-20mm)\hezi;
\draw (5mm,-20mm)\hezi;
\draw (10mm,-20mm)\hezi;
\draw (15mm,-20mm)\hezi;
\draw (20mm,-20mm)\judy;

\draw (0mm,-5mm)--(5mm,-5mm)--(10mm,0mm)--(25mm,0mm)[line width=1.5pt,draw=royalblue];
\node[cc] at (5mm,-5mm){};
\node at (-4mm,-5mm){\royalblue{$1$}};

\draw (0mm,-20mm)--(20mm,-20mm)--(25mm,-15mm)[line width=1.5pt,draw=royalblue];
\node[cc] at (20mm,-20mm){};
\node at (-4mm,-20mm){\royalblue{$1$}};

\node at (12.5mm,-27mm){$1_h^+$-paths in $P_\sigma$};

\begin{scope}[xshift=60mm]
\draw (0mm,0mm)\hezi;
\draw (5mm,0mm)\hezi;
\draw (10mm,0mm)\nbx;
\draw (15mm,0mm)\hezi;
\draw (20mm,0mm)\hezi;
\draw (0mm,-5mm)\hezi;
\draw (5mm,-5mm)\judy;
\draw (10mm,-5mm)\hezi;
\draw (15mm,-5mm)\hezi;
\draw (20mm,-5mm)\hezi;
\draw (0mm,-10mm)\nbx;
\draw (5mm,-10mm)\hezi;
\draw (10mm,-10mm)\hezi;
\draw (15mm,-10mm)\hezi;
\draw (20mm,-10mm)\hezi;
\draw (0mm,-15mm)\hezi;
\draw (5mm,-15mm)\hezi;
\draw (10mm,-15mm)\hezi;
\draw (15mm,-15mm)\nbx;
\draw (20mm,-15mm)\hezi;
\draw (0mm,-20mm)\hezi;
\draw (5mm,-20mm)\hezi;
\draw (10mm,-20mm)\hezi;
\draw (15mm,-20mm)\hezi;
\draw (20mm,-20mm)\judy;

\draw (0mm,5mm)--(10mm,5mm)--(15mm,0mm)--(25mm,0mm)[line width=1.5pt,draw=red];
\node[cc] at (0mm,5mm){};
\node at (-4mm,5mm){\red{$1$}};

\draw (0mm,-5mm)--(5mm,-10mm)--(15mm,-10mm)--(20mm,-15mm)--(25mm,-15mm)[line width=1.5pt,draw=red];
\node[cc] at (15mm,-10mm){};
\node at (-4mm,-5mm){\red{$1$}};

\draw (0mm,-20mm)--(25mm,-20mm)[line width=1.5pt,draw=red];
\node[cc] at (25mm,-20mm){};
\node at (-4mm,-20mm){\red{$1$}};

\node at (12.5mm,-27mm){$1_h^-$-paths in $P_\pi$};
\end{scope}
\end{tikzpicture}
\end{center}
It follows from Theorem \ref{thmpaths} that
the number of $1_h^+$-paths and $1_h^-$-paths in the grid $P_\pi$ with $\pi\in\mathcal{I}_{n-1,k-1}^B$
is $n-k$ and $n-k+1$, respectively.
Hence we deduce
\[
\left|\mathcal{E}^{(2)^+}_{n-1,k-1}\uplus\mathcal{E}^{(2)^-}_{n-1,k-1}\right|
=(n-k)I_{n-1,k-1}^\ob+(n-k+1)I_{n-1,k-1}^\ob
=(2n-2k+1)I_{n-1,k-1}^\ob
\]
that is the \textbf{second} term on the right side of \eqref{recIobnk}.
For any $(\pi,(i,j))\in \mathcal{E}^{(2)^+}_{n-1,k-1}\uplus\mathcal{E}^{(2)^-}_{n-1,k-1}$,
define
\begin{equation*}
\Psi((\pi, (i,j)))=(\sigma,i),
\end{equation*}
where
\begin{equation}\label{eqpfIn2nd}
\sigma=\left\{
      \begin{array}{ll}
        \varphi_{(i,j)}(\pi), & \hbox{if $(\pi,(i,j))\in\mathcal{E}^{(2)^+}_{n-1,k-1}$ ,} \\[6pt]
        \overline{\varphi}_{(i,j)}(\pi), & \hbox{if $(\pi,(i,j))\in\mathcal{E}^{(2)^-}_{n-1,k-1}$.}
      \end{array}
    \right.
\end{equation}
By Propositions \ref{prop-4corners} and \ref{lemij=ji},
the grid $(i,j)$ is of $d^+$-type or $d^-$-type $(1,1)$,
which implies $\des^B(\sigma)=\des^B(\pi)+1$.
And by Proposition \ref{propi+1i}, we have $\sigma\in\mathcal{I}_{n}^B$.
Thus $(\sigma,i)\in\mathcal{G}_{n,k}$.

Let
\begin{equation*}
\mathcal{E}^{(3)^+}_{n,k}=\{(\pi,(i,j))\mid \pi\in \mathcal{I}_{n,k}^\ob\text{ and $d^+(i,j)=(0,0)$ in $P_\pi$}\}
\end{equation*}
and
\begin{equation*}
\mathcal{E}^{(3)^-}_{n,k}=\{(\pi,(i,j))\mid \pi\in \mathcal{I}_{n,k}^\ob\text{ and $d^-(i,j)=(0,0)$ in $P_\pi$}\}.
\end{equation*}
By Theorem \ref{thm-enum-dtype}
we have
\begin{equation*}
\begin{aligned}
\left|\mathcal{E}^{(3)^+}_{n-2,k}\uplus\mathcal{E}^{(3)^-}_{n-2,k}\right|
&=\sum_{\pi\in\mathcal{I}_{n-2,k}^\ob}(k+1)^2-n(\pi)+(n-2)+k^2+n(\pi)\\
&=(n-1+2k(k+1))I_{n-2,k}^\ob,
\end{aligned}
\end{equation*}
which is equal to the \textbf{third} term on the right side of \eqref{recIobnk}.
For any $(\sigma,(i,j))\in \mathcal{E}^{(3)^+}_{n-2,k}\uplus\mathcal{E}^{(3)^-}_{n-2,k}$,
define
\begin{equation*}
\Psi((\pi, (i,j)))=(\sigma,\chi_{ij}),
\end{equation*}
where
\begin{equation*}
\chi_{ij}=\left\{
\begin{array}{ll}
i, & \hbox{if  $i\leq j$,} \\[6pt]
i+1, & \hbox{if $i>j$;}
\end{array}
    \right.
\end{equation*}
and
\begin{equation}\label{eqpfIn3rd}
\sigma=\left\{
\begin{array}{ll}
\xi_{(i,j)}(\pi), & \hbox{if $(\pi,(i,j))\in\mathcal{E}^{(3)^+}_{n-2,k}$ and $i\neq j$,} \\[6pt]
\eta_{i}(\pi), & \hbox{if $(\pi,(i,j))\in\mathcal{E}^{(3)^+}_{n-2,k}$ and $i=j$,} \\[6pt]
\overline{\xi}_{(i,j)}(\pi), & \hbox{if $(\pi,(i,j))\in\mathcal{E}^{(3)^-}_{n-2,k}$ and $i\neq j$,} \\[6pt]
\overline{\eta'}_{i}(\pi), & \hbox{if $(\pi,(i,j))\in\mathcal{E}^{(3)^-}_{n-2,k}$ and $i=j$.}
\end{array}
    \right.
\end{equation}
Then by Propositions \ref{lemxi} and \ref{lemetaeta'}, we see  $(\sigma,\chi_{ij})\in\mathcal{G}_{n,k}$.

Let
\begin{equation*}
  \begin{aligned}
  \mathcal{E}^{(4)^+}_{n,k}=&\{(\pi,(i,j))\mid\text{ $\pi\in\mathcal{I}_{n,k}^\ob$ and $d^+(i,j)=(1,0)$ or $(0,1)$ in $P_\pi$}\}\\
  &\quad\cup\{(\pi,(i,i))\mid\text{ $\pi\in\mathcal{I}_{n,k}^\ob$ and $d^+(i,i)=(0,0)$ or $(1,1)$ in $P_\pi$}\},
  \end{aligned}
\end{equation*}
and
\begin{equation*}
  \begin{aligned}
  \mathcal{E}^{(4)^-}_{n,k}=&\{(\pi,(i,j))\mid\text{ $\pi\in\mathcal{I}_{n,k}^\ob$ and $d^-(i,j)=(1,0)$ or $(0,1)$ in $P_\pi$}\}\\
  &\quad\cup\{(\pi,(i,i))\mid\text{ $\pi\in\mathcal{I}_{n,k}^\ob$ and $d^-(i,i)=(0,0)$ or $(1,1)$ in $P_\pi$}\}.
  \end{aligned}
\end{equation*}
From Proposition \ref{lemij=ji}, we know that the two subsets in the set $\mathcal{E}^{(4)^+}_{n,k}$ (resp. $\mathcal{E}^{(4)^-}_{n,k}$) are actually disjoint.
Thus, with Theorem \ref{thm-enum-dtype}, we see
\begin{equation*}
\begin{aligned}
&\left|\mathcal{E}^{(4)^+}_{n-2,k-1}\uplus\mathcal{E}^{(4)^-}_{n-2,k-1}\right|\\
=&\sum_{\pi\in\mathcal{I}_{n-2,k-1}^\ob}2(k(n-k-1)+n(\pi)-(n-2))+
2((k-1)(n-k)-n(\pi))+2(n-1)\\
=&(2(n-1)+4(n-k-1)(k-1)) I_{n-2, k-1}^\ob,
\end{aligned}
\end{equation*}
which is the \textbf{fourth} term on the right side of \eqref{recIobnk}.
For any $(\sigma,(i,j))\in \mathcal{E}^{(4)^+}_{n-2,k-1}\uplus\mathcal{E}^{(4)^-}_{n-2,k-1}$,
we define
\begin{equation*}
\Psi((\pi, (i,j)))=(\sigma,\chi_{ij}),
\end{equation*}
where
\begin{equation}\label{eqpfIn4th}
(\sigma,k)=\left\{
      \begin{array}{ll}
       \xi_{(i,j)}(\pi), & \hbox{if $(\pi,(i,j))\in\mathcal{E}^{(4)^+}_{n-2,k-1}$ and $i\neq j$,} \\[6pt]
      \eta'_{i}(\pi), & \hbox{if $(\pi,(i,i))\in\mathcal{E}^{(4)^+}_{n-2,k-1}$ and $d^+(i,i)=(0,0)$,} \\[6pt]
      \eta_{i}(\pi), & \hbox{if $(\pi,(i,i))\in\mathcal{E}^{(4)^+}_{n-2,k-1}$ and $d^+(i,i)=(1,1)$,} \\[6pt]
      \overline{\xi}_{(i,j)}(\pi), & \hbox{if $(\pi,(i,j))\in\mathcal{E}^{(4)^-}_{n-2,k-1}$ and $i\neq j$,} \\[6pt]
      \overline{\eta}_{i}(\pi), & \hbox{if $(\pi,(i,i))\in\mathcal{E}^{(4)^-}_{n-2,k-1}$ and $d^-(i,i)=(0,0)$,} \\[6pt]
      \overline{\eta'}_{i}(\pi), & \hbox{if $(\pi,(i,i))\in\mathcal{E}^{(4)^-}_{n-2,k-1}$ and $d^-(i,i)=(1,1)$.}
      \end{array}
    \right.
\end{equation}
Combining Propositions \ref{lemxi} and \ref{lemetaeta'}, we have $(\sigma,\chi_{ij})\in\mathcal{G}_{n,k}$.

Let
\begin{equation*}
\mathcal{E}^{(5)^+}_{n,k}=\{(\pi,(i,j))\mid \pi\in \mathcal{I}_{n,k}^\ob\text{ and $d^+(i,j)=(1,1)$ in $P_\pi$}\}
\end{equation*}
and
\begin{equation*}
\mathcal{E}^{(5)^-}_{n,k}=\{(\pi,(i,j))\mid \pi\in \mathcal{I}_{n,k}^\ob\text{ and $d^-(i,j)=(1,1)$ in $P_\pi$}\}.
\end{equation*}
It follows from Theorem \ref{thm-enum-dtype} that
we have
\begin{equation*}
\begin{aligned}
\left|\mathcal{E}^{(5)^+}_{n-2,k-2}\uplus\mathcal{E}^{(5)^-}_{n-2,k-2}\right|
&=\sum_{\pi\in\mathcal{I}_{n-2,k-2}^\ob}(n-k)^2-n(\pi)+(n-2)+(n-k+1)^2+n(\pi)\\
&=\left((n-k)^2+(n-k+1)^2+(n-2)\right) I_{n-2, k-2}^\ob,
\end{aligned}
\end{equation*}
which equals to the \textbf{fifth} term on the right side of \eqref{recIobnk} after a simple calculation.
For any $(\sigma,(i,j))\in \mathcal{E}^{(5)^+}_{n-2,k-2}\uplus\mathcal{E}^{(5)^-}_{n-2,k-2}$,
define
\begin{equation*}
\Psi((\pi, (i,j)))=(\sigma,\chi_{ij}),
\end{equation*}
where
\begin{equation}\label{eqpfIn5th}
\sigma=\left\{
      \begin{array}{ll}
       \xi_{(i,j)}(\pi), & \hbox{if $(\pi,(i,j))\in\mathcal{E}^{(5)^+}_{n-2,k-2}$ and $i\neq j$,} \\[6pt]
      \eta'_{i}(\pi), & \hbox{if $(\pi,(i,j))\in\mathcal{E}^{(5)^+}_{n-2,k-2}$ and $i=j$,} \\[6pt]
      \overline{\xi}_{(i,j)}(\pi), & \hbox{if $(\pi,(i,j))\in\mathcal{E}^{(5)^-}_{n-2,k-2}$ and $i\neq j$,} \\[6pt]
      \overline{\eta}_{i}(\pi), & \hbox{if $(\pi,(i,j))\in\mathcal{E}^{(5)^-}_{n-2,k-2}$ and $i=j$.}
      \end{array}
    \right.
\end{equation}
Hence we deduce $(\sigma,\chi_{ij})\in\mathcal{G}_{n,k}$ by Propositions \ref{lemxi} and \ref{lemetaeta'}.

It is trivial that the above sets $\mathcal{E}^{(1)^\pm}_{n-1,k}$,
$\mathcal{E}^{(2)^\pm}_{n-1,k-1}$, $\mathcal{E}^{(3)^\pm}_{n-2,k}$, $\mathcal{E}^{(4)^\pm}_{n-2,k-1}$, $\mathcal{E}^{(5)^\pm}_{n-2,k-2}$ are pairwise disjoint,
where $\mathcal{E}^{(r)\pm}_{n,k}$ denotes the disjoint union of $\mathcal{E}^{(r)+}_{n,k}$ and $\mathcal{E}^{(r)-}_{n,k}$ for $1\leq r\leq 5$.
Set
\begin{equation*}
{\mathcal{E}}_{n,k}=
\mathcal{E}^{(1)^\pm}_{n-1,k}\uplus\mathcal{E}^{(2)^\pm}_{n-1,k-1}\uplus\mathcal{E}^{(3)^\pm}_{n-2,k}
\uplus\mathcal{E}^{(4)^\pm}_{n-2,k-1}\uplus\mathcal{E}^{(5)^\pm}_{n-2,k-2},
\end{equation*}
then $|{\mathcal{E}}_{n,k}|$ is  exactly equal to the right side of \eqref{recIobnk}.
With the combination of \eqref{eqpfIn1st}--\eqref{eqpfIn5th},
we establish the mapping
\[
\Psi\colon{\mathcal{E}}_{n,k}\rightarrow{\mathcal{G}}_{n,k}.
\]
In order to complete the proof, the remaining is devoted to giving the inverse of $\Psi$.

For any $(\sigma, r)\in\mathcal{G}_{n,k}$, we set
\begin{equation*}
\Theta((\sigma,r))=(\pi,(i,j)),
\end{equation*}
where
\begin{equation}\label{theta}
(\pi,(i,j))=\left\{
 \begin{array}{ll}
 (\xi^{-1}_{\langle r,\sigma_r-1\rangle}(\sigma),(r,\sigma_r-1)), & r+1<\sigma_r\leq n, \\[8pt]
 (\xi^{-1}_{\langle r-1,\sigma_r\rangle}(\sigma),(r-1,\sigma_r)), & 1\leq \sigma_r<r-1, \\[8pt]
 (\eta'^{-1}_r(\sigma),(r,\sigma_r-1)), & \sigma_r=r+1, \\[8pt]
 (\varphi^{-1}_{\langle r,\sigma_r\rangle}(\sigma),(r,\sigma_r)), & \sigma_r=r-1, \\[8pt]
 (\varphi^{-1}_{\langle r,\sigma_r\rangle}(\sigma),(r,\sigma_r)), & \text{$\sigma_r=r$ for $r=1$, \emph{or} $r\geq 2$ with $\sigma_{r-1}\neq r-1$}, \\[8pt]
 (\eta^{-1}_{r-1}(\sigma),(r-1,\sigma_{r}-1)), & \text{$\sigma_r=r$ for $r\geq 2$ with $\sigma_{r-1}=r-1$}, \\[8pt]
 (\overline{\xi}^{-1}_{\langle r,|\sigma_r|-1\rangle}(\sigma),(r,|\sigma_r|-1)), & \overline{n}\leq \sigma_r< \overline{r+1}, \\[8pt]
 (\overline{\xi}^{-1}_{\langle r-1,|\sigma_r|\rangle}(\sigma),(r-1,|\sigma_r|)), & \overline{r-1}<\sigma_r\leq \overline{1}, \\[8pt]
 (\overline{\eta'}^{-1}_r(\sigma),(r,|\sigma_r|-1)), & \sigma_r=\overline{r+1}, \\[8pt]
 (\overline{\varphi}^{-1}_{\langle r,|\sigma_r|\rangle}(\sigma),(r,|\sigma_r|)), & \sigma_r=\overline{r-1}, \\[8pt]
 (\overline{\varphi}^{-1}_{\langle r,|\sigma_r|\rangle}(\sigma),(r,|\sigma_r|)), & \text{$\sigma_r=\overline{r}$ for $r=1$, \emph{or} $r \geq 2$ with $\sigma_r\neq \overline{r-1}$}, \\[8pt]
 (\overline{\eta}^{-1}_{r-1}(\sigma),(r-1,|\sigma_{r-1}|)), & \text{$\sigma_r=\overline{r}$ for $r \geq 2$ with $\sigma_r= \overline{r-1}$}.
 \end{array}
 \right.
\end{equation}
We show that the mapping $\Theta$ is the inverse of $\Psi$ based on the different  values of $\sigma_r$.

\noindent{\bf{Case 1: $\bm{r+1<|\sigma_r |\leq n}$.}}
By the definitions of $\xi^{-1}_{\ij}$ and $\overline{\xi}^{-1}_{\ij}$ with $i<j$ as given in \eqref{invxi(ij)} and \eqref{invbxi(ij)}, respectively, and Proposition \ref{lemxi}, we have
$$
\pi=\xi^{-1}_{{\langle r,\sigma_r-1\rangle}}(\sigma)\in\mathcal{I}^\ob_{n-2,k-p-q}
$$
if $\sigma_i >0$ and $(r,\sigma_r-1)$ is of $d^+$-type $(p,q)$ in $P_\pi$ for $p,q\in\{0,1\}$,
and
$$
\pi=\overline{\xi}^{-1}_{{\langle r,|\sigma_r|-1\rangle}}(\sigma)\in\mathcal{I}^\ob_{n-2,k-p'-q'}
$$
if $\sigma_i<0$ and $(r,|\sigma_r|-1)$ is of $d^-$-type $(p',q')$ in $P_\pi$ for $p',q'\in\{0,1\}$.
See the following figures for examples.
\begin{center}
\begin{tikzpicture}[scale = 0.75]
\def\hezi{-- +(5mm,0mm) -- +(5mm,5mm) -- +(0mm,5mm) -- cycle [line width=0.6pt]}
\def\judy{-- +(5mm,0mm) -- +(5mm,5mm) -- +(0mm,5mm) -- cycle [line width=0.6pt,fill=gainsboro]}
\def\nbx{-- +(5mm,0mm) -- +(5mm,5mm) -- +(0mm,5mm) -- cycle [line width=0.6pt,pattern={Lines[angle=45,distance=2.5pt, line width=0.6pt]}, pattern color=black]}
\tikzstyle{bdot}=[circle,fill=royalblue,draw=royalblue,inner sep=2]
\tikzstyle{rdot}=[circle,fill=red,draw=red,inner sep=2]

\draw (0mm,0mm)\hezi;
\draw (5mm,0mm)\hezi;
\draw (10mm,0mm)\hezi;
\draw (15mm,0mm)\judy;
\draw (20mm,0mm)\hezi;
\draw (0mm,-5mm)\hezi;
\draw (5mm,-5mm)\nbx;
\draw (10mm,-5mm)\hezi;
\draw (15mm,-5mm)\hezi;
\draw (20mm,-5mm)\hezi;
\draw (0mm,-10mm)\hezi;
\draw (5mm,-10mm)\hezi;
\draw (10mm,-10mm)\hezi;
\draw (15mm,-10mm)\hezi;
\draw (20mm,-10mm)\nbx;
\draw (0mm,-15mm)\judy;
\draw (5mm,-15mm)\hezi;
\draw (10mm,-15mm)\hezi;
\draw (15mm,-15mm)\hezi;
\draw (20mm,-15mm)\hezi;
\draw (0mm,-20mm)\hezi;
\draw (5mm,-20mm)\hezi;
\draw (10mm,-20mm)\nbx;
\draw (15mm,-20mm)\hezi;
\draw (20mm,-20mm)\hezi;

\draw (15mm,0mm)-- +(5mm,0mm) -- +(5mm,5mm) -- +(0mm,5mm) -- cycle [line width=1.5pt,draw=royalblue];

\draw (20mm,-10mm)-- +(5mm,0mm) -- +(5mm,5mm) -- +(0mm,5mm) -- cycle [line width=1.5pt,draw=red];

\node at (12.5mm,-26mm){$\sigma=4\overline{25}1\overline{3}$};


\begin{scope}[shift={(50mm,-5.5mm)},scale=1.25]
\draw (0mm,0mm)\nbx;
\draw (5mm,0mm)\hezi;
\draw (10mm,0mm)\hezi;
\draw (0mm,-5mm)\hezi;
\draw (5mm,-5mm)\hezi;
\draw (10mm,-5mm)\nbx;
\draw (0mm,-10mm)\hezi;
\draw (5mm,-10mm)\nbx;
\draw (10mm,-10mm)\hezi;

\node[bdot] at(10mm,5mm){};

\node at (7.5mm,-16mm){$\xi^{-1}_{\langle 1,3\rangle}(\sigma)$};

\end{scope}


\begin{scope}[shift={(100mm,-5mm)},scale=1.25]
\draw (0mm,0mm)\hezi;
\draw (5mm,0mm)\hezi;
\draw (10mm,0mm)\judy;
\draw (0mm,-5mm)\hezi;
\draw (5mm,-5mm)\nbx;
\draw (10mm,-5mm)\hezi;
\draw (0mm,-10mm)\judy;
\draw (5mm,-10mm)\hezi;
\draw (10mm,-10mm)\hezi;

\node[rdot] at(15mm,-5mm){};
\node at (7.5mm,-16mm){$\overline{\xi}^{-1}_{\langle 3,4\rangle}(\sigma)$};

\end{scope}
\end{tikzpicture}
\end{center}
Thus by \eqref{eqpfIn3rd}, \eqref{eqpfIn4th} and \eqref{eqpfIn5th},
we see $\Psi(\pi, (r,|\sigma_r|-1))=(\sigma,r)$.
%

\noindent{\bf{Case 2: $\bm{1\leq |\sigma_r|\leq r-1}$.}}
Similar to the Case 1, according to \eqref{invxi(ij)} and \eqref{invbxi(ij)} with the case $i>j$ and Proposition \ref{lemxi},
we arrive that
$$
\pi=\xi^{-1}_{\langle r-1,\sigma_r\rangle}(\sigma)\in\mathcal{I}^\ob_{n-2,k-p-q}
$$
if $\sigma_i>0$ and the $d^+$-type of the grid point $(r-1,\sigma_r)$ in $P_\pi$ is $(p,q)$ for $p,q\in\{0,1\}$, and
$$
\pi=\overline{\xi}^{-1}_{\langle r-1,|\sigma_r|\rangle}(\sigma)\in\mathcal{I}^\ob_{n-2,k-p'-q'}
$$
if $\sigma_i<0$ and the $d^-$-type of the grid point $(r-1,|\sigma_r|)$ is $(p',q')$ in $P_\pi$ for $p',q'\in\{0,1\}$.
One may refer the following figures for a clearer view.
\begin{center}
\begin{tikzpicture}[scale = 0.75]
\def\hezi{-- +(5mm,0mm) -- +(5mm,5mm) -- +(0mm,5mm) -- cycle [line width=0.6pt]}
\def\judy{-- +(5mm,0mm) -- +(5mm,5mm) -- +(0mm,5mm) -- cycle [line width=0.6pt,fill=gainsboro]}
\def\nbx{-- +(5mm,0mm) -- +(5mm,5mm) -- +(0mm,5mm) -- cycle [line width=0.6pt,pattern={Lines[angle=45,distance=2.5pt, line width=0.6pt]}, pattern color=black]}
\tikzstyle{bdot}=[circle,fill=royalblue,draw=royalblue,inner sep=2]
\tikzstyle{rdot}=[circle,fill=red,draw=red,inner sep=2]

\draw (0mm,0mm)\hezi;
\draw (5mm,0mm)\hezi;
\draw (10mm,0mm)\hezi;
\draw (15mm,0mm)\judy;
\draw (20mm,0mm)\hezi;
\draw (0mm,-5mm)\hezi;
\draw (5mm,-5mm)\nbx;
\draw (10mm,-5mm)\hezi;
\draw (15mm,-5mm)\hezi;
\draw (20mm,-5mm)\hezi;
\draw (0mm,-10mm)\hezi;
\draw (5mm,-10mm)\hezi;
\draw (10mm,-10mm)\hezi;
\draw (15mm,-10mm)\hezi;
\draw (20mm,-10mm)\nbx;
\draw (0mm,-15mm)\judy;
\draw (5mm,-15mm)\hezi;
\draw (10mm,-15mm)\hezi;
\draw (15mm,-15mm)\hezi;
\draw (20mm,-15mm)\hezi;
\draw (0mm,-20mm)\hezi;
\draw (5mm,-20mm)\hezi;
\draw (10mm,-20mm)\nbx;
\draw (15mm,-20mm)\hezi;
\draw (20mm,-20mm)\hezi;

\draw (0mm,-15mm)-- +(5mm,0mm) -- +(5mm,5mm) -- +(0mm,5mm) -- cycle [line width=1.5pt,draw=royalblue];

\draw (10mm,-20mm)-- +(5mm,0mm) -- +(5mm,5mm) -- +(0mm,5mm) -- cycle [line width=1.5pt,draw=red];

\node at (12.5mm,-26mm){$\sigma=4\overline{25}1\overline{3}$};


\begin{scope}[shift={(50mm,-5.5mm)},scale=1.25]
\draw (0mm,0mm)\nbx;
\draw (5mm,0mm)\hezi;
\draw (10mm,0mm)\hezi;
\draw (0mm,-5mm)\hezi;
\draw (5mm,-5mm)\hezi;
\draw (10mm,-5mm)\nbx;
\draw (0mm,-10mm)\hezi;
\draw (5mm,-10mm)\nbx;
\draw (10mm,-10mm)\hezi;

\node[bdot] at(0mm,-5mm){};

\node at (7.5mm,-16mm){$\xi^{-1}_{\langle 3,1\rangle}(\sigma)$};

\end{scope}


\begin{scope}[shift={(100mm,-5mm)},scale=1.25]
\draw (0mm,0mm)\hezi;
\draw (5mm,0mm)\hezi;
\draw (10mm,0mm)\judy;
\draw (0mm,-5mm)\hezi;
\draw (5mm,-5mm)\nbx;
\draw (10mm,-5mm)\hezi;
\draw (0mm,-10mm)\judy;
\draw (5mm,-10mm)\hezi;
\draw (10mm,-10mm)\hezi;

\node[rdot] at(10mm,-10mm){};
\node at (7.5mm,-16mm){$\overline{\xi}^{-1}_{\langle 4,3\rangle}(\sigma)$};

\end{scope}
\end{tikzpicture}
\end{center}
Hence $\Psi(\pi, (r-1,|\sigma_r|))=(\sigma,r)$
by \eqref{eqpfIn3rd}, \eqref{eqpfIn4th} and \eqref{eqpfIn5th}.


\noindent{\bf{Case 3: $\bm{|\sigma_r|=r+1}$.}}
By the symmetry of the grid $P_\sigma$ of signed involution $\sigma$,
$|\sigma_r|=r+1$ implies that $\langle r, r+1\rangle$ and $\langle r+1, r\rangle$
in $P_\sigma$ are both positive squares if $\sigma_r>0$
or negative squares if $\sigma_r<0$.
It follows from the definitions of $\eta'^{-1}_i$ and $\overline{\eta'}^{-1}_i$ in \eqref{inveta_i} and Proposition \ref{lemetaeta'} that
$$
\pi=\eta'^{-1}_r(\sigma)\in \mathcal{I}_{n-2,k-1-p}^\ob
$$
if $\sigma_r>0$ and the grid point $(r,\sigma_r-1)$ is of $d^+$-type $(p,p)$ in $P_\pi$ for $p\in\{0,1\}$, and
$$
\pi=\overline{\eta'}^{-1}_r(\sigma)\in \mathcal{I}_{n-2,k-p'}^\ob
$$
if $\sigma_r<0$ and the grid point $(r,|\sigma_r|-1)$ is of $d^-$-type $(p',p')$ in $P_\pi$ for $p'\in\{0,1\}$.
We can see the following figures as an example.
\begin{center}
\begin{tikzpicture}[scale = 0.75]
\def\hezi{-- +(5mm,0mm) -- +(5mm,5mm) -- +(0mm,5mm) -- cycle [line width=0.6pt]}
\def\judy{-- +(5mm,0mm) -- +(5mm,5mm) -- +(0mm,5mm) -- cycle [line width=0.6pt,fill=gainsboro]}
\def\nbx{-- +(5mm,0mm) -- +(5mm,5mm) -- +(0mm,5mm) -- cycle [line width=0.6pt,pattern={Lines[angle=45,distance=2.5pt, line width=0.6pt]}, pattern color=black]}
\tikzstyle{bdot}=[circle,fill=royalblue,draw=royalblue,inner sep=2]
\tikzstyle{rdot}=[circle,fill=red,draw=red,inner sep=2]

\draw (0mm,0mm)\judy;
\draw (5mm,0mm)\hezi;
\draw (10mm,0mm)\hezi;
\draw (15mm,0mm)\hezi;
\draw (20mm,0mm)\hezi;
\draw (0mm,-5mm)\hezi;
\draw (5mm,-5mm)\hezi;
\draw (10mm,-5mm)\judy;
\draw (15mm,-5mm)\hezi;
\draw (20mm,-5mm)\hezi;
\draw (0mm,-10mm)\hezi;
\draw (5mm,-10mm)\judy;
\draw (10mm,-10mm)\hezi;
\draw (15mm,-10mm)\hezi;
\draw (20mm,-10mm)\hezi;
\draw (0mm,-15mm)\hezi;
\draw (5mm,-15mm)\hezi;
\draw (10mm,-15mm)\hezi;
\draw (15mm,-15mm)\hezi;
\draw (20mm,-15mm)\nbx;
\draw (0mm,-20mm)\hezi;
\draw (5mm,-20mm)\hezi;
\draw (10mm,-20mm)\hezi;
\draw (15mm,-20mm)\nbx;
\draw (20mm,-20mm)\hezi;

\draw (10mm,-5mm)-- +(5mm,0mm) -- +(5mm,5mm) -- +(0mm,5mm) -- cycle [line width=1.5pt,draw=royalblue];

\draw (20mm,-15mm)-- +(5mm,0mm) -- +(5mm,5mm) -- +(0mm,5mm) -- cycle [line width=1.5pt,draw=red];

\node at (12.5mm,-26mm){$\sigma=132\overline{54}$};


\begin{scope}[shift={(50mm,-5.5mm)},scale=1.25]
\draw (0mm,0mm)\judy;
\draw (5mm,0mm)\hezi;
\draw (10mm,0mm)\hezi;
\draw (0mm,-5mm)\hezi;
\draw (5mm,-5mm)\hezi;
\draw (10mm,-5mm)\nbx;
\draw (0mm,-10mm)\hezi;
\draw (5mm,-10mm)\nbx;
\draw (10mm,-10mm)\hezi;

\node[bdot] at(5mm,0mm){};

\node at (7.5mm,-16mm){$\eta'^{-1}_{2}(\sigma)$};

\end{scope}


\begin{scope}[shift={(100mm,-5mm)},scale=1.25]
\draw (0mm,0mm)\judy;
\draw (5mm,0mm)\hezi;
\draw (10mm,0mm)\hezi;
\draw (0mm,-5mm)\hezi;
\draw (5mm,-5mm)\hezi;
\draw (10mm,-5mm)\judy;
\draw (0mm,-10mm)\hezi;
\draw (5mm,-10mm)\judy;
\draw (10mm,-10mm)\hezi;

\node[rdot] at(15mm,-10mm){};
\node at (7.5mm,-16mm){$\overline{\eta'}^{-1}_{4}(\sigma)$};

\end{scope}
\end{tikzpicture}
\end{center}
Thus by \eqref{eqpfIn3rd}, \eqref{eqpfIn4th} and \eqref{eqpfIn5th},
we obtain that $\Psi(\pi, (r,|\sigma_r|-1))=(\sigma,r)$.


\noindent{\bf{Case 4: $\bm {|\sigma_r|=r-1}$.}}
Similar to the previous case,
the squares $\langle r,r-1\rangle$ and $\langle r-1,r\rangle$ in the grid $P_\sigma$
are both positive squares if $\sigma_r>0$ and negative squares if $\sigma_r<0$.
From the example below,
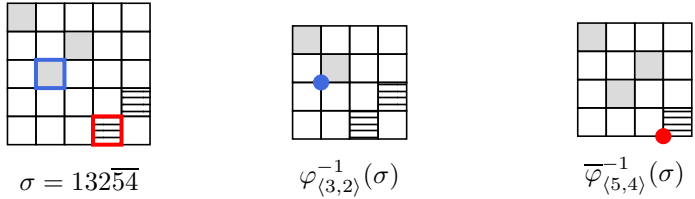
\begin{center}
\begin{tikzpicture}[scale = 0.75]
\def\hezi{-- +(5mm,0mm) -- +(5mm,5mm) -- +(0mm,5mm) -- cycle [line width=0.6pt]}
\def\judy{-- +(5mm,0mm) -- +(5mm,5mm) -- +(0mm,5mm) -- cycle [line width=0.6pt,fill=gainsboro]}
\def\nbx{-- +(5mm,0mm) -- +(5mm,5mm) -- +(0mm,5mm) -- cycle [line width=0.6pt,pattern={Lines[angle=45,distance=2.5pt, line width=0.6pt]}, pattern color=black]}
\tikzstyle{bdot}=[circle,fill=royalblue,draw=royalblue,inner sep=2]
\tikzstyle{rdot}=[circle,fill=red,draw=red,inner sep=2]

\draw (0mm,0mm)\judy;
\draw (5mm,0mm)\hezi;
\draw (10mm,0mm)\hezi;
\draw (15mm,0mm)\hezi;
\draw (20mm,0mm)\hezi;
\draw (0mm,-5mm)\hezi;
\draw (5mm,-5mm)\hezi;
\draw (10mm,-5mm)\judy;
\draw (15mm,-5mm)\hezi;
\draw (20mm,-5mm)\hezi;
\draw (0mm,-10mm)\hezi;
\draw (5mm,-10mm)\judy;
\draw (10mm,-10mm)\hezi;
\draw (15mm,-10mm)\hezi;
\draw (20mm,-10mm)\hezi;
\draw (0mm,-15mm)\hezi;
\draw (5mm,-15mm)\hezi;
\draw (10mm,-15mm)\hezi;
\draw (15mm,-15mm)\hezi;
\draw (20mm,-15mm)\nbx;
\draw (0mm,-20mm)\hezi;
\draw (5mm,-20mm)\hezi;
\draw (10mm,-20mm)\hezi;
\draw (15mm,-20mm)\nbx;
\draw (20mm,-20mm)\hezi;

\draw (5mm,-10mm)-- +(5mm,0mm) -- +(5mm,5mm) -- +(0mm,5mm) -- cycle [line width=1.5pt,draw=royalblue];

\draw (15mm,-20mm)-- +(5mm,0mm) -- +(5mm,5mm) -- +(0mm,5mm) -- cycle [line width=1.5pt,draw=red];

\node at (12.5mm,-26mm){$\sigma=132\overline{54}$};


\begin{scope}[shift={(50mm,-4mm)},scale=1]
\draw (0mm,0mm)\judy;
\draw (5mm,0mm)\hezi;
\draw (10mm,0mm)\hezi;
\draw (15mm,0mm)\hezi;
\draw (0mm,-5mm)\hezi;
\draw (5mm,-5mm)\judy;
\draw (10mm,-5mm)\hezi;
\draw (15mm,-5mm)\hezi;
\draw (0mm,-10mm)\hezi;
\draw (5mm,-10mm)\hezi;
\draw (10mm,-10mm)\hezi;
\draw (15mm,-10mm)\nbx;
\draw (0mm,-15mm)\hezi;
\draw (5mm,-15mm)\hezi;
\draw (10mm,-15mm)\nbx;
\draw (15mm,-15mm)\hezi;

\node[bdot] at(5mm,-5mm){};

\node at (10mm,-22mm){$\varphi^{-1}_{\langle 3,2 \rangle}(\sigma)$};

\end{scope}


\begin{scope}[shift={(100mm,-3.5mm)},scale=1]
\draw (0mm,0mm)\judy;
\draw (5mm,0mm)\hezi;
\draw (10mm,0mm)\hezi;
\draw (15mm,0mm)\hezi;
\draw (0mm,-5mm)\hezi;
\draw (5mm,-5mm)\hezi;
\draw (10mm,-5mm)\judy;
\draw (15mm,-5mm)\hezi;
\draw (0mm,-10mm)\hezi;
\draw (5mm,-10mm)\judy;
\draw (10mm,-10mm)\hezi;
\draw (15mm,-10mm)\hezi;
\draw (0mm,-15mm)\hezi;
\draw (5mm,-15mm)\hezi;
\draw (10mm,-15mm)\hezi;
\draw (15mm,-15mm)\nbx;

\node[rdot] at(15mm,-15mm){};
\node at (10mm,-22mm){$\overline{\varphi}^{-1}_{\langle 5,4\rangle}(\sigma)$};

\end{scope}
\end{tikzpicture}
\end{center}
by Proposition \ref{lemi+1i}, we deduce that for $\sigma_r>0$,
$$
\pi=\varphi^{-1}_{\langle r,\sigma_r\rangle}(\sigma)\in\mathcal{I}_{n-1,k-1}^\ob
$$
since $d^+(r,\sigma_r)=(1,1)$ in $P_\pi$, and for $\sigma_r<0$,
$$
\pi=\overline{\varphi}^{-1}_{\langle r,|\sigma_r|\rangle}(\pi)\in\mathcal{I}_{n-1,k}^\ob
$$
since $d^-(r,\sigma_r)=(0,0)$ in $P_\pi$.
And one can be aware of that the grid point $(r,|\sigma_r|)$ must be the point
where a $1^+_h$-path (resp. $0^-_v$-path) first touches
a positive (resp. negative) square $\langle r-1,|\sigma_r|\rangle$ on the main diagonal
of the grid $P_\pi$ if $\sigma_r>0$ (resp. $\sigma_r<0$).
Hence we see $\Psi(\pi, (r,|\sigma_r|))=(\sigma,r)$ by \eqref{eqpfIn1st} and \eqref{eqpfIn2nd}.


\noindent{\bf{Case 5: $\bm {\sigma_r=r}$ for $\bm {r=1}$, \emph{or} $\bm{r\geq 2}$ with $\bm{\sigma_r\neq r-1}$.}}
Deleting the square $\langle r,\sigma_r\rangle$ from the main diagonal of $P_\sigma$,
we have
\[
\pi=\varphi^{-1}_{\langle r,\sigma_r\rangle}(\sigma)\in \mathcal{I}_{n-1,k-p}^\ob
\]
if the grid point $(r,\sigma_r)$ is of $d^+$-type $(p,p)$ in $P_\pi$ for $p\in\{0,1\}$.
By theorem \ref{lempath},
for $p=1$, the point $(r,\sigma_r)$ is the point where a $1^+_h$-path first touches the main diagonal of $P_\pi$,
and for $p=0$,
the point $(r,\sigma_r)$ is the point
where a $0^+_h$-path first touches the main diagonal of $P_\pi$
since $r=1$ or $r\geq 2$ with $\sigma_{r-1}\neq r-1$ implies that $\langle r-1,r-1\rangle$ is not a positive square in $P_\pi$.
The following diagrams present an example for $\sigma=1\overline{32}45$.
\begin{center}
\begin{tikzpicture}[scale = 0.75]
\def\hezi{-- +(5mm,0mm) -- +(5mm,5mm) -- +(0mm,5mm) -- cycle [line width=0.6pt]}
\def\judy{-- +(5mm,0mm) -- +(5mm,5mm) -- +(0mm,5mm) -- cycle [line width=0.6pt,fill=gainsboro]}
\def\nbx{-- +(5mm,0mm) -- +(5mm,5mm) -- +(0mm,5mm) -- cycle [line width=0.6pt,pattern={Lines[angle=45,distance=2.5pt, line width=0.6pt]}, pattern color=black]}
\tikzstyle{bdot}=[circle,fill=royalblue,draw=royalblue,inner sep=2]
\tikzstyle{rdot}=[circle,fill=babyblue,draw=babyblue,inner sep=2]

\draw (0mm,0mm)\judy;
\draw (5mm,0mm)\hezi;
\draw (10mm,0mm)\hezi;
\draw (15mm,0mm)\hezi;
\draw (20mm,0mm)\hezi;
\draw (0mm,-5mm)\hezi;
\draw (5mm,-5mm)\hezi;
\draw (10mm,-5mm)\nbx;
\draw (15mm,-5mm)\hezi;
\draw (20mm,-5mm)\hezi;
\draw (0mm,-10mm)\hezi;
\draw (5mm,-10mm)\nbx;
\draw (10mm,-10mm)\hezi;
\draw (15mm,-10mm)\hezi;
\draw (20mm,-10mm)\hezi;
\draw (0mm,-15mm)\hezi;
\draw (5mm,-15mm)\hezi;
\draw (10mm,-15mm)\hezi;
\draw (15mm,-15mm)\judy;
\draw (20mm,-15mm)\hezi;
\draw (0mm,-20mm)\hezi;
\draw (5mm,-20mm)\hezi;
\draw (10mm,-20mm)\hezi;
\draw (15mm,-20mm)\hezi;
\draw (20mm,-20mm)\judy;

\draw (0mm,0mm)-- +(5mm,0mm) -- +(5mm,5mm) -- +(0mm,5mm) -- cycle [line width=1.5pt,draw=royalblue];

\draw (15mm,-15mm)-- +(5mm,0mm) -- +(5mm,5mm) -- +(0mm,5mm) -- cycle [line width=1.5pt,draw=babyblue];

\node at (12.5mm,-26mm){$\sigma=1\overline{32}45$};


\begin{scope}[shift={(50mm,-4mm)},scale=1]
\draw (0mm,0mm)\hezi;
\draw (5mm,0mm)\nbx;
\draw (10mm,0mm)\hezi;
\draw (15mm,0mm)\hezi;
\draw (0mm,-5mm)\nbx;
\draw (5mm,-5mm)\hezi;
\draw (10mm,-5mm)\hezi;
\draw (15mm,-5mm)\hezi;
\draw (0mm,-10mm)\hezi;
\draw (5mm,-10mm)\hezi;
\draw (10mm,-10mm)\judy;
\draw (15mm,-10mm)\hezi;
\draw (0mm,-15mm)\hezi;
\draw (5mm,-15mm)\hezi;
\draw (10mm,-15mm)\hezi;
\draw (15mm,-15mm)\judy;

\node[bdot] at(0mm,5mm){};

\node at (10mm,-22mm){$\varphi^{-1}_{\langle 1,1 \rangle}(\sigma)$};

\end{scope}


\begin{scope}[shift={(100mm,-3.5mm)},scale=1]
\draw (0mm,0mm)\judy;
\draw (5mm,0mm)\hezi;
\draw (10mm,0mm)\hezi;
\draw (15mm,0mm)\hezi;
\draw (0mm,-5mm)\hezi;
\draw (5mm,-5mm)\hezi;
\draw (10mm,-5mm)\nbx;
\draw (15mm,-5mm)\hezi;
\draw (0mm,-10mm)\hezi;
\draw (5mm,-10mm)\nbx;
\draw (10mm,-10mm)\hezi;
\draw (15mm,-10mm)\hezi;
\draw (0mm,-15mm)\hezi;
\draw (5mm,-15mm)\hezi;
\draw (10mm,-15mm)\hezi;
\draw (15mm,-15mm)\judy;

\node[rdot] at(15mm,-10mm){};
\node at (10mm,-22mm){$\varphi^{-1}_{\langle 4,4\rangle}(\sigma)$};

\end{scope}
\end{tikzpicture}
\end{center}
Therefore by \eqref{eqpfIn1st} and \eqref{eqpfIn2nd}, we derive that
$\Psi(\pi, (r,\sigma_r))=(\sigma,r)$.

\noindent{\bf{Case 6: $\bm {\sigma_r=\overline{r}}$ for $\bm {r=1}$, \emph{or} $\bm{r\geq 2}$ with $\bm{\sigma_r\neq \overline{r-1}}$.}}
We just need to replace all positive objects in Case 5, including squares, paths, $d$-types, and deleting operations, with negative ones to conclude that
\[
\pi=\overline{\varphi}^{-1}_{\langle r,\sigma_r\rangle}(\sigma)\in \mathcal{I}_{n-1,k-p}^\ob,
\]
and the grid point $(r,\sigma_r)$ is the first point where a $p^-_h$-path
touches the main diagonal of the grid $P_\pi$ for $p\in\{0,1\}$,
which implies that $\Psi(\pi, (r,|\sigma_r|))=(\sigma,r)$ by \eqref{eqpfIn1st} and \eqref{eqpfIn2nd}.

\noindent{\bf{Case 7: $\bm {\sigma_r=r}$ for $\bm {r\geq 2}$ with $\bm{\sigma_r= r-1}$.}}
Since the squares $\langle r-1,r-1\rangle$ and $\langle r, r\rangle$ compose a $12$-pair on the main diagonal of $P_\sigma$,
we employ the double deleting operation $\eta^{-1}_{r-1}$ to obtain
\begin{equation*}
  \pi=\eta^{-1}_{r-1}(\sigma)\in\mathcal{I}_{n-2,k-p}^\ob,
\end{equation*}
and the grid point $(r-1,\sigma_{r-1})$ in the grid $P_\pi$ is of $d^+$-type $(p,p)$ for $p\in\{0,1\}$ by Proposition \ref{lemetaeta'},
as in the following example for $\sigma=\overline{4}23\overline{1}5$.
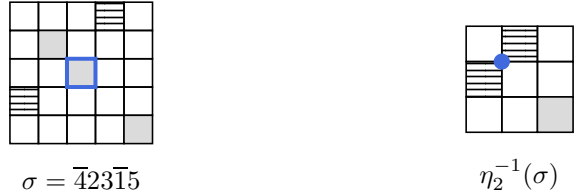
\begin{center}
\begin{tikzpicture}[scale = 0.75]
\def\hezi{-- +(5mm,0mm) -- +(5mm,5mm) -- +(0mm,5mm) -- cycle [line width=0.6pt]}
\def\judy{-- +(5mm,0mm) -- +(5mm,5mm) -- +(0mm,5mm) -- cycle [line width=0.6pt,fill=gainsboro]}
\def\nbx{-- +(5mm,0mm) -- +(5mm,5mm) -- +(0mm,5mm) -- cycle [line width=0.6pt,pattern={Lines[angle=45,distance=2.5pt, line width=0.6pt]}, pattern color=black]}
\tikzstyle{bdot}=[circle,fill=royalblue,draw=royalblue,inner sep=2]
\tikzstyle{rdot}=[circle,fill=babyblue,draw=babyblue,inner sep=2]

\draw (0mm,0mm)\hezi;
\draw (5mm,0mm)\hezi;
\draw (10mm,0mm)\hezi;
\draw (15mm,0mm)\nbx;
\draw (20mm,0mm)\hezi;
\draw (0mm,-5mm)\hezi;
\draw (5mm,-5mm)\judy;
\draw (10mm,-5mm)\hezi;
\draw (15mm,-5mm)\hezi;
\draw (20mm,-5mm)\hezi;
\draw (0mm,-10mm)\hezi;
\draw (5mm,-10mm)\hezi;
\draw (10mm,-10mm)\judy;
\draw (15mm,-10mm)\hezi;
\draw (20mm,-10mm)\hezi;
\draw (0mm,-15mm)\nbx;
\draw (5mm,-15mm)\hezi;
\draw (10mm,-15mm)\hezi;
\draw (15mm,-15mm)\hezi;
\draw (20mm,-15mm)\hezi;
\draw (0mm,-20mm)\hezi;
\draw (5mm,-20mm)\hezi;
\draw (10mm,-20mm)\hezi;
\draw (15mm,-20mm)\hezi;
\draw (20mm,-20mm)\judy;

\draw (10mm,-10mm)-- +(5mm,0mm) -- +(5mm,5mm) -- +(0mm,5mm) -- cycle [line width=1.5pt,draw=royalblue];


\node at (12.5mm,-26mm){$\sigma=\overline{4}23\overline{1}5$};


\begin{scope}[shift={(80mm,-5.5mm)},scale=1.25]
\draw (0mm,0mm)\hezi;
\draw (5mm,0mm)\nbx;
\draw (10mm,0mm)\hezi;
\draw (0mm,-5mm)\nbx;
\draw (5mm,-5mm)\hezi;
\draw (10mm,-5mm)\hezi;
\draw (0mm,-10mm)\hezi;
\draw (5mm,-10mm)\hezi;
\draw (10mm,-10mm)\judy;

\node[bdot] at(5mm,0mm){};

\node at (7.5mm,-16mm){$\eta^{-1}_{2}(\sigma)$};

\end{scope}

\end{tikzpicture}
\end{center}
Hence we have $\Psi(\pi, (r-1,\sigma_{r-1}))=(\sigma, r)$ by \eqref{eqpfIn3rd} and \eqref{eqpfIn4th}.

\noindent{\bf{Case 8: $\bm {\sigma_r=\overline{r}}$ for $\bm {r\geq 2}$ with $\bm{\sigma_r= \overline{r-1}}$.}}
As in Case 6 is to Case 5, substituting all positive objects in Case 7 with the corresponding negative ones to get
\begin{equation*}
  \pi=\overline{\eta}^{-1}_{r-1}(\sigma)\in\mathcal{I}_{n-2, k-1-p}^\ob,
\end{equation*}
where the grid point $(r-1,|\sigma_{r-1}|)$ in the grid $P_\pi$ is of $d^-$-type $(p,p)$ for $p\in\{0,1\}$ by Proposition \ref{lemetaeta'}.
Therefore, we have $\Psi(\pi, (r-1,|\sigma_{r-1}|))=(\sigma, r)$ by \eqref{eqpfIn4th} and \eqref{eqpfIn5th}.

In summary, Case 1 to Case 8 show that the mapping $\Theta$ defined by \eqref{theta} is indeed the inverse of the mapping $\Psi$,
which implies that
\[
\Psi\colon{\mathcal{E}}_{n,k}\leftrightarrow{\mathcal{G}}_{n,k}
\]
does give the desired bijection,
and thus we complete the combinatorial proof.
\qed

\subsection{Combinatorial proof of the recurrence of $J_{n,k}^\ob$}
For $1\leq k\leq n$, recall that $\mathcal{J}_{n,k}^\ob$ is the set of fixed-point free involutions in $\mathcal{J}_n^B$ with $k$ $\ob$-descents, and $|\mathcal{J}_{n,k}^\ob|=J_{n,k}^\ob$.

{\noindent{\emph{Combinatorial Proof of Theorem \ref{thmJobnk}.}\hskip 2pt}}
To emphasis the combinatorial significance of \eqref{recJobnk},
we multiply both sides of it by $2$, and prove
\begin{equation}\label{recJobnkpf}
\begin{aligned}
2n J_{2 n, k}^\ob&=2(k^2+k+n-1) J_{2 n-2, k}^\ob+4((k-1)(2 n-k-1)+n) J_{2 n-2, k-1}^\ob \\[3pt]
&\quad+2((2 n-k)(2 n-k+1)+(n-1)) J_{2 n-2, k-2}^\ob.
\end{aligned}
\end{equation}
Let
\begin{equation*}
  {\mathcal{H}}_{n,k}=
  \{(\pi,i)\mid \pi\in \mathcal{J}_{n,k}^\ob\text{ and }1\leq i\leq n \},
\end{equation*}
then we have $|\mathcal{H}_{2n,k}^\ob|=2nJ_{2n,k}^\ob$,
which is the left side of \eqref{recJobnkpf}.

Let
\begin{equation*}
\mathcal{F}^{(1)^+}_{n,k}=\{(\pi,(i,j))\mid \pi\in \mathcal{J}_{n,k}^\ob\text{ and $d^+(i,j)=(0,0)$ in $P_\pi$ with $i\neq j$}\}
\end{equation*}
and
\begin{equation*}
\mathcal{F}^{(1)^-}_{n,k}=\{(\pi,(i,j))\mid \pi\in \mathcal{J}_{n,k}^\ob\text{ and $d^-(i,j)=(0,0)$ in $P_\pi$ with $i\neq j$}\}.
\end{equation*}
Notice that there is no positive or negative square on the main diagonal of $P_\pi$
since $\pi$ is fixed-point free.
Thus by Theorem \ref{thmpaths}, the number of grid points $(i,i)$ of $d^+$-type $(0,0)$ and $d^-$-type $(0,0)$
equals to the number of $0_h^+$-paths and $0_h^-$-paths in $P_\pi$, respectively.
Let
\begin{equation*}
  \mathcal{F}^{(2)^-}_{n,k}=
  \{(\pi,(i,i))_1, (\pi,(i,i))_2 \mid \pi\in \mathcal{J}_{n,k}^\ob\text{ and $d^-(i,i)=(0,0)$ in $P_\pi$}  \},
\end{equation*}
where the subscripts $1$ and $2$ indicate the first and  second occurrences  of
the pair $(\pi,(i,i))$ in the set $\mathcal{F}^{(2)^-}_{n,k}$.
It follows from Theorems \ref{thmpaths} and \ref{thm-enum-dtype} that
\begin{equation*}
\begin{aligned}
&\left|\mathcal{F}^{(1)^+}_{2n-2,k}\uplus\mathcal{F}^{(1)^-}_{2n-2,k}\uplus\mathcal{F}^{(2)^-}_{2n-2,k}\right|\\
=&\sum_{\pi\in\mathcal{J}_{2n-2,k}^\ob}\left((k+1)^2-n(\pi)+(2n-2)-(k+1)\right)+\left(k^2+n(\pi)-k\right)+2k\\
=&2(k^2+k+n-1)J_{2n-2,k}^\ob,
\end{aligned}
\end{equation*}
which is the \textbf{first} term in the right side of \eqref{recJobnkpf}.
For any $(\pi, (i,j))\in \mathcal{F}^{(1)^+}_{2n-2,k}\uplus\mathcal{F}^{(1)^-}_{2n-2,k}\uplus\mathcal{F}^{(2)^-}_{2n-2,k}$,
define
\begin{equation*}
\Psi((\pi, (i,j)))=(\sigma,r)
\end{equation*}
by setting
\begin{equation}\label{eqpfJn1st}
(\sigma,r)=\left\{
      \begin{array}{ll}
      (\xi_{(i,j)}(\pi), \chi_{ij}), & \hbox{if $(\pi,(i,j))\in\mathcal{F}^{(1)^+}_{2n-2,k}$ ,} \\[6pt]
      (\overline{\xi}_{(i,j)}(\pi), \chi_{ij}), & \hbox{if $(\pi,(i,j))\in\mathcal{F}^{(1)^-}_{2n-2,k}$,}\\[6pt]
      (\overline{\eta'}_{i}(\pi),i+1), & \hbox{if $(\pi,(i,j))_1\in\mathcal{F}^{(2)^-}_{2n-2,k}$,}\\[6pt]
      (\overline{\eta'}_{i}(\pi),i), & \hbox{if $(\pi,(i,j))_2\in\mathcal{F}^{(2)^-}_{2n-2,k}$.}
      \end{array}
    \right.
\end{equation}

Let
\begin{equation*}
\mathcal{F}^{(5)^+}_{n,k}=\{(\pi,(i,j))\mid\text{ $\pi\in\mathcal{J}_{n,k}^\ob$ and $d^+(i,j)=(1,0)$ or $(0,1)$ in $P_\pi$}\}
\end{equation*}
and
\begin{equation*}
\mathcal{F}^{(5)^-}_{n,k}=\{(\pi,(i,j))\mid\text{ $\pi\in\mathcal{J}_{n,k}^\ob$ and $d^-(i,j)=(1,0)$ or $(0,1)$ in $P_\pi$}\}.
\end{equation*}
Clearly, by Proposition \ref{lemij=ji}, any grid point in the pairs from the above two sets is not on the main diagonal of $P_\pi$.
And by Theorem \ref{thm-enum-dtype}, we deduce that
\begin{equation*}
\begin{aligned}
&\left|\mathcal{F}^{(5)^+}_{2n-2,k-1}\uplus\mathcal{F}^{(5)^-}_{2n-2,k-1}\right|\\
=&\sum_{\pi\in\mathcal{J}_{2n-2,k-1}^\ob}2\left(k(2n-k-1)+n(\pi)-(2n-2)\right)
+2\left((k-1)(2n-k)-n(\pi)\right)\\
=&\sum_{\pi\in\mathcal{J}_{2n-2,k-1}^\ob}2k(2n-k-1)+2(k-1)(2n-k)-2(2n-2)\\
=&4(k-1)(2n-k-1)J_{2n-2,k-1}^\ob.
\end{aligned}
\end{equation*}
In addition, let
\begin{equation*}
\mathcal{F}^{(2)^+}_{n,k}=\{(\pi,(i,i))_1, (\pi,(i,i))_2 \mid \pi\in \mathcal{J}_{n,k}^\ob\text{ and $d^+(i,i)=(0,0)$ in $P_\pi$}  \}
\end{equation*}
and
\begin{equation*}
\mathcal{F}^{(4)^-}_{n,k}=\{(\pi,(i,i))_1, (\pi,(i,i))_2 \mid \pi\in \mathcal{J}_{n,k}^\ob\text{ and $d^-(i,i)=(1,1)$ in $P_\pi$}  \},
\end{equation*}
then from Theorem \ref{thmpaths}, we have
\begin{equation*}
\begin{aligned}
\left|\mathcal{F}^{(2)^+}_{2n-2,k-1}\uplus\mathcal{F}^{(4)^-}_{2n-2,k-1}\right|
=&\sum_{\pi\in\mathcal{J}_{2n-2,k-1}^\ob}2k
+2\left(2n-2-(k-1)+1\right)\\
=&(2k+2(2n-k))J_{2n-2,k-1}^\ob.
\end{aligned}
\end{equation*}
Thus,
\begin{equation*}
\begin{aligned}
&\left|\mathcal{F}^{(5)^+}_{2n-2,k-1}\uplus\mathcal{F}^{(5)^-}_{2n-2,k-1}\uplus\mathcal{F}^{(2)^+}_{2n-2,k-1}\uplus\mathcal{F}^{(4)^-}_{2n-2,k-1}\right|\\
=&4(k-1)(2n-k-1)J_{2n-2,k-1}^\ob+(2k+2(2n-k))J_{2n-2,k-1}^\ob\\
=&4((k-1)(2n-k-1)+n)J_{2n-2,k-1}^\ob,
\end{aligned}
\end{equation*}
which is the \textbf{second} term in the right side of \eqref{recJobnkpf}.
For any $(\pi, (i,j))\in \mathcal{F}^{(5)^+}_{2n-2,k-1}\uplus\mathcal{F}^{(5)^-}_{2n-2,k-1}\uplus\mathcal{F}^{(2)^+}_{2n-2,k-1}\uplus\mathcal{F}^{(4)^-}_{2n-2,k-1}$,
define
\begin{equation*}
\Psi((\pi, (i,j)))=(\sigma,r),
\end{equation*}
where
\begin{equation}\label{eqpfJn2nd}
(\sigma,r)=\left\{
      \begin{array}{ll}
      (\xi_{(i,j)}(\pi),\chi_{ij}), & \hbox{if $(\pi,(i,j))\in\mathcal{F}^{(5)^+}_{2n-2,k-1}$ ,} \\[6pt]
      (\overline{\xi}_{(i,j)}(\pi),\chi_{ij}), & \hbox{if $(\pi,(i,j))\in\mathcal{F}^{(5)^-}_{2n-2,k-1}$,}\\[6pt]
      (\eta'_{i}(\pi),i+1), & \hbox{if $(\pi,(i,j))_1\in\mathcal{F}^{(2)^+}_{2n-2,k-1}$,}\\[6pt]
      (\eta'_{i}(\pi),i), & \hbox{if $(\pi,(i,j))_2\in\mathcal{F}^{(2)^+}_{2n-2,k-1}$,}\\[6pt]
      (\overline{\eta'}_{i}(\pi),i+1), & \hbox{if $(\pi,(i,j))_1\in\mathcal{F}^{(4)^-}_{2n-2,k-1}$,}\\[6pt]
      (\overline{\eta'}_{i}(\pi),i), & \hbox{if $(\pi,(i,j))_2\in\mathcal{F}^{(4)^-}_{2n-2,k-1}$.}
      \end{array}
    \right.
\end{equation}

Let
\begin{equation*}
\mathcal{F}^{(3)^+}_{n,k}=\{(\pi,(i,j))\mid \pi\in \mathcal{J}_{n,k}^\ob\text{ and $d^+(i,j)=(1,1)$ in $P_\pi$ with $i\neq j$}\},
\end{equation*}
\begin{equation*}
\mathcal{F}^{(3)^-}_{n,k}=\{(\pi,(i,j))\mid \pi\in \mathcal{J}_{n,k}^\ob\text{ and $d^-(i,j)=(1,1)$ in $P_\pi$ with $i\neq j$}\},
\end{equation*}
and
\begin{equation*}
\mathcal{F}^{(4)^+}_{n,k}=\{(\pi,(i,i))_1, (\pi,(i,i))_2 \mid \pi\in \mathcal{J}_{n,k}^\ob\text{ and $d^+(i,i)=(1,1)$ in $P_\pi$}  \}.
\end{equation*}
By Theorems \ref{thmpaths} and \ref{thm-enum-dtype}, we obtain
\begin{equation*}
\begin{aligned}
&\left|\mathcal{F}^{(3)^+}_{2n-2,k-2}\uplus\mathcal{F}^{(3)^-}_{2n-2,k-2}\uplus\mathcal{F}^{(4)^+}_{2n-2,k-2}\right|\\
=&\sum_{\pi\in\mathcal{J}_{2n-2,k-2}^\ob}\left((2n-k)(2n-k-1)-n(\pi)+(2n-2)\right)+
\left((2n-k+1)(2n-k)+n(\pi)\right)+2(2n-k)\\
=&2((2n-k)(2n-k+1)+n-1)J_{2n-2,k-2}^\ob,
\end{aligned}
\end{equation*}
which is the \textbf{third} term in the right side of \eqref{recJobnkpf}.
For any $(\pi, (i,j))\in \mathcal{F}^{(3)^+}_{2n-2,k-2}\uplus\mathcal{F}^{(3)^-}_{2n-2,k-2}\uplus
\mathcal{F}^{(4)^+}_{2n-2,k-2}$,
we set
\begin{equation*}
\Psi((\pi, (i,j)))=(\sigma,r),
\end{equation*}
where
\begin{equation}\label{eqpfJn3rd}
(\sigma,r)=\left\{
      \begin{array}{ll}
      (\xi_{(i,j)}(\pi),\chi_{ij}), & \hbox{if $(\pi,(i,j))\in\mathcal{F}^{(3)^+}_{2n-2,k-2}$ ,} \\[6pt]
      (\overline{\xi}_{(i,j)}(\pi),\chi_{ij}), & \hbox{if $(\pi,(i,j))\in\mathcal{F}^{(3)^-}_{2n-2,k-2}$,}\\[6pt]
      (\eta'_{i}(\pi),i+1), & \hbox{if $(\pi,(i,j))_1\in\mathcal{F}^{(4)^+}_{2n-2,k-2}$,}\\[6pt]
      (\eta'_{i}(\pi),i) & \hbox{if $(\pi,(i,j))_2\in\mathcal{F}^{(4)^+}_{2n-2,k-2}$.}
      \end{array}
    \right.
\end{equation}

Thus, the cardinality of the set
\begin{equation*}
{\mathcal{F}}_{n,k}=
\mathcal{F}^{(1)^\pm}_{n-2,k}\uplus\mathcal{F}^{(2)+}_{n-2,k-1}\uplus\mathcal{F}^{(2)-}_{n-2,k}
\uplus\mathcal{F}^{(3)^\pm}_{n-2,k-2}\uplus\mathcal{F}^{(4)+}_{n-2,k-2}\uplus\mathcal{F}^{(4)-}_{n-2,k-1}
\uplus\mathcal{F}^{(5)^\pm}_{n-2,k-1}
\end{equation*}
equals to the right side of \eqref{recJobnkpf}.
It can be directly checked that by Propositions \ref{lemxi} and \ref{lemetaeta'},
we establish the mapping
\begin{equation}\label{PsiJnkob}
\Psi: \mathcal{F}_{2n,k}\rightarrow \mathcal{H}_{2n,k}
\end{equation}
in terms of \eqref{eqpfJn1st}--\eqref{eqpfJn3rd}.

To complete the proof, for any chosen $(\sigma,r)\in {\mathcal{H}}_{2n,k}$,
we show that the mapping
\begin{equation*}
\Theta((\sigma, r))=(\pi, (i,j))
\end{equation*}
with
\begin{equation}\label{Jtheta}
(\pi, (i,j))=\left\{
        \begin{array}{ll}
          (\xi^{-1}_{\langle r,\sigma_r-1\rangle}(\sigma),(r,\sigma_r-1)), & r+1< \sigma_r\leq 2n, \\[8pt]
          (\xi^{-1}_{\langle r-1,\sigma_r\rangle}(\sigma),(r-1,\sigma_r)), & 1\leq \sigma_r< r-1, \\[8pt]
          (\eta'^{-1}_{r-1}(\sigma),(r-1,\sigma_r))_1, & \sigma_r=r-1,\\[8pt]
          (\eta'^{-1}_{r}(\sigma),(r,\sigma_r-1))_2, & \sigma_r=r+1, \\[8pt]
          (\overline{\xi}^{-1}_{\langle r,|\sigma_r|-1\rangle}(\sigma),(r,|\sigma_r|-1)), & \overline{2n}\leq \sigma_r<\overline{r+1}, \\[8pt]
          (\overline{\xi}^{-1}_{\langle r-1,|\sigma_r|\rangle}(\sigma),(r-1,|\sigma_r|)), & \overline{r-1}<\sigma_r\leq \overline{1}, \\[8pt]
          (\overline{\eta}'^{-1}_{r-1}(\sigma),(r-1,|\sigma_r|))_1, & \sigma_r=\overline{r-1}, \\[8pt]
          (\overline{\eta}'^{-1}_{r}(\sigma),(r,|\sigma_r|-1))_2, & \sigma_r=\overline{r+1}
        \end{array}
      \right.
\end{equation}
is the inverse of $\Psi$.

For $r+1<|\sigma_r |\leq 2n$ and $1<|\sigma_r |\leq r-1$, to conclude $(\pi,(i,j))\in\mathcal{F}_{2n,k}$,
one can refer to Case $1$ and Case $2$ in the proof of Theorem \ref{thmIobnk}.
If $\sigma_r=r-1$, the positive squares $\langle r-1, \sigma_{r-1}\rangle$ and $\langle r,\sigma_r\rangle$ form
a $21$-pair on the main diagonal of $P_\sigma$,
and the square  $\langle r,\sigma_r\rangle$ plays the role of $1$ in the $21$-pair.
By noticing that only pairs $(\pi, (i,i))$ in the set $\mathcal{F}_{n,k}^{(l)^\ast}$ for $l\in\{2,4\}$ and $\ast\in\{+,-\}$ appear twice,
we can think of the grid point $(r-1,r-1)$ as the first occurrence.
Thus we see
$$(\eta'^{-1}_{r-1}(\sigma),(r-1,\sigma_r))_1\in \mathcal{F}^{(4)^+}_{2n-2,k-2}$$
if the grid point $(r-1,r-1)$ is of $d$-type $(1,1)$,
and
$$(\eta'^{-1}_{r-1}(\sigma),(r-1,\sigma_r))_1\in \mathcal{F}^{(2)^+}_{2n-2,k-1}$$
if the grid point $(r-1,r-1)$ is of $d$-type $(0,0)$.
If $\sigma_r=r+1$, then on the main diagonal in $P_\sigma$,
there must exist a $21$-pair
consists of the positive squares $\langle r,\sigma_r \rangle$ and $\langle r+1, \sigma_{r+1} \rangle$.
Since the square $\langle r,\sigma_r\rangle$ plays the part of $2$ in the $21$-pair,
we designate the grid point $(r,r)$ as the second occurrence.
Hence we derive that
$$(\eta'^{-1}_r(\sigma),(r,\sigma_r-1))_2\in \mathcal{F}^{(4)^+}_{2n-2,k-2}$$
if the grid point $(r,r)$ is of $d$-type $(1,1)$,
and
$$(\eta'^{-1}_r(\sigma),(r,\sigma_r-1))_2\in \mathcal{F}^{(2)^+}_{2n-2,k-1}$$
if the grid point $(r,r)$ is of $d$-type $(0,0)$.
The following example of $\sigma=\overline{5}643\overline{1}2$
helps to better understand,
where the dark blue point in the grid of $\eta'^{-1}_3(\sigma)$ represents the first occurrence of $(3,3)$,
while the light blue one indicates the second occurrence.
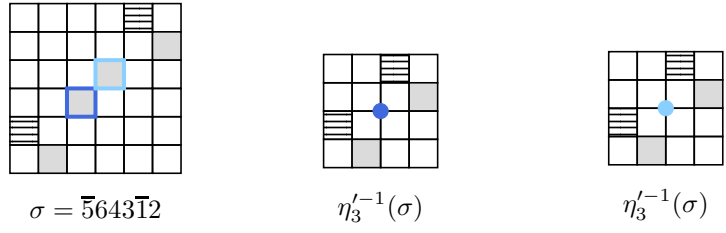
\begin{center}
\begin{tikzpicture}[scale = 0.75]
\def\hezi{-- +(5mm,0mm) -- +(5mm,5mm) -- +(0mm,5mm) -- cycle [line width=0.6pt]}
\def\judy{-- +(5mm,0mm) -- +(5mm,5mm) -- +(0mm,5mm) -- cycle [line width=0.6pt,fill=gainsboro]}
\def\nbx{-- +(5mm,0mm) -- +(5mm,5mm) -- +(0mm,5mm) -- cycle [line width=0.6pt,pattern={Lines[angle=45,distance=2.5pt, line width=0.6pt]}, pattern color=black]}
\tikzstyle{bdot}=[circle,fill=royalblue,draw=royalblue,inner sep=2]
\tikzstyle{rdot}=[circle,fill=babyblue,draw=babyblue,inner sep=2]

\draw (0mm,5mm)\hezi;
\draw (5mm,5mm)\hezi;
\draw (10mm,5mm)\hezi;
\draw (15mm,5mm)\hezi;
\draw (20mm,5mm)\nbx;
\draw (25mm,5mm)\hezi;

\draw (0mm,0mm)\hezi;
\draw (5mm,0mm)\hezi;
\draw (10mm,0mm)\hezi;
\draw (15mm,0mm)\hezi;
\draw (20mm,0mm)\hezi;
\draw (25mm,0mm)\judy;
\draw (0mm,-5mm)\hezi;
\draw (5mm,-5mm)\hezi;
\draw (10mm,-5mm)\hezi;
\draw (15mm,-5mm)\judy;
\draw (20mm,-5mm)\hezi;
\draw (25mm,-5mm)\hezi;
\draw (0mm,-10mm)\hezi;
\draw (5mm,-10mm)\hezi;
\draw (10mm,-10mm)\judy;
\draw (15mm,-10mm)\hezi;
\draw (20mm,-10mm)\hezi;
\draw (25mm,-10mm)\hezi;
\draw (0mm,-15mm)\nbx;
\draw (5mm,-15mm)\hezi;
\draw (10mm,-15mm)\hezi;
\draw (15mm,-15mm)\hezi;
\draw (20mm,-15mm)\hezi;
\draw (25mm,-15mm)\hezi;
\draw (0mm,-20mm)\hezi;
\draw (5mm,-20mm)\judy;
\draw (10mm,-20mm)\hezi;
\draw (15mm,-20mm)\hezi;
\draw (20mm,-20mm)\hezi;
\draw (25mm,-20mm)\hezi;

\draw (10mm,-10mm)-- +(5mm,0mm) -- +(5mm,5mm) -- +(0mm,5mm) -- cycle [line width=1.5pt,draw=royalblue];

\draw (15mm,-5mm)-- +(5mm,0mm) -- +(5mm,5mm) -- +(0mm,5mm) -- cycle [line width=1.5pt,draw=babyblue];

\node at (15mm,-26mm){$\sigma=\overline{5}643\overline{1}2$};


\begin{scope}[shift={(55mm,-4mm)},scale=1]
\draw (0mm,0mm)\hezi;
\draw (5mm,0mm)\hezi;
\draw (10mm,0mm)\nbx;
\draw (15mm,0mm)\hezi;
\draw (0mm,-5mm)\hezi;
\draw (5mm,-5mm)\hezi;
\draw (10mm,-5mm)\hezi;
\draw (15mm,-5mm)\judy;
\draw (0mm,-10mm)\nbx;
\draw (5mm,-10mm)\hezi;
\draw (10mm,-10mm)\hezi;
\draw (15mm,-10mm)\hezi;
\draw (0mm,-15mm)\hezi;
\draw (5mm,-15mm)\judy;
\draw (10mm,-15mm)\hezi;
\draw (15mm,-15mm)\hezi;

\node[bdot] at(10mm,-5mm){};

\node at (10mm,-22mm){$\eta'^{-1}_{3}(\sigma)$};

\end{scope}


\begin{scope}[shift={(105mm,-3.5mm)},scale=1]
\draw (0mm,0mm)\hezi;
\draw (5mm,0mm)\hezi;
\draw (10mm,0mm)\nbx;
\draw (15mm,0mm)\hezi;
\draw (0mm,-5mm)\hezi;
\draw (5mm,-5mm)\hezi;
\draw (10mm,-5mm)\hezi;
\draw (15mm,-5mm)\judy;
\draw (0mm,-10mm)\nbx;
\draw (5mm,-10mm)\hezi;
\draw (10mm,-10mm)\hezi;
\draw (15mm,-10mm)\hezi;
\draw (0mm,-15mm)\hezi;
\draw (5mm,-15mm)\judy;
\draw (10mm,-15mm)\hezi;
\draw (15mm,-15mm)\hezi;

\node[rdot] at(10mm,-5mm){};
\node at (10mm,-22mm){$\eta'^{-1}_{3}(\sigma)$};

\end{scope}
\end{tikzpicture}
\end{center}
For $\sigma_r=\overline{r-1}$ and $\sigma_r=\overline{r+1}$,
the proof is similar,
except that the positive squares, pairs, double deleting operations, and $d$-type are replaced with the corresponding negative ones.

In conclusion, we establish the bijection
\[
\Psi: \mathcal{F}_{2n,k}\leftrightarrow \mathcal{H}_{2n,k}
\]
by combining \eqref{PsiJnkob} and \eqref{Jtheta},
which verifies \eqref{recJobnkpf} and completes the proof.
\qed

\section{Geometric properties of $\ub$-descents and inverse $\ub$-descents}\label{Sec-geounder}

In the subsequent contexts, the definitions and the terminologies
is under the $r$-order, except for special mentions.
Since techniques used in the proof are similar to the previous sections,
we shall omit most of the detailed process
and instead give the corresponding conclusion directly.

We still start by investigating the des-pairs, asc-pairs, and $d$-types.
Comparing to the natural order \eqref{norder},
we note that the only difference is the order between the two negative integers.
Thus, the following are the eight types of des-pairs and asc-pairs under the $r$-order.
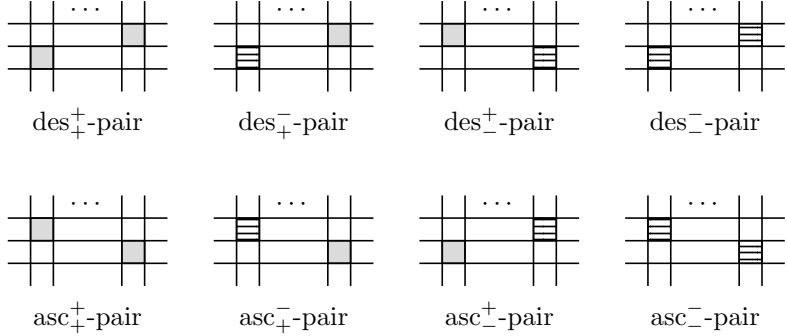
\begin{center}
\begin{tikzpicture}[scale = 0.6]
\def\hezi{-- +(5mm,0mm) -- +(5mm,5mm) -- +(0mm,5mm) -- cycle [line width=0.6pt]}
\def\judy{-- +(5mm,0mm) -- +(5mm,5mm) -- +(0mm,5mm) -- cycle [line width=0.6pt,fill=gainsboro]}
\tikzstyle{bdot}=[circle,fill=black,draw=black,inner sep=1.2]
\tikzstyle{gdot}=[circle,fill=royalblue,draw=royalblue,inner sep=1.5]
\draw (5mm,-5mm)\judy;
\draw (25mm,0mm)\judy;
\draw[line width=0.6pt] (0mm,5mm)--(35mm,5mm);
\draw[line width=0.6pt] (0mm,0mm)--(35mm,0mm);
\draw[line width=0.6pt] (0mm,-5mm)--(35mm,-5mm);
\draw[line width=0.6pt] (5mm,-10mm)--(5mm,10mm);
\draw[line width=0.6pt] (10mm,-10mm)--(10mm,10mm);
\draw[line width=0.6pt] (25mm,-10mm)--(25mm,10mm);
\draw[line width=0.6pt] (30mm,-10mm)--(30mm,10mm);
\node at (17.5mm,8mm){$\cdots$};
\node at (17.5mm,-17mm){des$^+_+$-pair};
\end{tikzpicture}
\hspace{1em}
\begin{tikzpicture}[scale = 0.6]
\def\hezi{-- +(5mm,0mm) -- +(5mm,5mm) -- +(0mm,5mm) -- cycle [line width=0.6pt]}
\def\judy{-- +(5mm,0mm) -- +(5mm,5mm) -- +(0mm,5mm) -- cycle [line width=0.6pt,fill=gainsboro]}
\def\nbx{-- +(5mm,0mm) -- +(5mm,5mm) -- +(0mm,5mm) -- cycle [line width=0.6pt,pattern={Lines[angle=45,distance=2.5pt, line width=0.6pt]}, pattern color=black]}
\tikzstyle{bdot}=[circle,fill=black,draw=black,inner sep=1.2]
\tikzstyle{gdot}=[circle,fill=royalblue,draw=royalblue,inner sep=1.5]
\draw (5mm,-5mm)\nbx;
\draw (25mm,0mm)\judy;
\draw[line width=0.6pt] (0mm,5mm)--(35mm,5mm);
\draw[line width=0.6pt] (0mm,0mm)--(35mm,0mm);
\draw[line width=0.6pt] (0mm,-5mm)--(35mm,-5mm);
\draw[line width=0.6pt] (5mm,-10mm)--(5mm,10mm);
\draw[line width=0.6pt] (10mm,-10mm)--(10mm,10mm);
\draw[line width=0.6pt] (25mm,-10mm)--(25mm,10mm);
\draw[line width=0.6pt] (30mm,-10mm)--(30mm,10mm);
\node at (17.5mm,8mm){$\cdots$};
\node at (17.5mm,-17mm){des$^-_+$-pair};
\end{tikzpicture}
\hspace{1em}
\begin{tikzpicture}[scale = 0.6]
\def\hezi{-- +(5mm,0mm) -- +(5mm,5mm) -- +(0mm,5mm) -- cycle [line width=0.6pt]}
\def\judy{-- +(5mm,0mm) -- +(5mm,5mm) -- +(0mm,5mm) -- cycle [line width=0.6pt,fill=gainsboro]}
\def\nbx{-- +(5mm,0mm) -- +(5mm,5mm) -- +(0mm,5mm) -- cycle [line width=0.6pt,pattern={Lines[angle=45,distance=2.5pt, line width=0.6pt]}, pattern color=black]}
\tikzstyle{bdot}=[circle,fill=black,draw=black,inner sep=1.2]
\tikzstyle{gdot}=[circle,fill=royalblue,draw=royalblue,inner sep=1.5]
\draw (5mm,0mm)\judy;
\draw (25mm,-5mm)\nbx;
\draw[line width=0.6pt] (0mm,5mm)--(35mm,5mm);
\draw[line width=0.6pt] (0mm,0mm)--(35mm,0mm);
\draw[line width=0.6pt] (0mm,-5mm)--(35mm,-5mm);
\draw[line width=0.6pt] (5mm,-10mm)--(5mm,10mm);
\draw[line width=0.6pt] (10mm,-10mm)--(10mm,10mm);
\draw[line width=0.6pt] (25mm,-10mm)--(25mm,10mm);
\draw[line width=0.6pt] (30mm,-10mm)--(30mm,10mm);
\node at (17.5mm,8mm){$\cdots$};
\node at (17.5mm,-17mm){des$^+_-$-pair};
\end{tikzpicture}
\hspace{1em}
\begin{tikzpicture}[scale = 0.6]
\def\hezi{-- +(5mm,0mm) -- +(5mm,5mm) -- +(0mm,5mm) -- cycle [line width=0.6pt]}
\def\nbx{-- +(5mm,0mm) -- +(5mm,5mm) -- +(0mm,5mm) -- cycle [line width=0.6pt,pattern={Lines[angle=45,distance=2.5pt, line width=0.6pt]}, pattern color=black]}
\tikzstyle{bdot}=[circle,fill=black,draw=black,inner sep=1.2]
\tikzstyle{gdot}=[circle,fill=royalblue,draw=royalblue,inner sep=1.5]
\draw (5mm,-5mm)\nbx;
\draw (25mm,0mm)\nbx;
\draw[line width=0.6pt] (0mm,5mm)--(35mm,5mm);
\draw[line width=0.6pt] (0mm,0mm)--(35mm,0mm);
\draw[line width=0.6pt] (0mm,-5mm)--(35mm,-5mm);
\draw[line width=0.6pt] (5mm,-10mm)--(5mm,10mm);
\draw[line width=0.6pt] (10mm,-10mm)--(10mm,10mm);
\draw[line width=0.6pt] (25mm,-10mm)--(25mm,10mm);
\draw[line width=0.6pt] (30mm,-10mm)--(30mm,10mm);
\node at (17.5mm,8mm){$\cdots$};
\node at (17.5mm,-17mm){des$^-_-$-pair};
\end{tikzpicture}
\end{center}
\begin{center}
\begin{tikzpicture}[scale = 0.6]
\def\hezi{-- +(5mm,0mm) -- +(5mm,5mm) -- +(0mm,5mm) -- cycle [line width=0.6pt]}
\def\judy{-- +(5mm,0mm) -- +(5mm,5mm) -- +(0mm,5mm) -- cycle [line width=0.6pt,fill=gainsboro]}
\tikzstyle{bdot}=[circle,fill=black,draw=black,inner sep=1.2]
\tikzstyle{gdot}=[circle,fill=royalblue,draw=royalblue,inner sep=1.5]
\draw (5mm,0mm)\judy;
\draw (25mm,-5mm)\judy;
\draw[line width=0.6pt] (0mm,5mm)--(35mm,5mm);
\draw[line width=0.6pt] (0mm,0mm)--(35mm,0mm);
\draw[line width=0.6pt] (0mm,-5mm)--(35mm,-5mm);
\draw[line width=0.6pt] (5mm,-10mm)--(5mm,10mm);
\draw[line width=0.6pt] (10mm,-10mm)--(10mm,10mm);
\draw[line width=0.6pt] (25mm,-10mm)--(25mm,10mm);
\draw[line width=0.6pt] (30mm,-10mm)--(30mm,10mm);
\node at (17.5mm,8mm){$\cdots$};
\node at (17.5mm,-17mm){asc$^+_+$-pair};
\end{tikzpicture}
\hspace{1em}
\begin{tikzpicture}[scale = 0.6]
\def\hezi{-- +(5mm,0mm) -- +(5mm,5mm) -- +(0mm,5mm) -- cycle [line width=0.6pt]}
\def\judy{-- +(5mm,0mm) -- +(5mm,5mm) -- +(0mm,5mm) -- cycle [line width=0.6pt,fill=gainsboro]}
\def\nbx{-- +(5mm,0mm) -- +(5mm,5mm) -- +(0mm,5mm) -- cycle [line width=0.6pt,pattern={Lines[angle=45,distance=2.5pt, line width=0.6pt]}, pattern color=black]}
\tikzstyle{bdot}=[circle,fill=black,draw=black,inner sep=1.2]
\tikzstyle{gdot}=[circle,fill=royalblue,draw=royalblue,inner sep=1.5]
\draw (5mm,0mm)\nbx;
\draw (25mm,-5mm)\judy;
\draw[line width=0.6pt] (0mm,5mm)--(35mm,5mm);
\draw[line width=0.6pt] (0mm,0mm)--(35mm,0mm);
\draw[line width=0.6pt] (0mm,-5mm)--(35mm,-5mm);
\draw[line width=0.6pt] (5mm,-10mm)--(5mm,10mm);
\draw[line width=0.6pt] (10mm,-10mm)--(10mm,10mm);
\draw[line width=0.6pt] (25mm,-10mm)--(25mm,10mm);
\draw[line width=0.6pt] (30mm,-10mm)--(30mm,10mm);
\node at (17.5mm,8mm){$\cdots$};
\node at (17.5mm,-17mm){asc$^-_+$-pair};
\end{tikzpicture}
\hspace{1em}
\begin{tikzpicture}[scale = 0.6]
\def\hezi{-- +(5mm,0mm) -- +(5mm,5mm) -- +(0mm,5mm) -- cycle [line width=0.6pt]}
\def\judy{-- +(5mm,0mm) -- +(5mm,5mm) -- +(0mm,5mm) -- cycle [line width=0.6pt,fill=gainsboro]}
\def\nbx{-- +(5mm,0mm) -- +(5mm,5mm) -- +(0mm,5mm) -- cycle [line width=0.6pt,pattern={Lines[angle=45,distance=2.5pt, line width=0.6pt]}, pattern color=black]}
\tikzstyle{bdot}=[circle,fill=black,draw=black,inner sep=1.2]
\tikzstyle{gdot}=[circle,fill=royalblue,draw=royalblue,inner sep=1.5]
\draw (5mm,-5mm)\judy;
\draw (25mm,0mm)\nbx;
\draw[line width=0.6pt] (0mm,5mm)--(35mm,5mm);
\draw[line width=0.6pt] (0mm,0mm)--(35mm,0mm);
\draw[line width=0.6pt] (0mm,-5mm)--(35mm,-5mm);
\draw[line width=0.6pt] (5mm,-10mm)--(5mm,10mm);
\draw[line width=0.6pt] (10mm,-10mm)--(10mm,10mm);
\draw[line width=0.6pt] (25mm,-10mm)--(25mm,10mm);
\draw[line width=0.6pt] (30mm,-10mm)--(30mm,10mm);
\node at (17.5mm,8mm){$\cdots$};
\node at (17.5mm,-17mm){asc$^+_-$-pair};
\end{tikzpicture}
\hspace{1em}
\begin{tikzpicture}[scale = 0.6]
\def\hezi{-- +(5mm,0mm) -- +(5mm,5mm) -- +(0mm,5mm) -- cycle [line width=0.6pt]}
\def\nbx{-- +(5mm,0mm) -- +(5mm,5mm) -- +(0mm,5mm) -- cycle [line width=0.6pt,pattern={Lines[angle=45,distance=2.5pt, line width=0.6pt]}, pattern color=black]}
\tikzstyle{bdot}=[circle,fill=black,draw=black,inner sep=1.2]
\tikzstyle{gdot}=[circle,fill=royalblue,draw=royalblue,inner sep=1.5]
\draw (5mm,0mm)\nbx;
\draw (25mm,-5mm)\nbx;
\draw[line width=0.6pt] (0mm,5mm)--(35mm,5mm);
\draw[line width=0.6pt] (0mm,0mm)--(35mm,0mm);
\draw[line width=0.6pt] (0mm,-5mm)--(35mm,-5mm);
\draw[line width=0.6pt] (5mm,-10mm)--(5mm,10mm);
\draw[line width=0.6pt] (10mm,-10mm)--(10mm,10mm);
\draw[line width=0.6pt] (25mm,-10mm)--(25mm,10mm);
\draw[line width=0.6pt] (30mm,-10mm)--(30mm,10mm);
\node at (17.5mm,8mm){$\cdots$};
\node at (17.5mm,-17mm){asc$^-_-$-pair};
\end{tikzpicture}
\end{center}
And the following are the eight types of ides-pairs and iasc-pairs.
\begin{center}
\begin{tikzpicture}[scale = 0.6]
\def\hezi{-- +(5mm,0mm) -- +(5mm,5mm) -- +(0mm,5mm) -- cycle [line width=0.6pt]}
\def\judy{-- +(5mm,0mm) -- +(5mm,5mm) -- +(0mm,5mm) -- cycle [line width=0.6pt,fill=gainsboro]}
\tikzstyle{bdot}=[circle,fill=black,draw=black,inner sep=1.2]
\tikzstyle{gdot}=[circle,fill=royalblue,draw=royalblue,inner sep=1.5]
\draw (0mm,0mm)\judy;
\draw (5mm,15mm)\judy;
\draw[line width=0.6pt] (5mm,-5mm)--(5mm,25mm);
\draw[line width=0.6pt] (0mm,-5mm)--(0mm,25mm);
\draw[line width=0.6pt] (10mm,-5mm)--(10mm,25mm);
\draw[line width=0.6pt] (-5mm,20mm)--(15mm,20mm);
\draw[line width=0.6pt] (-5mm,15mm)--(15mm,15mm);
\draw[line width=0.6pt] (-5mm,5mm)--(15mm,5mm);
\draw[line width=0.6pt] (-5mm,0mm)--(15mm,0mm);
\node at (-3mm,11mm){$\vdots$};
\node at (5mm,-15mm){ides$^+_+$-pair};
\end{tikzpicture}
\hspace{1em}
\begin{tikzpicture}[scale = 0.6]
\def\hezi{-- +(5mm,0mm) -- +(5mm,5mm) -- +(0mm,5mm) -- cycle [line width=0.6pt]}
\def\judy{-- +(5mm,0mm) -- +(5mm,5mm) -- +(0mm,5mm) -- cycle [line width=0.6pt,fill=gainsboro]}
\def\nbx{-- +(5mm,0mm) -- +(5mm,5mm) -- +(0mm,5mm) -- cycle [line width=0.6pt,pattern={Lines[angle=45,distance=2.5pt, line width=0.6pt]}, pattern color=black]}
\tikzstyle{bdot}=[circle,fill=black,draw=black,inner sep=1.2]
\tikzstyle{gdot}=[circle,fill=royalblue,draw=royalblue,inner sep=1.5]
\draw (0mm,0mm)\judy;
\draw (5mm,15mm)\nbx;
\draw[line width=0.6pt] (5mm,-5mm)--(5mm,25mm);
\draw[line width=0.6pt] (0mm,-5mm)--(0mm,25mm);
\draw[line width=0.6pt] (10mm,-5mm)--(10mm,25mm);
\draw[line width=0.6pt] (-5mm,20mm)--(15mm,20mm);
\draw[line width=0.6pt] (-5mm,15mm)--(15mm,15mm);
\draw[line width=0.6pt] (-5mm,5mm)--(15mm,5mm);
\draw[line width=0.6pt] (-5mm,0mm)--(15mm,0mm);
\node at (-3mm,11mm){$\vdots$};
\node at (5mm,-15mm){ides$^-_+$-pair};
\end{tikzpicture}
\hspace{1em}
\begin{tikzpicture}[scale = 0.6]
\def\hezi{-- +(5mm,0mm) -- +(5mm,5mm) -- +(0mm,5mm) -- cycle [line width=0.6pt]}
\def\judy{-- +(5mm,0mm) -- +(5mm,5mm) -- +(0mm,5mm) -- cycle [line width=0.6pt,fill=gainsboro]}
\def\nbx{-- +(5mm,0mm) -- +(5mm,5mm) -- +(0mm,5mm) -- cycle [line width=0.6pt,pattern={Lines[angle=45,distance=2.5pt, line width=0.6pt]}, pattern color=black]}
\tikzstyle{bdot}=[circle,fill=black,draw=black,inner sep=1.2]
\tikzstyle{gdot}=[circle,fill=royalblue,draw=royalblue,inner sep=1.5]
\draw (5mm,0mm)\nbx;
\draw (0mm,15mm)\judy;
\draw[line width=0.6pt] (5mm,-5mm)--(5mm,25mm);
\draw[line width=0.6pt] (0mm,-5mm)--(0mm,25mm);
\draw[line width=0.6pt] (10mm,-5mm)--(10mm,25mm);
\draw[line width=0.6pt] (-5mm,20mm)--(15mm,20mm);
\draw[line width=0.6pt] (-5mm,15mm)--(15mm,15mm);
\draw[line width=0.6pt] (-5mm,5mm)--(15mm,5mm);
\draw[line width=0.6pt] (-5mm,0mm)--(15mm,0mm);
\node at (-3mm,11mm){$\vdots$};
\node at (5mm,-15mm){ides$^+_-$-pair};
\end{tikzpicture}
\hspace{1em}
\begin{tikzpicture}[scale = 0.6]
\def\hezi{-- +(5mm,0mm) -- +(5mm,5mm) -- +(0mm,5mm) -- cycle [line width=0.6pt]}
\def\nbx{-- +(5mm,0mm) -- +(5mm,5mm) -- +(0mm,5mm) -- cycle [line width=0.6pt,pattern={Lines[angle=45,distance=2.5pt, line width=0.6pt]}, pattern color=black]}
\tikzstyle{bdot}=[circle,fill=black,draw=black,inner sep=1.2]
\tikzstyle{gdot}=[circle,fill=royalblue,draw=royalblue,inner sep=1.5]
\draw (0mm,0mm)\nbx;
\draw (5mm,15mm)\nbx;
\draw[line width=0.6pt] (5mm,-5mm)--(5mm,25mm);
\draw[line width=0.6pt] (0mm,-5mm)--(0mm,25mm);
\draw[line width=0.6pt] (10mm,-5mm)--(10mm,25mm);
\draw[line width=0.6pt] (-5mm,20mm)--(15mm,20mm);
\draw[line width=0.6pt] (-5mm,15mm)--(15mm,15mm);
\draw[line width=0.6pt] (-5mm,5mm)--(15mm,5mm);
\draw[line width=0.6pt] (-5mm,0mm)--(15mm,0mm);
\node at (-3mm,11mm){$\vdots$};
\node at (5mm,-15mm){ides$^-_-$-pair};
\end{tikzpicture}
\end{center}
\begin{center}
\begin{tikzpicture}[scale = 0.6]
\def\hezi{-- +(5mm,0mm) -- +(5mm,5mm) -- +(0mm,5mm) -- cycle [line width=0.6pt]}
\def\judy{-- +(5mm,0mm) -- +(5mm,5mm) -- +(0mm,5mm) -- cycle [line width=0.6pt,fill=gainsboro]}
\tikzstyle{bdot}=[circle,fill=black,draw=black,inner sep=1.2]
\tikzstyle{gdot}=[circle,fill=royalblue,draw=royalblue,inner sep=1.5]
\draw (5mm,0mm)\judy;
\draw (0mm,15mm)\judy;
\draw[line width=0.6pt] (5mm,-5mm)--(5mm,25mm);
\draw[line width=0.6pt] (0mm,-5mm)--(0mm,25mm);
\draw[line width=0.6pt] (10mm,-5mm)--(10mm,25mm);
\draw[line width=0.6pt] (-5mm,20mm)--(15mm,20mm);
\draw[line width=0.6pt] (-5mm,15mm)--(15mm,15mm);
\draw[line width=0.6pt] (-5mm,5mm)--(15mm,5mm);
\draw[line width=0.6pt] (-5mm,0mm)--(15mm,0mm);
\node at (-3mm,11mm){$\vdots$};
\node at (5mm,-15mm){iasc$^+_+$-pair};
\end{tikzpicture}
\hspace{1em}
\begin{tikzpicture}[scale = 0.6]
\def\hezi{-- +(5mm,0mm) -- +(5mm,5mm) -- +(0mm,5mm) -- cycle [line width=0.6pt]}
\def\judy{-- +(5mm,0mm) -- +(5mm,5mm) -- +(0mm,5mm) -- cycle [line width=0.6pt,fill=gainsboro]}
\def\nbx{-- +(5mm,0mm) -- +(5mm,5mm) -- +(0mm,5mm) -- cycle [line width=0.6pt,pattern={Lines[angle=45,distance=2.5pt, line width=0.6pt]}, pattern color=black]}
\tikzstyle{bdot}=[circle,fill=black,draw=black,inner sep=1.2]
\tikzstyle{gdot}=[circle,fill=royalblue,draw=royalblue,inner sep=1.5]
\draw (5mm,0mm)\judy;
\draw (0mm,15mm)\nbx;
\draw[line width=0.6pt] (5mm,-5mm)--(5mm,25mm);
\draw[line width=0.6pt] (0mm,-5mm)--(0mm,25mm);
\draw[line width=0.6pt] (10mm,-5mm)--(10mm,25mm);
\draw[line width=0.6pt] (-5mm,20mm)--(15mm,20mm);
\draw[line width=0.6pt] (-5mm,15mm)--(15mm,15mm);
\draw[line width=0.6pt] (-5mm,5mm)--(15mm,5mm);
\draw[line width=0.6pt] (-5mm,0mm)--(15mm,0mm);
\node at (-3mm,11mm){$\vdots$};
\node at (5mm,-15mm){iasc$^-_+$-pair};
\end{tikzpicture}
\hspace{1em}
\begin{tikzpicture}[scale = 0.6]
\def\hezi{-- +(5mm,0mm) -- +(5mm,5mm) -- +(0mm,5mm) -- cycle [line width=0.6pt]}
\def\judy{-- +(5mm,0mm) -- +(5mm,5mm) -- +(0mm,5mm) -- cycle [line width=0.6pt,fill=gainsboro]}
\def\nbx{-- +(5mm,0mm) -- +(5mm,5mm) -- +(0mm,5mm) -- cycle [line width=0.6pt,pattern={Lines[angle=45,distance=2.5pt, line width=0.6pt]}, pattern color=black]}
\tikzstyle{bdot}=[circle,fill=black,draw=black,inner sep=1.2]
\tikzstyle{gdot}=[circle,fill=royalblue,draw=royalblue,inner sep=1.5]
\draw (0mm,0mm)\nbx;
\draw (5mm,15mm)\judy;
\draw[line width=0.6pt] (5mm,-5mm)--(5mm,25mm);
\draw[line width=0.6pt] (0mm,-5mm)--(0mm,25mm);
\draw[line width=0.6pt] (10mm,-5mm)--(10mm,25mm);
\draw[line width=0.6pt] (-5mm,20mm)--(15mm,20mm);
\draw[line width=0.6pt] (-5mm,15mm)--(15mm,15mm);
\draw[line width=0.6pt] (-5mm,5mm)--(15mm,5mm);
\draw[line width=0.6pt] (-5mm,0mm)--(15mm,0mm);
\node at (-3mm,11mm){$\vdots$};
\node at (5mm,-15mm){iasc$^+_-$-pair};
\end{tikzpicture}
\hspace{1em}
\begin{tikzpicture}[scale = 0.6]
\def\hezi{-- +(5mm,0mm) -- +(5mm,5mm) -- +(0mm,5mm) -- cycle [line width=0.6pt]}
\def\nbx{-- +(5mm,0mm) -- +(5mm,5mm) -- +(0mm,5mm) -- cycle [line width=0.6pt,pattern={Lines[angle=45,distance=2.5pt, line width=0.6pt]}, pattern color=black]}
\tikzstyle{bdot}=[circle,fill=black,draw=black,inner sep=1.2]
\tikzstyle{gdot}=[circle,fill=royalblue,draw=royalblue,inner sep=1.5]
\draw (5mm,0mm)\nbx;
\draw (0mm,15mm)\nbx;
\draw[line width=0.6pt] (5mm,-5mm)--(5mm,25mm);
\draw[line width=0.6pt] (0mm,-5mm)--(0mm,25mm);
\draw[line width=0.6pt] (10mm,-5mm)--(10mm,25mm);
\draw[line width=0.6pt] (-5mm,20mm)--(15mm,20mm);
\draw[line width=0.6pt] (-5mm,15mm)--(15mm,15mm);
\draw[line width=0.6pt] (-5mm,5mm)--(15mm,5mm);
\draw[line width=0.6pt] (-5mm,0mm)--(15mm,0mm);
\node at (-3mm,11mm){$\vdots$};
\node at (5mm,-15mm){iasc$^-_-$-pair};
\end{tikzpicture}
\end{center}

By changing only $\des^B$ to $\des_B$ in Definitions \ref{defd+type} and \ref{defd-type},
we get the definitions and related notations of $d_h^\pm$-types, $d_v^\pm$-types and $d^\pm$-types
under the $r$-order.
The following propositions describe the distributions of the $d_h^\pm$-types and $d_v^\pm$-types on the grid lines.

\begin{prop}\label{prop-dhdes-r}
The $d_h$-types of the grid points on the middle horizontal line of
des-pairs and asc-pairs in the grid  are indicated in the figure below,
where the numbers above are $d_h^+$-types,
and the numbers below are $d_h^-$-types.
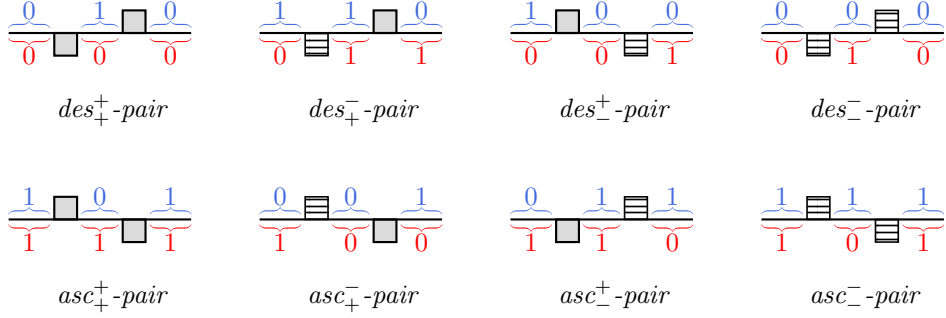
\begin{center}

\end{center}
\end{prop}

\begin{prop}\label{prop-dviasc-r}
The $d_v$-types of the grid points on the middle vertical line of  ides-pairs and iasc-pairs in the grid are indicated in the figure below,
where the numbers on the left are $d_v^+$-types,
and the numbers on the right are $d_v^-$-types.
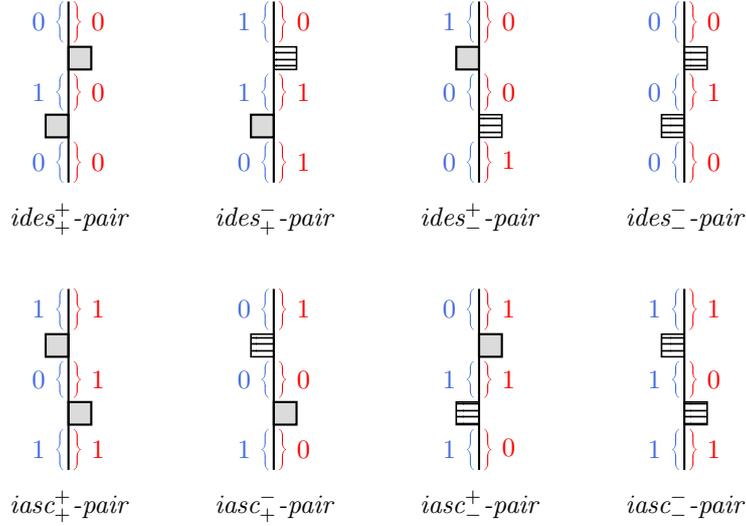
\begin{center}
%
\end{center}
\end{prop}

\begin{prop}\label{prop-dtype-rim-r}
The $d_h$-types of the grid points on the top and bottom boundaries of the grid $P_\pi$
are indicated in the figure below,
where the numbers above are $d_h^+$-types,
and the numbers below are $d_h^-$-types.
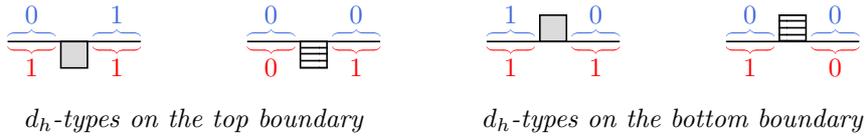
\begin{center}
%
\end{center}
And the $d_v$-types of the grid points on the left and right boundaries of the grid are indicated in the figure below,
where the numbers on the left are $d_v^+$-types,
and the numbers on the right are $d_v^-$-types.
\begin{center}
%
\end{center}
\end{prop}

Combining Propositions \ref{prop-dhdes-r}, \ref{prop-dviasc-r} and \ref{prop-dtype-rim-r},
we see the $d$-types of the four corners of the filled squares
in the permutation grid have fixed values.

\begin{prop}\label{prop-4corners-r}
Let $\pi\in\B_n$ and $P_\pi$ be its grid, then for $1\leq i,j\leq n$, we have
\[
d^+(i,j)=d^+(i+1,j+1)=(0,0)\quad\mbox{and}\quad d^+(i,j+1)=d^+(i+1,j)=(1,1),
\]
if $\ij$ is a positive square in $P_\pi$, and
\[
d^-(i,j)=d^-(i+1,j+1)=(0,0)\quad\mbox{and}\quad d^-(i,j+1)=d^-(i+1,j)=(1,1).
\]
and if $\ij$ is a negative square in $P_\pi$.

\begin{center}
\begin{tikzpicture}[scale=1.5]
\def\judy{-- +(5mm,0mm) -- +(5mm,5mm) -- +(0mm,5mm) -- cycle [line width=1pt,fill=gainsboro]}
\tikzstyle{gdot}=[circle,fill=royalblue,draw=royalblue,inner sep=2]

\draw(15mm,0mm)\judy;
\draw[line width=1pt,dashed](10mm,5mm)--(25mm,5mm);
\draw[line width=1pt,dashed](10mm,0mm)--(25mm,0mm);
\draw[line width=1pt,dashed] (15mm,-5mm)--(15mm,10mm);
\draw[line width=1pt,dashed] (20mm,-5mm)--(20mm,10mm);

\node[gdot] at (15mm,5mm){};
\node at (12mm,7mm){\royalblue{$(0,0)$}};

\node[gdot] at (20mm,0mm){};
\node at (23mm,-2mm){\royalblue{$(0,0)$}};

\node[gdot] at (15mm,0mm){};
\node at (12mm,-2mm){\royalblue{$(1,1)$}};

\node[gdot] at (20mm,5mm){};
\node at (23mm,7mm){\royalblue{$(1,1)$}};

\node at (17.5mm,-10mm){$d^+$-types on a positive square};

\end{tikzpicture}
\hspace{3em}
\begin{tikzpicture}[scale =1.5]
\def\nbx{-- +(5mm,0mm) -- +(5mm,5mm) -- +(0mm,5mm) -- cycle [line width=1pt,pattern={Lines[angle=45,distance=4pt, line width=1pt]}, pattern color=black]}
\tikzstyle{rdot}=[circle,fill=red,draw=red,inner sep=2]

\draw(15mm,0mm)\nbx;
\draw[line width=1pt,dashed](10mm,5mm)--(25mm,5mm);
\draw[line width=1pt,dashed](10mm,0mm)--(25mm,0mm);
\draw[line width=1pt,dashed] (15mm,-5mm)--(15mm,10mm);
\draw[line width=1pt,dashed] (20mm,-5mm)--(20mm,10mm);

\node[rdot] at (15mm,5mm){};
\node at (12mm,7mm){\red{$(0,0)$}};

\node[rdot] at (20mm,0mm){};
\node at (23mm,-2mm){\red{$(0,0)$}};

\node[rdot] at (15mm,0mm){};
\node at (12mm,-2mm){\red{$(1,1)$}};

\node[rdot] at (20mm,5mm){};
\node at (23mm,7mm){\red{$(1,1)$}};

\node at (17.5mm,-10mm){$d^-$-types on a negative square};

\end{tikzpicture}
\end{center}
\end{prop}

Consistent with Definition \ref{def-path}, for $p,q\in\{0,1\}$,
we can construct $p_h^\pm$-paths and $q_v^\pm$-paths
based on $d_h^\pm$-types and $d_v^\pm$-types under the $r$-order.

\begin{example}
  The figures below present each one of $p_x^\ast$-paths in the permutation grid of $\pi=31\bar{6}5\bar{2}4$, where $p\in\{0,1\},x\in\{h,v\}$ and $\ast\in\{+,-\}$.
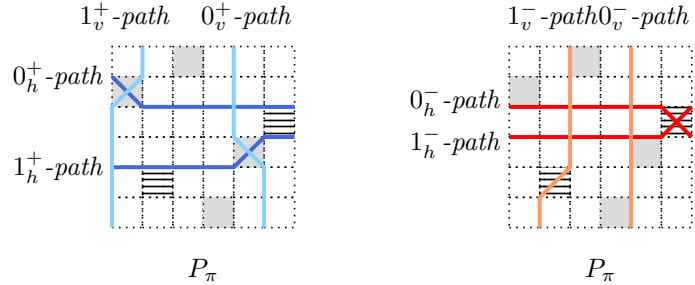
\begin{center}
\begin{tikzpicture}[scale =0.8]
\def\hezi{-- +(5mm,0mm) -- +(5mm,5mm) -- +(0mm,5mm) -- cycle [line width=0.6pt,dotted]}
\def\judy{-- +(5mm,0mm) -- +(5mm,5mm) -- +(0mm,5mm) -- cycle [line width=0.6pt,fill=gainsboro,dotted]}
\def\nbx{-- +(5mm,0mm) -- +(5mm,5mm) -- +(0mm,5mm) -- cycle [line width=0.6pt,pattern={Lines[angle=45,distance=2.5pt, line width=0.6pt]}, pattern color=black,dotted]}
\tikzstyle{cc}=[circle,draw=black,fill=yellow, line width=0.5pt, inner sep=1.5]
\draw (0mm,0mm)\hezi;
\draw (5mm,0mm)\hezi;
\draw (10mm,0mm)\judy;
\draw (15mm,0mm)\hezi;
\draw (20mm,0mm)\hezi;
\draw (25mm,0mm)\hezi;
\draw (0mm,-5mm)\judy;
\draw (5mm,-5mm)\hezi;
\draw (10mm,-5mm)\hezi;
\draw (15mm,-5mm)\hezi;
\draw (20mm,-5mm)\hezi;
\draw (25mm,-5mm)\hezi;
\draw (0mm,-10mm)\hezi;
\draw (5mm,-10mm)\hezi;
\draw (10mm,-10mm)\hezi;
\draw (15mm,-10mm)\hezi;
\draw (20mm,-10mm)\hezi;
\draw (25mm,-10mm)\nbx;
\draw (0mm,-15mm)\hezi;
\draw (5mm,-15mm)\hezi;
\draw (10mm,-15mm)\hezi;
\draw (15mm,-15mm)\hezi;
\draw (20mm,-15mm)\judy;
\draw (25mm,-15mm)\hezi;
\draw (0mm,-20mm)\hezi;
\draw (5mm,-20mm)\nbx;
\draw (10mm,-20mm)\hezi;
\draw (15mm,-20mm)\hezi;
\draw (20mm,-20mm)\hezi;
\draw (25mm,-20mm)\hezi;
\draw (0mm,-25mm)\hezi;
\draw (5mm,-25mm)\hezi;
\draw (10mm,-25mm)\hezi;
\draw (15mm,-25mm)\judy;
\draw (20mm,-25mm)\hezi;
\draw (25mm,-25mm)\hezi;
\node at (-9mm,0mm){$0_h^+$-path};
\draw(0mm,0mm)--(5mm,-5mm)--(25mm,-5mm)--(30mm,-5mm)[line width=1.5pt,draw=royalblue];

\node at (-9mm,-15mm){$1_h^+$-path};
\draw (0mm,-15mm)--(20mm,-15mm)--(25mm,-10mm)--(30mm,-10mm)[line width=1.5pt,draw=royalblue];
\node at (2mm,10mm){$1_v^+$-path};
\draw (5mm,5mm)--(5mm,0mm)--(0mm,-5mm)--(0mm,-25mm)[line width=1.5pt,draw=babyblue];
\node at (22mm,10mm){$0_v^+$-path};
\draw(20mm,5mm)--(20mm,-10mm)--(25mm,-15mm)--(25mm,-25mm)[line width=1.5pt,draw=babyblue];

\node at (15mm,-32mm){$P_\pi$};
\end{tikzpicture}
\hspace{3em}
\begin{tikzpicture}[scale =0.8]
\def\hezi{-- +(5mm,0mm) -- +(5mm,5mm) -- +(0mm,5mm) -- cycle [line width=0.6pt,dotted]}
\def\judy{-- +(5mm,0mm) -- +(5mm,5mm) -- +(0mm,5mm) -- cycle [line width=0.6pt,fill=gainsboro,dotted]}
\def\nbx{-- +(5mm,0mm) -- +(5mm,5mm) -- +(0mm,5mm) -- cycle [line width=0.6pt,pattern={Lines[angle=45,distance=2.5pt, line width=0.6pt]}, pattern color=black,dotted]}
\tikzstyle{cc}=[circle,draw=black,fill=yellow, line width=0.5pt, inner sep=1.5]
\draw (0mm,0mm)\hezi;
\draw (5mm,0mm)\hezi;
\draw (10mm,0mm)\judy;
\draw (15mm,0mm)\hezi;
\draw (20mm,0mm)\hezi;
\draw (25mm,0mm)\hezi;
\draw (0mm,-5mm)\judy;
\draw (5mm,-5mm)\hezi;
\draw (10mm,-5mm)\hezi;
\draw (15mm,-5mm)\hezi;
\draw (20mm,-5mm)\hezi;
\draw (25mm,-5mm)\hezi;
\draw (0mm,-10mm)\hezi;
\draw (5mm,-10mm)\hezi;
\draw (10mm,-10mm)\hezi;
\draw (15mm,-10mm)\hezi;
\draw (20mm,-10mm)\hezi;
\draw (25mm,-10mm)\nbx;
\draw (0mm,-15mm)\hezi;
\draw (5mm,-15mm)\hezi;
\draw (10mm,-15mm)\hezi;
\draw (15mm,-15mm)\hezi;
\draw (20mm,-15mm)\judy;
\draw (25mm,-15mm)\hezi;
\draw (0mm,-20mm)\hezi;
\draw (5mm,-20mm)\nbx;
\draw (10mm,-20mm)\hezi;
\draw (15mm,-20mm)\hezi;
\draw (20mm,-20mm)\hezi;
\draw (25mm,-20mm)\hezi;
\draw (0mm,-25mm)\hezi;
\draw (5mm,-25mm)\hezi;
\draw (10mm,-25mm)\hezi;
\draw (15mm,-25mm)\judy;
\draw (20mm,-25mm)\hezi;
\draw (25mm,-25mm)\hezi;
\node at (-9mm,-4mm){$0_h^-$-path};
\draw (0mm,-5mm)--(25mm,-5mm)--(30mm,-10mm)[line width=1.5pt,draw=red];
\node at (-9mm,-11mm){$1_h^-$-path};
\draw(0mm,-10mm)--(25mm,-10mm)--(30mm,-5mm)[line width=1.5pt,draw=red];

\node at (22mm,10mm){$0_v^-$-path};
\draw (20mm,5mm)--(20mm,-25mm)[line width=1.5pt,draw=orange];
\node at (7mm,10mm){$1_v^-$-path};
\draw (10mm,5mm)--(10mm,-15mm)--(5mm,-20mm)--(5mm,-25mm)[line width=1.5pt,draw=orange];

\node at (15mm,-32mm){$P_\pi$};
\end{tikzpicture}
\end{center}
\end{example}

These paths defined by $d$-types under the $r$-order
also exhibit similar behavior as those defined under the natural order
when they encounter filled squares as described in Theorem \ref{lempath}.

\begin{thm}\label{lempathr}
Let $\pi\in\mathfrak{B}_n$, then in the signed permutation grid $P_\pi$,
\begin{enumerate}[(a)]
\item each $0_h^+$-path (resp. $1_h^+$-path) goes from the left boundary to the right boundary along the horizontal grid lines except for carrying out a southeast (resp. northeast) step when encountering a positive square;
\item each $0_v^+$-path (resp. $1_v^+$-path) goes from the top boundary to the bottom boundary along the vertical grid lines except for carrying out a southeast (resp. southwest) step when encountering a positive square;
\item each $0_h^-$-path (resp. $1_h^-$-path) goes from the left boundary to the right boundary along the horizontal grid lines except for carrying out a southeast (resp. northeast) step when encountering a negative square;
\item each $0_v^-$-path (resp. $1_v^-$-path) goes from the top boundary to the bottom boundary along the vertical grid lines except for carrying out a southeast (resp. southwest) step when encountering a negative square.
\end{enumerate}
\end{thm}

Following the proof of Theorem \ref{thmpaths} but dealing with more complicated cases,
we give the counting formulas for the paths in terms of $\des_B(\pi)$ and $\ides_B(\pi)$ with any given $\pi\in\B_n$.

\begin{thm}\label{thmpathsr}
Let $\pi\in\mathfrak{B}_n$, then in the signed permutation grid $P_\pi$,
\begin{enumerate}[(a)]
    \item the number of $0_h^+$-paths is  $\des_B(\pi)+1$,
    \item the number of $0_h^-$-paths is  $\des_B(\pi)$,
    \item the number of $1_h^+$-paths is  $n-\des_B(\pi)$,
    \item the number of $1_h^-$-paths is  $n-\des_B(\pi)+1$;
\end{enumerate}
and
\begin{enumerate}[(a)]
\addtocounter{enumi}{4}
    \item the number of $0_v^+$-paths is  $\ides_B(\pi)+1$,
    \item the number of $0_v^-$-paths is  $\ides_B(\pi)$,
    \item the number of $1_v^+$-paths is  $n-\ides_B(\pi)$,
    \item the number of $1_v^-$-paths is  $n-\ides_B(\pi)+1$.
\end{enumerate}
\end{thm}

By Theorems \ref{lempathr} and \ref{thmpathsr}, in any grid $P_\pi$,
we can use the statistics $\des_B(\pi)$, $\ides_B(\pi)$ and $n(\pi)$
to determine the number of grid points of given $d^+$-type or $d^-$-type.

\begin{thm}\label{thm-enum-dtype-r}
Let $\pi\in\B_n$, then in the grid $P_\pi$, there are
\begin{enumerate}[(a)]
    \item $(\des_B(\pi)+1)(\ides_B(\pi)+1)-n(\pi)+n$ grid points of $d^+$-type$ (0,0)$,
    \item $(\ides_B(\pi)+1)(n-\des_B(\pi))+n(\pi)-n$ grid points of $d^+$-type $(1,0)$,
    \item $(\des_B(\pi)+1)(n-\ides_B(\pi))+n(\pi)-n$ grid points of $d^+$-type $(0,1)$,
    \item $(n-\des_B(\pi))(n-\ides_B(\pi))-n(\pi)+n$ grid points of $d^+$-type $(1,1)$;
\end{enumerate}
and
\begin{enumerate}[(a)]
\addtocounter{enumi}{4}
    \item $\des_B(\pi)\ides_B(\pi)+n(\pi)$ grid points of $d^-$-type$ (0,0)$,
    \item $\ides_B(\pi)(n-\des_B(\pi)+1)-n(\pi)$ grid points of $d^-$-type $(1,0)$,
    \item $\des_B(\pi)(n-\ides_B(\pi)+1)-n(\pi)$ grid points of $d^-$-type $(0,1)$,
    \item $(n-\des_B(\pi)+1)(n-\ides_B(\pi)+1)+n(\pi)$ grid points of $d^-$-type $(1,1)$.
\end{enumerate}
\end{thm}

\begin{exam}\label{examdtype-r}
Let $\pi=4\overline{3}1\overline{2}5\in \B_{5}$ with $des_B(\pi)=2$, $ides_B(\pi)=3$ and $n(\pi)=2$,
then the following diagrams show the $d^+$-type and $d^-$-type of each grid point,
where $\#d^\ast(p,q)$ denote the number of grid points with $d^\ast$-type $(p,q)$
for $\ast\in\{+,-\}$ and $p,q\in\{0,1\}$.
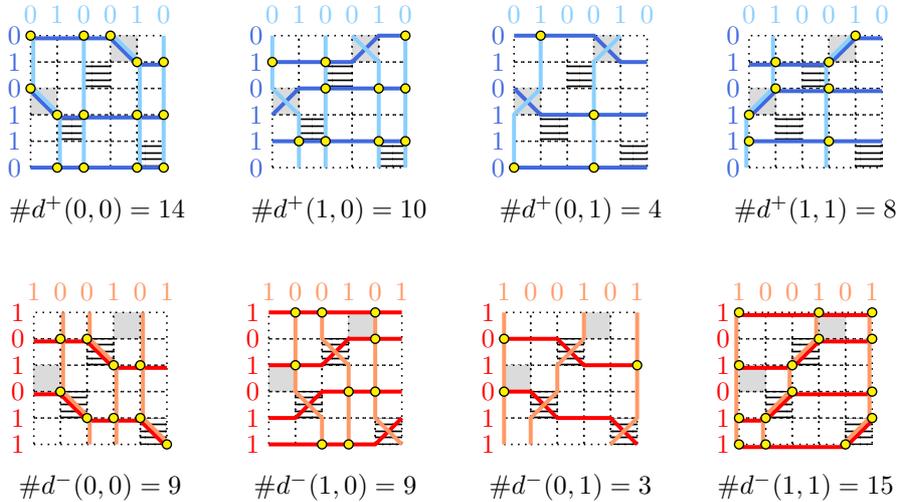
\begin{center}

\end{center}
\end{exam}

\section{Recurrences of $\underline{b}_{n,i,j}$, $I_{n,k}^\ub$ and $J_{n,k}^\ub$ under the $r$-order}\label{Sec-recunderall}

In this section, we present combinatorial proofs of Theorems \ref{thmof_bnij}, \ref{thmIubnk} and \ref{thmJubnk}.
Note that some  symbols in the proofs of Theorems \ref{thmof^obnij}, \ref{thmIobnk} and \ref{thmJobnk} are reused ,
but keep in mind that the objects involved in these symbols are defined under the $r$-order.

{\noindent{\emph{Combinatorial Proof of Theorem \ref{thmof_bnij}.\hskip 2pt}}
We have
$$
\underline{\mathcal{B}}_{n,i,j}=\{\sigma\in\B_n\mid \mbox{ $\des_B(\sigma)=i$ and $\ides_B(\sigma)=j$}\}
$$
and define
$$
\mathcal{C}_{n,i,j}=\{(\sigma, k) \mid \pi\in\underline{\mathcal{B}}_{n, i, j} \text{ and }1\leq k\leq n \}
$$
whose cardinality $|\mathcal{C}_{n,i,j}|=n\underline{b}_{n,i,j}$
is the left side of \eqref{rec_bnij}.
For $\ast\in\{+,-\}$ and $p,q\in\{0,1\}$, let
\begin{equation*}
\mathcal{D}_{n,i,j}^{(p,q)^\ast}=\left\{(\pi,(r,s))\mid \text{ $\pi\in\underline{\mathcal{B}}_{n,i,j}$ and $d^\ast(r,s)=(p,q)$ in $P_\pi$}\right\}.
\end{equation*}
Then it follows from Theorem \ref{thm-enum-dtype-r} that
the cardinality of the set
\begin{equation*}
  \mathcal{D}_{n-1,i,j}=\biguplus_{p,q\in\{0,1\}}\mathcal{D}_{n-1,i-p,j-q}^{(p,q)^+}\uplus\mathcal{D}_{n-1,i-p,j-q}^{(p,q)^-}
\end{equation*}
equals  the right side of \eqref{rec_bnij}.
Define
\begin{equation*}
\Psi((\pi,(r,s)))=(\sigma, r),
\end{equation*}
where for $p,q\in\{0,1\}$,
\begin{equation*}
\sigma=\left\{
      \begin{array}{ll}
        \varphi_{(r,s)}(\pi), & \hbox{if $(\pi,(r,s))\in\mathcal{D}_{n-1,i-p,j-q}^{(p,q)^+}$ ,} \\[6pt]
        \overline{\varphi}_{(r,s)}(\pi), & \hbox{if $(\pi,(r,s))\in\mathcal{D}_{n-1,i-p,j-q}^{(p,q)^-}$ .}
      \end{array}
    \right.
\end{equation*}
Clearly $\Psi$ is a bijection between ${\mathcal{D}}_{n-1,i,j}$ and $\mathcal{C}_{n,i,j}$,
which completes the proof. \qed

{\noindent{\emph{Combinatorial Proof of Theorem \ref{thmIubnk}.\hskip 2pt}}
Recall that we denote $|\mathcal{I}_{n,k}^\ub|$ by $I_{n,k}^\ub$,
where $\mathcal{I}_{n,k}^\ub$ is the set of
signed involutions in $\mathcal{I}^B_n$ with $k$ $\ub$-descents.
Let
\[
\mathcal{G}_{n,k}=\{(\sigma,i)\mid \pi\in \mathcal{I}_{n,k}^\ub\text{ and }1\leq i\leq n \},
\]
then we see $|\mathcal{G}_{n,k}|=nI_{n,k}^B$
that is equal to the left side of \eqref{recIubnk}.
Let
\begin{equation*}
  \begin{aligned}
\mathcal{E}^{(1)^+}_{n,k}=\{(\pi,(i,j))\mid&\text{ $\pi\in\mathcal{I}^\ub_{n,k}$ and $(i,j)$ is the grid point where}\\
    &\text{ a $0_h^+$-path first touches the main diagonal of $P_\pi$}\}
  \end{aligned}
\end{equation*}
and
\begin{equation*}
  \begin{aligned}
\mathcal{E}^{(1)^-}_{n,k}=\{(\pi,(i,j))\mid&\text{ $\pi\in\mathcal{I}^\ub_{n,k}$ and $(i,j)$ is the grid point where }\\
    &\text{ a $0_h^-$-path first touches the main diagonal of $P_\pi$}\}.
  \end{aligned}
\end{equation*}
By Theorem \ref{thmpathsr},
we have
$$\left|\mathcal{E}^{(1)^+}_{n-1,k}\uplus\mathcal{E}^{(1)^-}_{n-1,k}\right|
=(2k+1)I_{n-1,k}^\ub,$$
which is the \textbf{first} term in the right side of \eqref{recIubnk}.
For  $(\pi,(i,j))\in\mathcal{E}^{(1)^+}_{n-1,k}\uplus \mathcal{E}^{(1)^-}_{n-1,k}$,
define
\begin{equation*}
  \Psi((\pi,(i,j)))=(\sigma,i)
\end{equation*}
with
\begin{equation}\label{eqpfIn1st-r}
\sigma=\left\{
      \begin{array}{ll}
        \varphi_{(i,j)}(\pi), & \hbox{if $(\pi,(i,j))\in\mathcal{E}^{(1)^+}_{n-1,k}$,} \\[6pt]
        \overline{\varphi}_{(i,j)}(\pi), & \hbox{if $(\pi,(i,j))\in\mathcal{E}^{(1)^-}_{n-1,k}$.}
      \end{array}
    \right.
\end{equation}

Let
\begin{equation*}
  \begin{aligned}
\mathcal{E}^{(2)^+}_{n,k}=\{(\pi,(i,j))\mid& \text{ $\pi\in\mathcal{I}_{n,k}^\ub$ and $(i,j)$ is the grid point where a $1_h^+$-path}\\
&\text{ touches the main diagonal of $P_\pi$ \emph{or} a positive square on it}\},
  \end{aligned}
\end{equation*}
and
\begin{equation*}
  \begin{aligned}
\mathcal{E}^{(2)^-}_{n,k}=\{(\pi,(i,j))\mid& \text{$\pi\in\mathcal{I}_{n,k}^\ub$ and $(i,j)$ is the grid point where  a $1_h^-$-path}\\
&\text{ touches the main diagonal of $P_\pi$  \emph{or} a negative square on it}\}.
  \end{aligned}
\end{equation*}
Then we obtain $$\left|\mathcal{E}^{(2)^+}_{n-1,k-1}\uplus\mathcal{E}^{(2)^-}_{n-1,k-1}\right|
=(2n-2k+1)I_{n-1,k-1}^\ub$$ by Theorem \ref{thmpathsr},
and $(2n-2k+1)I_{n-1,k-1}^\ub$ is the \textbf{second} term in the right side of \eqref{recIubnk}.
For  $(\pi,(i,j))\in \mathcal{E}^{(2)^+}_{n-1,k-1}\uplus\mathcal{E}^{(2)^-}_{n-1,k-1}$,
define
\begin{equation*}
\Psi((\pi, (i,j)))=(\sigma,i)
\end{equation*}
with
\begin{equation}\label{eqpfIn2nd-r}
\sigma=\left\{
      \begin{array}{ll}
        \varphi_{(i,j)}(\pi), & \hbox{if $(\pi,(i,j))\in\mathcal{E}^{(2)^+}_{n-1,k-1}$ ,} \\[6pt]
        \overline{\varphi}_{(i,j)}(\pi), & \hbox{if $(\pi,(i,j))\in\mathcal{E}^{(2)^-}_{n-1,k-1}$.}
      \end{array}
    \right.
\end{equation}

Let
\begin{equation*}
\mathcal{E}^{(3)^+}_{n,k}=\{(\pi,(i,j))\mid\text{ $\pi\in\mathcal{I}^\ub_{n,k}$ and $d^+(i,j)=(0,0)$ in $P_\pi$}\}.
\end{equation*}
and
\begin{equation*}
\mathcal{E}^{(3)^-}_{n,k}=\{(\pi,(i,j))\mid\text{ $\pi\in\mathcal{I}^\ub_{n,k}$ and $d^-(i,j)=(0,0)$ in $P_\pi$}\}.
\end{equation*}
By Theorem \ref{thm-enum-dtype-r},
we see
$$\left|\mathcal{E}^{(3)^+}_{n-2,k}\uplus\mathcal{E}^{(3)^-}_{n-2,k}\right|=
\left((n-1)+2k(k+1)\right) I_{n-2,k}^\ub,$$
which is the \textbf{third} term in the right side of \eqref{recIubnk}.
For $(\sigma,(i,j))\in \mathcal{E}^{(3)^+}_{n-2,k}\uplus\mathcal{E}^{(3)^-}_{n-2,k}$, set
\begin{equation*}
\Psi((\pi, (i,j)))=(\sigma,\chi_{ij})
\end{equation*}
with
\begin{equation}\label{eqpfIn3rd-r}
\sigma=\left\{
\begin{array}{ll}
\xi_{(i,j)}(\pi), & \hbox{if $(\pi,(i,j))\in\mathcal{E}^{(3)^+}_{n-2,k}$ and $i\neq j$,} \\[6pt]
\eta_{i}(\pi), & \hbox{if $(\pi,(i,j))\in\mathcal{E}^{(3)^+}_{n-2,k}$ and $i=j$,} \\[6pt]
\overline{\xi}_{(i,j)}(\pi), & \hbox{if $(\pi,(i,j))\in\mathcal{E}^{(3)^-}_{n-2,k}$ and $i\neq j$,} \\[6pt]
\overline{\eta}_{i}(\pi), & \hbox{if $(\pi,(i,j))\in\mathcal{E}^{(3)^-}_{n-2,k}$ and $i=j$.}
\end{array}
    \right.
\end{equation}

Let
\begin{equation*}
  \begin{aligned}
  \mathcal{E}^{(4)^+}_{n,k}=\{&(\pi,(i,j))\mid\text{ $\pi\in\mathcal{I}^\ub_{n,k}$ and $d^+(i,j)=(0,1)$ or $(1,0)$ in $P_\pi$}\\
  &\quad\cup\{(\pi,(i,i))\mid\text{ $\pi\in\mathcal{I}_{n,k}^\ub$ and $d^+(i,i)=(0,0)$ or $(1,1)$ in $P_\pi$}\},
  \end{aligned}
\end{equation*}
and
\begin{equation*}
  \begin{aligned}
  \mathcal{E}^{(4)^-}_{n,k}=\{&(\pi,(i,j))\mid\text{ $\pi\in\mathcal{I}^\ub_{n,k}$ and $d^-(i,j)=(0,1)$ or $(1,0)$ in $P_\pi$}\\
  &\quad\cup\{(\pi,(i,i))\mid\text{ $\pi\in\mathcal{I}_{n,k}^\ub$ and $d^-(i,i)=(0,0)$ or $(1,1)$ in $P_\pi$}\}.
  \end{aligned}
\end{equation*}
It follows from Theorem \ref{thm-enum-dtype-r} that $$\left|\mathcal{E}^{(4)^+}_{n-2,k-1}\uplus\mathcal{E}^{(4)^-}_{n-2,k-1}\right|
=(2(n-1)+4(k-1)(n-k-1)) I_{n-2,k-1}^\ub,$$
which is the \textbf{fourth} term in the right side of \eqref{recIubnk}.
For $(\sigma,(i,j))\in \mathcal{E}^{(4)^+}_{n-2,k-1}\uplus\mathcal{E}^{(4)^-}_{n-2,k-1}$,
let
\begin{equation*}
\Psi((\pi, (i,j)))=(\sigma,\chi_{ij}),
\end{equation*}
with
\begin{equation}\label{eqpfIn4th-r}
(\sigma,k)=\left\{
      \begin{array}{ll}
       \xi_{(i,j)}(\pi), & \hbox{if $(\pi,(i,j))\in\mathcal{E}^{(4)^+}_{n-2,k-1}$ and $i\neq j$,} \\[6pt]
      \eta'_{i}(\pi), & \hbox{if $(\pi,(i,i))\in\mathcal{E}^{(4)^+}_{n-2,k-1}$ and $d^+(i,i)=(0,0)$,} \\[6pt]
      \eta_{i}(\pi), & \hbox{if $(\pi,(i,i))\in\mathcal{E}^{(4)^+}_{n-2,k-1}$ and $d^+(i,i)=(1,1)$,} \\[6pt]
      \overline{\xi}_{(i,j)}(\pi), & \hbox{if $(\pi,(i,j))\in\mathcal{E}^{(4)^-}_{n-2,k-1}$ and $i\neq j$,} \\[6pt]
      \overline{\eta'}_{i}(\pi), & \hbox{if $(\pi,(i,i))\in\mathcal{E}^{(4)^-}_{n-2,k-1}$ and $d^-(i,i)=(0,0)$,} \\[6pt]
      \overline{\eta}_{i}(\pi), & \hbox{if $(\pi,(i,i))\in\mathcal{E}^{(4)^-}_{n-2,k-1}$ and $d^-(i,i)=(1,1)$.}
      \end{array}
    \right.
\end{equation}

Let
\begin{equation*}
\mathcal{E}^{(5)^+}_{n,k}=\{(\pi,(i,j))\mid\text{ $\pi\in\mathcal{I}^\ub_{n,k}$ and $d^+(i,j)=(1,1)$ in $P_\pi$}\}.
\end{equation*}
and
\begin{equation*}
\mathcal{E}^{(5)^-}_{n,k}=\{(\pi,(i,j))\mid\text{ $\pi\in\mathcal{I}^\ub_{n,k}$ and $d^-(i,j)=(1,1)$ in $P_\pi$}\}.
\end{equation*}
By Theorem \ref{thm-enum-dtype-r}, we have
$$\left|\mathcal{E}^{(5)^+}_{n-2,k-2}\uplus
\mathcal{E}^{(5)^-}_{n-2,k-2}\right|
=\left((2n-3)(n-1)+2(k-2)(k-2n+1)\right)
I_{n-2,k-2}^\ub$$
that is the \textbf{fifth} term in the right side of \eqref{recIubnk}.
For $(\sigma,(i,j))\in \mathcal{E}^{(5)^+}_{n-2,k-2}\uplus\mathcal{E}^{(5)^-}_{n-2,k-2}$,
define
\begin{equation*}
\Psi((\pi, (i,j)))=(\sigma,\chi_{ij}),
\end{equation*}
with
\begin{equation}\label{eqpfIn5th-r}
\sigma=\left\{
      \begin{array}{ll}
       \xi_{(i,j)}(\pi), & \hbox{if $(\pi,(i,j))\in\mathcal{E}^{(5)^+}_{n-2,k-2}$ and $i\neq j$,} \\[6pt]
      \eta'_{i}(\pi), & \hbox{if $(\pi,(i,j))\in\mathcal{E}^{(5)^+}_{n-2,k-2}$ and $i=j$,} \\[6pt]
      \overline{\xi}_{(i,j)}(\pi), & \hbox{if $(\pi,(i,j))\in\mathcal{E}^{(5)^-}_{n-2,k-2}$ and $i\neq j$,} \\[6pt]
      \overline{\eta'}_{i}(\pi), & \hbox{if $(\pi,(i,j))\in\mathcal{E}^{(5)^-}_{n-2,k-2}$ and $i=j$.}
      \end{array}
    \right.
\end{equation}

Thus, the cardinality of the set
\begin{equation*}
{\mathcal{E}}_{n,k}=
\mathcal{E}^{(1)^\pm}_{n-1,k}\uplus\mathcal{E}^{(2)^\pm}_{n-1,k-1}\uplus\mathcal{E}^{(3)^\pm}_{n-2,k}
\uplus\mathcal{E}^{(4)^\pm}_{n-2,k-1}\uplus\mathcal{E}^{(5)^\pm}_{n-2,k-2}
\end{equation*}
is counted by the right side of \eqref{recIobnk}.
And referring the proof of Theorem \ref{thmIobnk},
we see that the combination of \eqref{eqpfIn1st-r}--\eqref{eqpfIn5th-r}
establishes the bijection
\[
\Psi\colon{\mathcal{E}}_{n,k}\leftrightarrow{\mathcal{G}}_{n,k}.
\]
\qed

{\noindent{\emph{Combinatorial Proof of Theorem \ref{thmJubnk}.}\hskip 2pt}}
We proof the following equivalent form of \eqref{recJubnk}:
\begin{equation}\label{recJubnkpf}
\begin{aligned}
2n J^\ub_{2 n, k}=& {2(k^2+n-1) J^\ub_{2 n-2, k}+2(2(k-1)(2n-k)+1) J^\ub_{2 n-2, k-1} } \\[3pt]
&+2((k-2)(k-4n)+4n^2-3n) J^\ub_{2 n-2, k-2}.
\end{aligned}
\end{equation}
For $1\leq k\leq n$, since $\mathcal{J}_{n,k}^\ub$ is the set of fixed-point free signed involutions in $\mathcal{J}_n^B$ with $k$ $\ub$-descents and $|\mathcal{J}_{n,k}^\ub|=J_{n,k}^\ub$,
by constructing
\begin{equation*}
  \mathcal{H}_{n,k}=
  \{(\pi,i)\mid \pi\in \mathcal{J}^\ub_{n,k}\text{ and }1\leq i\leq n \},
\end{equation*}
we have $\left|\mathcal{H}_{2n,k}\right|=2n J_{2n,k}^\ub$
that equals the left side of \eqref{recJubnkpf}.

Let
\begin{equation*}
 \mathcal{F}^{(1)^+}_{n,k}=\{(\pi,(i,j))\mid \pi\in\mathcal{J}^\ub_{n,k}\text{ and $d^+(i,j)=(0,0)$ in $P_\pi$ with $i\neq j$}
  \}
\end{equation*}
and
\begin{equation*}
  \mathcal{F}^{(1)^-}_{n,k}=\{(\pi,(i,j))\mid \pi\in\mathcal{J}^\ub_{n,k}\text{ and  $d^-(i,j)=(0,0)$ in $P_\pi$ with $i\neq j$}\},
\end{equation*}
then we derive that
\begin{equation*}
 \left|\mathcal{F}^{(1)^+}_{2n-2,k}\uplus\mathcal{F}^{(1)^-}_{2n-2,k}\right|
 =(2k^2+2n-2) J_{2n-2,k}^\ub
\end{equation*}
by Theorems \ref{thmpathsr} and \ref{thm-enum-dtype-r},
which is the \textbf{first} term in the right side of \eqref{recJubnkpf}.
For  $(\pi, (i,j))\in \mathcal{F}^{(1)^+}_{2n-2,k}\uplus\mathcal{F}^{(1)^-}_{2n-2,k}$,
define
\begin{equation*}
\Psi((\pi, (i,j)))=(\sigma,r)
\end{equation*}
by
\begin{equation}\label{eqpfJn1str}
(\sigma,r)=\left\{
      \begin{array}{ll}
      (\xi_{(i,j)}(\pi), \chi_{ij}), & \hbox{if $(\pi,(i,j))\in\mathcal{F}^{(1)^+}_{2n-2,k}$ ,} \\[6pt]
      (\overline{\xi}_{(i,j)}(\pi), \chi_{ij}), & \hbox{if $(\pi,(i,j))\in\mathcal{F}^{(1)^-}_{2n-2,k}$.}
      \end{array}
    \right.
\end{equation}

Let
\begin{equation*}
  \mathcal{F}^{(5)^+}_{n,k}=\{
  (\pi,(i,j))\mid \pi\in\mathcal{J}^\ub_{n,k}\text{ and $d^+(i,j)=(1,0)$ or $(0,1)$ in $P_\pi$ }\},
\end{equation*}
\begin{equation*}
  \mathcal{F}^{(5)^-}_{n,k}=\{
  (\pi,(i,j))\mid \pi\in\mathcal{J}^\ub_{n,k}\text{ and $d^-(i,j)=(1,0)$ or $(0,1)$ in $P_\pi$ }\},
\end{equation*}
and
\begin{equation*}
  \mathcal{F}^{(2)^+}_{n,k}=
  \{
  (\pi,(i,i))_1,(\pi,(i,i))_2\mid \pi\in\mathcal{J}^\ub_{n,k} \text{ and $d^+(i,i)=(0,0)$ in $P_\pi$ }
  \},
\end{equation*}
\begin{equation*}
  \mathcal{F}^{(2)^-}_{n,k}=
  \{
  (\pi,(i,i))_1,(\pi,(i,i))_2\mid \pi\in\mathcal{J}^\ub_{n,k} \text{ and $d^-(i,i)=(0,0)$ in $P_\pi$ }
  \}.
\end{equation*}
Then by Theorems \ref{thmpathsr} and \ref{thm-enum-dtype-r}, we have
\begin{equation*}
\left|\mathcal{F}^{(5)^+}_{2n-2,k-1}\uplus\mathcal{F}^{(5)^-}_{2n-2,k-1}\right|
=4(k-1)(2n-k-1) J_{2n-2,k-1}^\ub
\end{equation*}
and
\begin{equation*}
\left|\mathcal{F}^{(2)^+}_{2n-2,k-1}\uplus\mathcal{F}^{(2)^-}_{2n-2,k-1}\right|
=2(k+(k-1))J_{2n-2,k-1}^\ub,
\end{equation*}
which leads to
\begin{equation*}
\left|\mathcal{F}^{(5)^+}_{2n-2,k-1}\uplus\mathcal{F}^{(5)^-}_{2n-2,k-1}
\uplus\mathcal{F}^{(2)^+}_{2n-2,k-1}\uplus\mathcal{F}^{(2)^-}_{2n-2,k-1}\right|
=2(2(k-1)(2 n-k)+1) J_{2n-2,k-1}^\ub
\end{equation*}
that is the \textbf{second} term in the right side of \eqref{recJubnkpf}.
For $(\pi, (i,j))\in \mathcal{F}^{(5)^+}_{2n-2,k-1}\uplus\mathcal{F}^{(5)^-}_{2n-2,k-1}\uplus\mathcal{F}^{(2)^+}_{2n-2,k-1}\uplus\mathcal{F}^{(2)^-}_{2n-2,k-1}$,
define
\begin{equation*}
\Psi((\pi, (i,j)))=(\sigma,r)
\end{equation*}
by
\begin{equation}\label{eqpfJn2ndr}
(\sigma,r)=\left\{
      \begin{array}{ll}
      (\xi_{(i,j)}(\pi),\chi_{ij}), & \hbox{if $(\pi,(i,j))\in\mathcal{F}^{(5)^+}_{2n-2,k-1}$ ,} \\[6pt]
      (\overline{\xi}_{(i,j)}(\pi),\chi_{ij}), & \hbox{if $(\pi,(i,j))\in\mathcal{F}^{(5)^-}_{2n-2,k-1}$,}\\[6pt]
      (\eta'_{i}(\pi),i+1), & \hbox{if $(\pi,(i,j))_1\in\mathcal{F}^{(2)^+}_{2n-2,k-1}$,}\\[6pt]
      (\eta'_{i}(\pi),i), & \hbox{if $(\pi,(i,j))_2\in\mathcal{F}^{(2)^+}_{2n-2,k-1}$,}\\[6pt]
      (\overline{\eta'}_{i}(\pi),i+1), & \hbox{if $(\pi,(i,j))_1\in\mathcal{F}^{(2)^-}_{2n-2,k-1}$,}\\[6pt]
      (\overline{\eta'}_{i}(\pi),i), & \hbox{if $(\pi,(i,j))_2\in\mathcal{F}^{(2)^-}_{2n-2,k-1}$.}
      \end{array}
    \right.
\end{equation}

Let
\begin{equation*}
 \mathcal{F}^{(3)^+}_{n,k}=\{
  (\pi,(i,j))\mid \pi\in\mathcal{J}^\ub_{n,k} \text{ and $d^+(i,j)=(1,1)$ in $P_\pi$ with $i\neq j$}\},
\end{equation*}
\begin{equation*}
   \mathcal{F}^{(3)^-}_{n,k}=\{
  (\pi,(i,j))\mid \pi\in\mathcal{J}^\ub_{n,k} \text{ and $d^-(i,j)=(1,1)$ in $P_\pi$ with $i\neq j$}\}.
\end{equation*}
and
\begin{equation*}
  \mathcal{F}^{(4)^+}_{n,k}=\{
  (\pi,(i,i))_1,(\pi,(i,i))_2\mid \pi\in\mathcal{J}^\ub_{n,k}
  \text{ and $d^+(i,i)=(1,1)$ in $P_\pi$ }\},
\end{equation*}
\begin{equation*}
  \mathcal{F}^{(4)^-}_{n,k}=\{
  (\pi,(i,i))_1,(\pi,(i,i))_2\mid \pi\in\mathcal{J}^\ub_{n,k}
   \text{ and $d^-(i,i)=(1,1)$ in $P_\pi$ }\}.
\end{equation*}
From Theorems \ref{thmpaths} and \ref{thm-enum-dtype}, we obtain
\[
\left|\mathcal{F}^{(3)^+}_{2n-2,k-2}\uplus\mathcal{F}^{(3)^-}_{2n-2,k-2}\right|
=2((2n-k)^2+n-1)J_{2n-2,k-2}^\ub,
\]
and
\[
\left|\mathcal{F}^{(4)^+}_{2n-2,k-2}\uplus\mathcal{F}^{(4)^-}_{2n-2,k-2}\right|
=2(4n-2k+1)J_{2n-2,k-2}^\ub.
\]
Hence
\begin{equation*}
\left|\mathcal{F}^{(3)^+}_{2n-2,k-2}\uplus\mathcal{F}^{(3)^-}_{2n-2,k-2}
\uplus\mathcal{F}^{(4)^+}_{2n-2,k-2}\uplus\mathcal{F}^{(4)^-}_{2n-2,k-2}\right|
=2((k-2)(k-4n)+4 n^2-3n) J_{2n-2,k-2}^\ub
\end{equation*}
is the \textbf{third} term in the right side of \eqref{recJubnkpf}.
For $(\pi, (i,j))\in \mathcal{F}^{(3)^+}_{2n-2,k-2}\uplus\mathcal{F}^{(3)^-}_{2n-2,k-2}\uplus
\mathcal{F}^{(4)^+}_{2n-2,k-2}\uplus\mathcal{F}^{(4)^-}_{2n-2,k-2}$,
define
\begin{equation*}
\Psi((\pi, (i,j)))=(\sigma,r),
\end{equation*}
by
\begin{equation}\label{eqpfJn3rdr}
(\sigma,r)=\left\{
      \begin{array}{ll}
      (\xi_{(i,j)}(\pi),\chi_{ij}), & \hbox{if $(\pi,(i,j))\in\mathcal{F}^{(3)^+}_{2n-2,k-2}$,} \\[6pt]
      (\overline{\xi}_{(i,j)}(\pi),\chi_{ij}), & \hbox{if $(\pi,(i,j))\in\mathcal{F}^{(3)^-}_{2n-2,k-2}$,}\\[6pt]
      (\eta'_{i}(\pi),i+1), & \hbox{if $(\pi,(i,j))_1\in\mathcal{F}^{(4)^+}_{2n-2,k-2}$,}\\[6pt]
      (\eta'_{i}(\pi),i), & \hbox{if $(\pi,(i,j))_2\in\mathcal{F}^{(4)^+}_{2n-2,k-2}$,}\\[6pt]
      (\overline{\eta'}_{i}(\pi),i+1), & \hbox{if
      $(\pi,(i,j))_1\in\mathcal{F}^{(4)^-}_{2n-2,k-2}$,}\\[6pt]
      (\overline{\eta'}_{i}(\pi),i), & \hbox{if $(\pi,(i,j))_2\in\mathcal{F}^{(4)^-}_{2n-2,k-2}$.}
      \end{array}
    \right.
\end{equation}

Therefore, by setting
\begin{equation*}
{\mathcal{F}}_{n,k}=
\mathcal{F}^{(1)^\pm}_{n-2,k}\uplus\mathcal{F}^{(2)\pm}_{n-2,k-1}
\uplus\mathcal{F}^{(3)^\pm}_{n-2,k-2}\uplus\mathcal{F}^{(4)\pm}_{n-2,k-2}
\uplus\mathcal{F}^{(5)^\pm}_{n-2,k-1},
\end{equation*}
we see $|\mathcal{F}_{n,k}|$ is equal to the right side of \eqref{recJubnkpf}.
Similar to the proof of Theorem \ref{thmJobnk},
by combining \eqref{eqpfJn1str}--\eqref{eqpfJn3rdr},
we establish the desired bijection
\begin{equation*}
  \Psi\colon \mathcal{F}_{2n,k}\leftrightarrow \mathcal{H}_{2n,k},
\end{equation*}
and completes the proof.
\qed

\vspace{2ex}

\noindent{\bf Acknowledgements}:
The second author is supported by the National Natural Science Foundation of China (12001078),
the Natural Science Foundation of Chongqing (CSTB2022NSCQ-MSX0465),
and the Doctoral Research Start-up Funding of CQUPT (A2020-31).
The last author is supported by the National Natural Science Foundation of China (11901074).

\end{document}